\documentclass[a4paper,12pt]{amsart}
\usepackage{amssymb}
\usepackage{amsthm}
\usepackage{amscd}
\usepackage{amsmath}
\usepackage[dvipdfm]{graphicx}
\theoremstyle{plain}
\newtheorem{lemma}{Lemma}[subsection]
\newtheorem{prop}[lemma]{Proposition}
\newtheorem{thm}[lemma]{Theorem}
\newtheorem{cor}[lemma]{Corollary}
\newtheorem{aplemma}{Lemma~A.\hspace{-1.5mm}}
\newtheorem{approp}{Proposition~A.\hspace{-1.5mm}}
\newtheorem{apthm}{Theorem~A.\hspace{-1.5mm}}
\newtheorem{apcor}{Corollary~A.\hspace{-1.5mm}}
\newtheorem{intthm}{Theorem}

\usepackage[all]{xy}

\newtheorem{conj}[lemma]{Conjecture}

\theoremstyle{definition}

\newtheorem{rem}{Remark}
\newtheorem{rema}[lemma]{Remark}

\newtheorem{remb}{Remark}

\newtheorem{defi}[lemma]{Definition}
\newtheorem{exa}[lemma]{Example}
\newtheorem{aprem}{Remark~A.\hspace{-1.5mm}}
\newtheorem{apdefi}{Definition~A.\hspace{-1.5mm}}
\newcommand{\bde}{\begin{defi}}
\newcommand{\ede}{\end{defi}\vspace{1mm}}
\newcommand{\ble}{\begin{lemma}}
\newcommand{\ele}{\end{lemma}}
\newcommand{\bpr}{\begin{prop}}
\newcommand{\epr}{\end{prop}}
\newcommand{\bt}{\begin{thm}}
\newcommand{\et}{\end{thm}}
\newcommand{\bco}{\begin{cor}}
\newcommand{\eco}{\end{cor}}
\newcommand{\bre}{\begin{rem}}
\newcommand{\ere}{\end{rem}}
\newcommand{\brea}{\begin{rema}}
\newcommand{\erea}{\end{rema}\vspace{1mm}}
\newcommand{\breb}{\begin{remb}}
\newcommand{\ereb}{\end{remb}\vspace{1mm}}
\newcommand{\bex}{\begin{exa}}
\newcommand{\eex}{\end{exa}}
\newcommand{\bpf}{\begin{proof}}
\newcommand{\epf}{\end{proof}\vspace{1mm}}

\newcommand{\bade}{\begin{apdefi}}
\newcommand{\eade}{\end{apdefi}}
\newcommand{\bale}{\begin{aplemma}}
\newcommand{\eale}{\end{aplemma}}
\newcommand{\bapr}{\begin{approp}}
\newcommand{\eapr}{\end{approp}}
\newcommand{\bat}{\begin{apthm}}
\newcommand{\eat}{\end{apthm}}
\newcommand{\baco}{\begin{apcor}}
\newcommand{\eaco}{\end{apcor}}
\newcommand{\bare}{\begin{aprem}}
\newcommand{\eare}{\end{aprem}}

\newcommand{\be}{\begin{enumerate}}
\newcommand{\ee}{\end{enumerate}}
\newcommand{\bcd}{\[\begin{CD}}
\newcommand{\ecd}{\end{CD}\]}
\newcommand{\bit}{\begin{itemize}}
\newcommand{\eit}{\end{itemize}}
\newcommand{\bq}{\begin{quote}}
\newcommand{\eq}{\end{quote}}
\newcommand{\ba}{\begin{array}}
\newcommand{\ea}{\end{array}}

\newcommand{\mcD}{\mathcal{D}}
\newcommand{\mcE}{\mathcal{E}}
\newcommand{\mcF}{\mathcal{F}}

\newcommand{\mcH}{\mathcal{H}}
\newcommand{\mcI}{\mathcal{I}}
\newcommand{\mcJ}{\mathcal{J}}
\newcommand{\mcK}{\mathcal{K}}
\newcommand{\mcL}{\mathcal{L}}
\newcommand{\mcM}{\mathcal{M}}
\newcommand{\mcN}{\mathcal{N}}
\newcommand{\mcO}{\mathcal{O}}

\newcommand{\mcS}{\mathcal{S}}
\newcommand{\mcT}{\mathcal{T}}

\newcommand{\mbA}{\mathbb{A}}

\newcommand{\mbC}{\mathbb{C}}

\newcommand{\mbR}{\mathbb{R}}

\newcommand{\mbU}{\mathbb{U}}

\newcommand{\mbX}{\mathbb{X}}
\newcommand{\mbY}{\mathbb{Y}}
\newcommand{\mbZ}{\mathbb{Z}}

\newcommand{\mfC}{\mathfrak{C}}

\newcommand{\mfE}{\mathfrak{E}}

\newcommand{\mfM}{\mathfrak{M}}

\newcommand{\mfS}{\mathfrak{S}}

\newcommand{\mfX}{\mathfrak{X}}
\newcommand{\mfY}{\mathfrak{Y}}

\newcommand{\mfc}{\mathfrak{c}}

\newcommand{\mfh}{\mathfrak{h}}
\newcommand{\mfi}{\mathfrak{i}}

\newcommand{\mfk}{\mathfrak{k}}

\newcommand{\mfn}{\mathfrak{n}}

\newcommand{\mfp}{\mathfrak{p}}

\newcommand{\mft}{\mathfrak{t}}

\newcommand{\migi}{\rightarrow}

\newcommand{\isom}{\stackrel{\sim}{\migi}}

\newcommand{\migiincl}{\hookrightarrow}

\newcommand{\migisurj}{\twoheadrightarrow}

\newcommand{\mr}{\mathrm}
\newcommand{\hidden}[1]{\,}

\begin{document}

\title[Moduli  of  log twisted $\mcN =1$  SUSY curves]{Moduli  of log twisted $\mathcal{N} =1$  SUSY curves}
\author{Yasuhiro Wakabayashi}
\date{}
\markboth{Yasuhiro Wakabayashi}{}
\maketitle
\footnotetext{Y. Wakabayashi: Graduate School of Mathematical Sciences, The University of Tokyo, 3-8-1 Komaba, Meguro, Tokyo,  153-8914, Japan;}
\footnotetext{e-mail: {\tt wkbysh@ms.u-tokyo.ac.jp};}
\begin{abstract}
The goal  of the present paper is to construct a   smooth compactification of the moduli superstack classifying  pointed $\mathcal{N} =1$ SUSY (= $\text{SUSY}_1$)  curves.
This construction is  based on the Abramovich-Jarvis-Chiodo compactification of the moduli stack classifying  spin curves.
First, we give a general framework of a theory of  log superschemes (or more generally, log superstacks).
Then, we introduce  the notion of a pointed (stable)  log twisted $\text{SUSY}_1$ curve; it may be thought of as  a  logarithmic and twisted  generalization of  the classical notion of a pointed $\text{SUSY}_1$ curve, as well as   
  a supersymmetric analogue of the notion of a pointed  (log) twisted curve.
The main result of the present paper asserts that the moduli superstack classifying pointed stable  log twisted $\text{SUSY}_1$ curves may be represented by  a log superstack whose underlying superstack is a superproper and supersmooth Deligne-Mumford  superstack.
Consequently,  this moduli superstack forms 
a smooth  compactification different from the compactification  proposed  by P. Deligne.

\end{abstract}
\tableofcontents 
\section*{Introduction}

\vspace{3mm}
\subsection*{0.1} \label{S01}

 The goal  of the present paper is   {\it to provide a  rigorous construction of a  smooth compactification of the moduli superstack classifying  pointed $\mathcal{N} =1$ SUSY  curves}.
Throughout the present paper, we  abbreviate ``$\mathcal{N} =1$ SUSY" to ``$\text{SUSY}_1$" for simplicity.
Recall that $\text{SUSY}_1$ curves and their analytic  counterparts, called super Riemann surfaces, 
 have been widely considered (intensively in the 1980s) in the physical literature on supersymmetry.
Super Riemann surfaces are 
defined to be 
 complex supermanifolds of superdimension $1|1$ satisfying  an additional superconformal condition (cf. e.g., ~\cite{LR}). 
Super Riemann surfaces 
play a role of
  the correct supersymmetric analogue of  Riemann surfaces, and their moduli superspace plays a role  analogous to the role of moduli space classifying Riemann surfaces  in bosonic string theory (cf. ~\cite{Frie}). 
Indeed, just as the world sheet of a bosonic string carries the structure of a Riemann surface, the world sheet in superstring theory is a super Riemann surface.
Also, perturbative calculations in superstring theory are carried out by integration over this moduli superspace. 
 Besides having such physical applications, the theory of super Riemann surfaces and their moduli  is interesting on its own from the mathematical viewpoint.
In order to achieve a deep  understanding of this theory (from the mathematical viewpoint or another),  
 it will be  worth asking the following question regarding their global structure that should be answered:
\begin{quote}\textit{What is  a natural (smooth) compactification of the moduli superspace classifying  super Riemann surfaces (or  more generally,  pointed $\text{SUSY}_1$ curves)?} \end{quote}

\vspace{3mm}
\subsection*{0.2} \label{S02}

To obtain an  answer of   this question, 
  P. Deligne constructed,  in his letter to Y. Manin  (cf. ~\cite{Del}),   a smooth compactification of the moduli superstack classifying (unpointed) $\text{SUSY}_1$ curves; it may be thought of as
 an  analogue of  the Deligne-Mumford compactification of the moduli stack classifying  proper smooth algebraic curves, and   obtained by adding certain divisors at infinity  parametrizing $\text{SUSY}_1$ curves with nodes.
 The main deference from the bosonic case is that,  in constructing  the compactification, we need to allow 
 two different types of degeneration of $\text{SUSY}_1$ curves, called Neveu-Schwarz and Ramond degenerations.
  We refer to   ~\cite{Witten1}, \S\,6,  for a detailed exposition, including the physical viewpoint,  of  the  moduli superspace classifying  $\text{SUSY}_1$ curves  (with marked points).
In the present  paper, we consider a  smooth compactification different from the compactification constructed by P. Deligne, including the case of pointed  $\text{SUSY}_1$ curves.

\vspace{3mm}
\subsection*{0.3} \label{S03}

Let  us describe the main theorem of the present paper.
Let $S_0$ be a noetherian  affine scheme over $\mbZ [\frac{1}{2}]$, $\lambda$ a positive even integer which is invertible in $S_0$, and $(g,r)$ a pair of nonnegative integers such that $r$ is even and  $2g-2 +r>0$.
Write  $\mfS \mfc \mfh_{/S_0}^{\circledS \mr{log}}$ (cf. (\ref{EE13}))  for  the category of fs log superschemes (cf. Definition \ref{d2b}) over $S_0$. 
Also, write
\begin{equation}
{^{\S_1} \overline{\mfM}}_{g,r, \lambda}^{\circledS \mr{log}}
\end{equation} 
(cf. (\ref{EE12})) for the category fibered in groupoids  over $\mfS \mfc \mfh_{/S_0}^{\circledS \mr{log}}$
classifying families of  {\it stable  log twisted $\text{SUSY}_1$ curves of type $(g, r, \lambda)$} (cf. Definition \ref{D03}) parametrized 
by log superschems  in $\mfS \mfc \mfh_{/S_0}^{\circledS \mr{log}}$.
Then, our  main result is the following theorem.
 (See Definitions  \ref{d089} (ii), \ref{d0299}, \ref{d08gg9}, and \ref{d2}  for the  definitions of various notions  appeared in the statement.)
\begin{intthm} 
  \label{y019}
\leavevmode\\
 \ \ \ 
 ${^{\S_1} \overline{\mfM}}_{g,r, \lambda}^{\circledS \mr{log}}$ may be represented by 
 a log superstack whose underlying superstack is a superproper and supersmooth Deligne-Mumford superstack over $S_0$ of relative superdimension $3g-3+r | 2g-2+ \frac{r}{2}$.
\end{intthm}

\vspace{3mm}
\subsection*{0.4} \label{S04}



Let us make a remark on the  main result just described. 
Denote by ${^{\S_1} \mfM}_{g,r}^\circledS$ the moduli superstack classifying  $r$-pointed  (supersmooth) $\text{SUSY}_1$ curves of genus $g$ (in the classical sense) over log superschemes in $\mfS \mfc \mfh^{\circledS \mr{log}}_{/S_0}$.
Since any pointed  (supersmooth) $\text{SUSY}_1$ curve  is a stable log twisted $\text{SUSY}_1$ curve, 
we have  
a natural inclusion ${^{\S_1} \mfM}_{g,r}^\circledS \migiincl {^{\S_1} \overline{\mfM}}_{g,r, \lambda}^{\circledS \mr{log}}$; it  is an open immersion whose  image is dense  and coincides with the locus in which the log structure of ${^{\S_1} \overline{\mfM}}_{g,r, \lambda}^{\circledS \mr{log}}$ becomes  trivial.
Thus,  the moduli  log superstack ${^{\S_1} \overline{\mfM}}_{g,r, \lambda}^{\circledS \mr{log}}$, being our central character,  forms a smooth compactification of ${^{\S_1} \mfM}_{g,r}^\circledS$ (different from the compactification by P. Deligne).
The feature of our compactifcation is 
that we add  divisors at infinity which parametrize pointed $\text{SUSY}_1$ curves (equipped with a logarithmic structure) admitting  at most 
a single type of degeneration.

Next, recall the discussion of the  non-projectedness (cf. ~\cite{Witten2}, \S\,2), as well as the non-splitness   (cf. Definition \ref{d0d99ff}), of ${^{\S_1} \mfM}_{g,r}^\circledS$ considered  in ~\cite{Witten2}.
By applying (an argument  similar to) the argument in {\it loc.\,cit.}\,to our situation, we will be able to verify 
the non-projectedness (and hence, the non-splitness) of (the underlying superstack of) ${^{\S_1} \overline{\mfM}}_{g,r, \lambda}^{\circledS \mr{log}}$    even when $S_0$ is not necessarily $\mr{Spec} (\mbC)$ (where $\mbC$ denotes the field of complex numbers).
Indeed, 
 since the non-projectedness  of ${^{\S_1} \mfM}_{g,r}^\circledS$ implies the non-projectedness of 
$ {^{\S_1} \overline{\mfM}}_{g,r, \lambda}^{\circledS \mr{log}}$, 
it follows from  Theorems 1.1-1.3 in {\it loc.\,cit.}\,that $ {^{\S_1} \overline{\mfM}}_{g,r, \lambda}^{\circledS \mr{log}}$ is non-projected for many cases of $(g,r)$.
This means that ${^{\S_1} \overline{\mfM}}_{g,r, \lambda}^{\circledS \mr{log}}$ cannot  be reconstructed from purely bosonic moduli stacks
 in an elementally fashion and, in some sense, needs to  be studied independently.

\vspace{3mm}
\subsection*{0.5} \label{S05}

Let us  briefly explain the points of our discussion and  the organization of the present paper.
In \S\,1, we give (and recall) a general framework of a theory of superschemes, and more generally, superstacks.
Then, we  define,  in \S\,2,   a logarithmic  structure on a superscheme, as well as a superstack.
The motivation for  introducing the notion of a log superstack (i.e., a superstack equipped with a logarithmic  structure) is
to consider  a  supersymmetric analogue of  (stable) pointed   twisted curves with a canonical log structure.
(We refer to  ~\cite{AV1}, Definition 4.3.1, and ~\cite{Chi1}, Definition 2.4.1, for the definition of a pointed twisted curve,  and to ~\cite{O1}, Theorem 3.5, for the canonical  logarithmic structure defined on a pointed twisted curve.)
 By means of the various notions defined   in \S\S\,1-2, we present, in \S\,3, the definition of  a (stable) pointed log twisted $\text{SUSY}_1$ curve as, roughly speaking, a certain pointed  
  log superstack  of superdimension $1|1$ equipped  with an additional   superconformal structure (cf. Definition \ref{D02} and Definition \ref{D03}).
Thus, for  a  suitable triple $(g,r, \lambda)$ of nonnegative integers,  one may obtain the category  ${^{\S_1} \overline{\mfM}}_{g,r, \lambda}^{\circledS \mr{log}}$ fibered in groupoids, as we introduced above, 
 classifying stable log twisted $\text{SUSY}_1$ curves of type $(g,r, \lambda)$.
Denote by $({^{\S_1} \overline{\mfM}}_{g,r, \lambda})_t^\mr{log}$ 
the restriction   of  ${^{\S_1} \overline{\mfM}}_{g,r, \lambda}^{\circledS \mr{log}}$ to the full subcategory $\mfS \mfc \mfh_{/S_0}^\mr{log} \subseteq  \mfS \mfc \mfh_{/S_0}^{\circledS \mr{log}}$ consisting of fs log schemes (in the classical  sense) .
As discussed in \S\,4,  the kay point is  that 
to giving a family of pointed log twisted $\text{SUSY}_1$ curves parametrized by an fs  log scheme  (i.e., an object in $({^{\S_1} \overline{\mfM}}_{g,r, \lambda})_t^\mr{log}$)  is  equivalent to giving a family of pointed log twisted curves  equipped with an additional data called a pointed   spin structure (cf. Definition \ref{De2}).
(This observation for  the case of unpointed smooth $\text{SUSY}_1$ curves is classical and well-known.)
This implies  (cf. Proposition \ref{P66}) that 
$({^{\S_1} \overline{\mfM}}_{g,r, \lambda})_t^\mr{log}$  is canonically isomorphic to
the moduli stack ${^{\mr{tw}} \overline{\mfM}}_{g,r, \lambda, \mr{spin}}^\mr{log}$  (cf. (\ref{EE33})) classifying  $\lambda$-stable log  twisted  curves of type $(g,r)$ equipped with a pointed spin structure.
On the other hand, D.  Abramovich, T. J. Jarvis, and A. Chiodo  proved (cf. ~\cite{AJ1}, Theorem 1.5.1 and ~\cite{Chi1}, Corollary 4.11)
that  ${^{\mr{tw}} \overline{\mfM}}_{g,r, \lambda, \mr{spin}}^\mr{log}$ ($\cong ({^{\S_1} \overline{\mfM}}_{g,r, \lambda})_t^\mr{log}$) may be represented by a proper smooth Deligne-Mumford stack, which  forms a  compatification  of the moduli stack classifying pointed  smooth spin curves (in the classical sense).
Thus, by thickening  this Deligne-Mumford  stack  in the fermionic directions in a way that a  universal stable log twisted $\text{SUSY}_1$ curve  exists (uniquely), we construct, in \S\,5,  a log superstack  representing ${^{\S_1} \overline{\mfM}}_{g,r, \lambda}^{\circledS \mr{log}}$ and satisfying the  desired conditions described  in Theorem A.

\vspace{5mm}
\hspace{-4mm}{\bf Acknowledgement} \leavevmode\\
 \ \ \ The author cannot express enough his sincere and  deep gratitude to  all those who give  the opportunity or   impart  the great joy of studying
  mathematics to him.
The author means the present paper for a gratitude letter to them.
The author was partially  supported by the Grant-in-Aid for Scientific Research (KAKENHI No.\,15J02721) and the FMSP program at the Graduate School of Mathematical Sciences of the University of Tokyo.


\vspace{10mm}
\section{Superschemes and superstacks} \vspace{3mm}

The  aim of this section is to  give a brief introduction to the theory of superschemes (or more generally, superstacks).
We first  recall  the notion of a superscheme (cf. Definition of \ref{d3})  and then, discuss basic properties of  superschemes and   morphisms between them.
In particular, we define a  super\'{e}tale morphism (cf. Definition \ref{d2}), which is a supersymmetric analogue of an \'{e}tale morphism in the classical sense.
By means of this sort  of morphism, one obtains a category of superschemes equipped with 
the Grothendieck 
 pretopology
  which will be denoted by  $\mfS \mfc \mfh_{/S_0}^\circledS$ (cf. (\ref{FF01})), and moreover,  obtains the definition of a (Deligne-Mumford) superstack (cf. Definitions \ref{D029} and \ref{d0299}).
Finally, we show (cf. Proposition \ref{p0207}) that
any Deligne-Mumford stack admits a complete versal family  which is isomorphic to a split (cf. Definition \ref{d0d99ff}) and supersmooth (cf. Definition \ref{d0819}) superscheme.
Basic references for the notion of a {\it superscheme} are, e.g.,   ~\cite{CR},  ~\cite{G}, and ~\cite{Man}.

Let $R_0$ be a noetherian ring over $\mbZ [\frac{1}{2}]$.
Throughout the present paper,
{\it all schemes are assumed to be locally noetherian schemes over  the  affine  scheme $S_0 := \mr{Spec}(R_0)$ and all morphisms of schemes are assumed to be locally of finite presentation.}

\vspace{5mm}
\subsection{Superschemes} \label{S11}
\leavevmode\\
\vspace{-4mm}

First, recall the definition of a superscheme as follows.
\vspace{3mm}
\bde \label{d3}\leavevmode\\
 \vspace{-5mm}
\begin{itemize}
\item[(i)]
 A {\bf superscheme} (over $S_0$) is a pair  $X^\circledS := (X_b, \mcO_{X^\circledS})$ consisting of a (locally noetherian) scheme $X_b$ over $S_0$ and a {\it coherent} sheaf of superalgebras $\mcO_{X^\circledS}$  over $\mcO_{X_b}$ 
 such that the natural morphism $\mcO_{X_b} \migi \mcO_{X^\circledS}$ is injective and its image coincides with 
 the bosonic (i.e., even) part of $\mcO_{X^\circledS}$.
 We  write $\mcO_{X_f}$ for the fermionic (i.e., odd) part of $\mcO_{X^\circledS}$ and identify $\mcO_{X_b}$ with the bosonic part via the injection $\mcO_{X_b} \migiincl \mcO_{X^\circledS}$ (hence, $\mcO_{X^\circledS} = \mcO_{X_b} \oplus \mcO_{X_f}$).
 \item[(ii)]
 Let $X^\circledS := (X_b, \mcO_{X^\circledS})$ and $Y^\circledS := (Y_b, \mcO_{Y^\circledS})$ be two superschemes (over $S_0$).
 A {\bf morphism of superschemes (over $S_0$)} from $Y^\circledS$  to $X^\circledS$  is a pair $f^\circledS := (f_b, f^\flat)$ consisting of a  morphism $f_b : Y_b \migi X_b$ of schemes (over $S_0$, which is locally of finite presentation) and a morphism of superalgebras $f^\flat :  f_b^*(\mcO_{X^\circledS}) \ (:= \mcO_{Y_b} \otimes_{f_b^{-1}(\mcO_{X_b})} f_b^{-1}(\mcO_{X^\circledS})) \migi \mcO_{Y^\circledS}$ over $\mcO_{Y_b}$. 
 \end{itemize}
 \ede
\vspace{3mm}

We always identify  any  scheme $X_b$ (over $S_0$) with  a superscheme $X^\circledS := (X_b, \mcO_{X^\circledS})$ with $\mcO_{X_f} =0$.  

\vspace{1mm}
\begin{rema} \label{r486} \leavevmode\\
 \ \ \ 
 The definition of a superscheme may differ from the usual definition in the sense that  the fermionic part of the structure sheaf of a superscheme  is  assumed not to be coherent.
 In fact,  that condition is usually regarded as an additional condition on a superscheme which is, in  ~\cite{G}, Definition 2.6,  referred to as ``{\it fermionically of finite presentation}".
 But, in the present paper, we only deal with superschemes fermionically of finite presentation in order that  super\'{e}tale morphisms, as well as  supersmooth morphisms,  defined later are  well-behaved.
 This is why we define the notion of a superscheme as above.
 \end{rema}
\vspace{3mm}

 Let $X^\circledS$ be 
a superscheme and  $\mcF$   a left $\mcO_{X^\circledS}$-supermodule.
We  write 
\begin{align}
\mcF_b \ \ \  (\text{resp.},  \ \mcF_f)
\end{align}
 for the bosonic (resp., fermionic) part of $\mcF$ (hence $\mcF= \mcF_b \oplus \mcF_f$).
 $\mcF$ may be considered as 
a right $\mcO_{X^\circledS}$-supermodule  equipped with an 
 $\mcO_{X^\circledS}$-action 
 given 
by $m \cdot a := (-1)^{|m| \cdot |a|}  a \cdot m$ for homogeneous local sections $a \in \mcO_{X^\circledS}$ ,  $m \in \mcF$ (where $|-|$ denotes the parity function).
By an {\bf $\mcO_{X^\circledS}$-supermodule}, we shall means simply  a left $\mcO_{X^\circledS}$-supermodule,  which   is often treated as   a right $\mcO_{X^\circledS}$-supermodule by this consideration.
Also,  
by a {\bf  supervector bundle (of superrank $m|n$)} on $X^\circledS$ (where both $m$ and  $n$ are nonnegative integers), we mean a locally free (left)  $\mcO_{X^\circledS}$-supermodule (of superrank $m|n$).


Let   $f^\circledS := (f_b, f^\flat): Y^\circledS \migi X^\circledS$ be  a morphism of superschemes.
If we are given an  $\mcO_{X^\circledS}$-supermodule (resp., an $\mcO_{Y^\circledS}$-supermodule) $\mcF$,
then one may define, via the natural morphism $f_b^{-1}(\mcO_{X^\circledS}) \migi \mcO_{Y^\circledS}$, the pull-back (resp.,   direct image) of $\mcF$ to be the  $\mcO_{Y^\circledS}$-supermodule (resp., the  $\mcO_{X^\circledS}$-supermodule)
\begin{align}
f^{\circledS *}(\mcF) := \mcO_{Y^\circledS} \otimes_{f^{-1}_b (\mcO_{X^\circledS})} f^{-1}_b(\mcF) \ \ 
(\text{resp.,} \ f_*^\circledS (\mcF) := f_{b*}(\mcF)).
\end{align}

\vspace{3mm}
\bde \label{d32}\leavevmode\\
 \ \ \ Let $S^\circledS$ be a superscheme.
  \begin{itemize}
  \item[(i)]
  Let 
  $X^\circledS$ and $Y^\circledS$ are superschemes over  $S^\circledS$ and  
$f^\circledS \ (:= (f_b, f^\flat)) : Y^\circledS \migi X^\circledS$   a morphism of superschemes over $S^\circledS$.
We shall say that  $f^\circledS$  is a  {\bf closed immersion (over $S^\circledS$)} if  $f_b : Y_b \migi X_b$ is a closed immersion
  and $f^\flat : f_b^*(\mcO_{X^\circledS}) \migi \mcO_Y$ is surjective.
\item[(ii)]
Let $X^\circledS$ be a superscheme over $S^\circledS$.
A {\bf closed subsuperscheme} of $X^\circledS$  is an equivalence class of closed immersions  into $X^\circledS$, where two morphisms $f_1^\circledS: Y^\circledS_1 \migi X^\circledS$, $f^\circledS_2 : Y^\circledS_2 \migi X^\circledS$ over $S^\circledS$  are {\bf equivalent} if there exists an isomorphism $\iota^\circledS : Y^\circledS_1 \isom Y^\circledS_2$  satisfying that $f_2^\circledS \circ \iota^\circledS = f^\circledS_1$.
If $f^\circledS : Y^\circledS \migi X^\circledS$ is a closed immersion,
then we shall write $[f^\circledS]$ for the closed subsuperscheme of $X^\circledS$  represented by $f^\circledS$. 
 \end{itemize}
 \ede
\vspace{3mm}

Let  $X^\circledS := (X_b, \mcO_{X^\circledS})$ be a superscheme.
By means of  the morphism 
\begin{equation} \label{e22}
\beta^\circledS_X : X^\circledS \migi X_b
\end{equation}
 corresponding to  
the inclusion $\mcO_{X_b} \migiincl  \mcO_{X^\circledS}$,
$X^\circledS$ may be thought of as a superscheme over the  scheme $X_b$.
The construction of $\beta_X^\circledS$ is evidently  functorial in $X^\circledS$, that is to say,
$\beta^\circledS_X \circ f^\circledS = f_b \circ \beta_Y^\circledS$ for any superscheme $Y^\circledS$ and any morphism
 $f^\circledS \ (:= (f_b, f^\flat))  : Y^\circledS  \migi X^\circledS$.

Denote by 
\begin{align}
\mcN_{X^\circledS}
\end{align}
the superideal of $\mcO_{X^\circledS}$ 
generated by $\mcO_{X_f}$.
We shall write 
\begin{align} \label{e480}
\tau^\circledS_X : X_t \migi X^\circledS
\end{align}
 for the closed
 immersion corresponding to the quotient $\mcO_{X^\circledS} \migisurj \mcO_{X^\circledS}/\mcN_{X^\circledS}$.
Hence, $X_t$ forms  a scheme,  and 
the composite
\begin{align} \label{B01}
\gamma_X := \beta_X^\circledS \circ \tau_X^\circledS : X_t \migi X_b
\end{align}
forms  a closed immersion of schemes  corresponding to the quotient  $\mcO_{X_b} \migisurj \mcO_{X_b}/\mcO_{X_f}^2$ ($= \mcO_{X^\circledS}/\mcN_{X^\circledS}$) by the nilpotent ideal $\mcO_{X_f}^2 \subseteq  \mcO_{X_b}$.
Any morphism  $f^\circledS : Y^\circledS \migi X^\circledS$  induces a morphism $f_t : Y_t \migi X_t$ of schemes satisfying that $f^\circledS \circ \tau^\circledS_Y = \tau_X^\circledS \circ f_t$ and 
$f_b \circ \tau_Y = \tau_X \circ f_t$.
In particular, any morphism $Z \migi X^\circledS$ (where $Z$ is a scheme) decomposes  as  $Z \migi X_t \stackrel{\tau_X^\circledS}{\migi} X^\circledS$ for a unique morphism $Z \migi X_t$ of schemes.

Finally, for each nonnegative  integer $n$, we write 
\begin{align} \label{e467}
\mr{gr}_{X^\circledS}^n  := \mcN^n_{X^\circledS} /\mcN_{X^\circledS}^{n+1},
\end{align}
 which may be thought of as an $\mcO_{X_t}$-module.

\vspace{5mm}
\subsection{Morphisms of  superschemes} \label{S12}
\leavevmode\\
\vspace{-4mm}

We shall consider analogues of flat morphism and \'{e}tale morphisms to superschemes.
 Let $f^\circledS  \  (:= (f_b, f^\flat)) : Y^\circledS \migi X^\circledS$ be a morphism of superschemes.

\bde \label{d302d}\leavevmode\\
 \ \ \ 
 We shall say that  $f^\circledS$ is {\bf bosonic} if for any scheme $Z$ together with a morphism $Z \migi X^\circledS$, the fiber product
 $Y^\circledS \times_{f^\circledS, X^\circledS} Z$ is a scheme.
 (Here, we note that the superschemes and morphisms between them form a category, in which the fiber products 
  exist. See  ~\cite{CR}, Corollary 10.3.9.)
  \ede

\bde \label{d312d}\leavevmode\\
 \ \ \ 
 We shall say that  $f^\circledS$ is {\bf superflat} if for any point $y$ of $Y_b$  the homomorphism $\mcO_{X^\circledS, f_b (y)} \migi \mcO_{Y^\circledS, y}$ of local rings  induced by $f^\circledS$ is flat.
 \ede

\begin{rema} \label{r48} \leavevmode\\
 \ \ \ 
Suppose that  $f^\circledS$ is  bosonic and superflat.
According to ~\cite{G}, Lemma 2.7 and Proposition 2.1 (cf. Remark \ref{r486}),  the following properties hold 
(where  although the results of {\it loc.\,cit.} is assumed that $S_0 = \mr{Spec}(\mbC)$,  one may prove the same assertion for  our general case):
\begin{itemize}
\item[(i)]
The homomorphism $f^\flat$ induces, by restriction,   isomorphisms 
\begin{align} \label{E1}
f_b^*(\mcO_{X_f}) \isom  \mcO_{Y_f} \ \ \text{and} \ \   f_b^*(\mcO^2_{X_f}) \isom  \mcO^2_{Y_f}
\end{align}
(hence,  we have $ f_b^*(\mcN_{X^\circledS}) \isom  \mcN_{Y^\circledS}$). 
In particular, the natural morphisms 
\begin{align} \label{E3}
Y^\circledS \migi Y_b \times_{f_b, X_b, \beta^\circledS_X} X^\circledS \ \ \text{and} \  \ Y_t \migi Y^\circledS \times_{f^\circledS, X^\circledS, \tau^\circledS_X} X_t
\end{align}
 are isomorphisms.
\item[(ii)]
The underlying morphism  $f_b : Y_b \migi X_b$ is flat (in the classical sense).
\end{itemize}
 \end{rema}
\vspace{3mm}
\bde \label{d2}\leavevmode\\
 \ \ \ 
We shall say that  $f^\circledS$  is {\bf  super\'{e}tale} if 
$f^\circledS$ is bosonic and superflat,  and the flat morphism $f_b : Y_b \migi X_b$ (cf. Remark \ref{r48} (ii) above) is  unramified.
 \ede
\vspace{3mm}

\bpr \label{p0607} \leavevmode\\
 \ \ \
For a superscheme $Z^\circledS$ over $S_0$, we shall denote by  $\mfE \mft_{/Z^\circledS}$   the category defined as follows:
\begin{itemize}
\item[$\bullet$]
The {\it objects}  are  super\'{e}tale  morphisms $W^\circledS \migi Z^\circledS$ of superschemes  to $Z^\circledS$;
\item[$\bullet$]
The {\it morphisms}  from $W_1^\circledS \migi Z^\circledS$ to $W_2^\circledS \migi Z^\circledS$ (where both $W_1^\circledS \migi Z^\circledS$ and $W_2^\circledS \migi Z^\circledS$ are objects of this category) are morphisms $W_1^\circledS \migi W_2^\circledS$ of superschemes over $Z^\circledS$.
\end{itemize}
 Then, 
 the functor
   \begin{align} \label{FF02}
   \mfE \mft_{/X^\circledS} \isom \mfE \mft_{/X_t}.
   \end{align}
determined by
   base-change  $Y^\circledS  \mapsto Y^\circledS \times_{X^\circledS, \tau_X^\circledS} X_t$ is  an equivalence of categories. 
   In particular, if ${X'}^\circledS$ and ${X''}^\circledS$ are  superschemes  over $S_0$ such that
  $(X'_t)_\mr{red} \cong (X''_t)_\mr{red}$ (where $(-)_{\mr{red}}$ denotes the  reduced scheme associated with the scheme $(-)$), then we have $\mfE \mft_{/{X'}^\circledS} \cong \mfE \mft_{/{X''}^\circledS}$.
 \epr
\begin{proof}
We shall construct a functor $\mfE \mft_{/X_t} \migi \mfE \mft_{/X^\circledS}$.
Let $Y_0 \migi X_t$ be an object in $\mfE \mft_{/ X_t}$ (i.e.,  $Y_0 \migi X_t$ is \'{e}tale in the classical sense).
Since $X_b$ is a nilpotent thickening of $X_t$ (via the  closed immersion $\gamma_X$),
$Y_0$ extends
 {\it uniquely}  to an \'{e}tale scheme $Y_1$ over $X_b$.
The  superscheme  $Y_1 \times_{X_b} X^\circledS$ (together with the projection to $X^\circledS$) is an object of $\mfE \mft_{/X^\circledS}$ whose image of the functor (\ref{FF02}) is isomorphic to $Y_0$.
The assignment $Y_0 \mapsto Y_1 \times_{X_b} X^\circledS$ is well-defined and  functorial with respect to $Y_0$, and hence, determines a functor $\mfE \mft_{/X_t} \migi  \mfE \mft_{/X^\circledS}$.
This functor  is verified to be   the inverse to 
the functor (\ref{FF02}).
This completes the proof  of Proposition \ref{p0607}.
\end{proof}

\vspace{5mm}
\subsection{The category of superschemes} \label{S13}
\leavevmode\\
\vspace{-4mm}


Write
\begin{align} \label{FF01}
\mfS \mfc \mfh_{/S_0} \ \  (\text{resp.,} \ \mfS \mfc \mfh_{/S_0}^\circledS) 
\end{align}
for the category whose {\it objects} are 
  schemes (resp.,  superschemes)  and  whose {\it morphisms} are morphisms of schemes  (resp., morphisms of  superschemes).
By the natural inclusion $\mfS \mfc \mfh_{/S_0} \migiincl \mfS \mfc \mfh_{/S_0}^\circledS$, we identify $\mfS \mfc \mfh_{/S_0}$ with a full subcategory of  $\mfS \mfc \mfh_{/S_0}^\circledS$.
The fiber products and finite coproducts exist in $\mfS \mfc \mfh_{/S_0}^\circledS$,
and the  inclusion $\mfS \mfc \mfh_{/S_0} \migiincl \mfS \mfc \mfh_{/S_0}^\circledS$ preserves the  fiber products and finite coproducts.
When there is fear of confusion, by a {\bf stack (over $S_0$)}, we  mean a stack over the site $\mfS \mfc \mfh^{}_{/S_0}$ with respect to the \'{e}tale pretopology.
We shall equip   $\mfS \mfc \mfh^{\circledS}_{/S_0}$ with the Grothendieck pretopology 
consisting of coverings $\{ U^\circledS_i \migi X^\circledS \}_{i \in I}$, where each $U^\circledS_i \migi X^\circledS$  is a  super\'{e}tale morphism such that (the underlying morphism between schemes of) $\coprod_{i \in I} U^\circledS_i \migi X^\circledS$ is surjective;
we 
shall refer to this pretopology  as the {\bf  super\'{e}tale pretopology}.
One verifies that the property on a morphism in $\mfS \mfc \mfh_{/S_0}^\circledS$ of being  bosonic  (resp., superflat; resp.,  super\'{e}tale)  is closed under composition and base-change, and satisfies descent for super\'{e}tale coverings.

\vspace{3mm}
\subsection{Split superschemes and affine superschemes} \label{S14}
\leavevmode\\
\vspace{-4mm}

Let $\underline{X}$ be 
a scheme and
  $\mcE$ a coherent  $\mcO_{\underline{X}}$-module.
Consider the exterior algebra $\bigwedge^\bullet_{\mcO_{\underline{X}}} \mcE$ ($:= \bigoplus_{i \geq 0} \bigwedge^i_{\mcO_{\underline{X}}} \mcE$) associated with $\mcE$ over $\mcO_{\underline{X}}$.
The $\mcO_{\underline{X}}$-subalgebra $\bigwedge^{\mr{even}}_{\mcO_{\underline{X}}} \mcE$ ($:= \bigoplus_{i : \mr{even}} \bigwedge^i_{\mcO_{\underline{X}}} \mcE$) defines (since it is commutative) a relative affine space $\mcS pec (\bigwedge^{\mr{even}}_{\mcO_{\underline{X}}} \mcE)$ over $\underline{X}$.
Also,   $\bigwedge^\bullet_{\mcO_{\underline{X}}} \mcE$  may be thought of as a coherent $\mcO_{\mcS pec (\bigwedge^{\mr{even}}_{\mcO_{\underline{X}}} \mcE)}$-module.
Thus, we obtain a  superscheme
\begin{align} \label{E12}
\langle \underline{X}, \mcE \rangle^\circledS := (\mcS pec ({\bigwedge}_{\mcO_{\underline{X}}}^{\mr{even}} \mcE), {\bigwedge}_{\mcO_{\underline{X}}}^\bullet \mcE).
\end{align}
The inclusion $\mcO_{\underline{X}} \ (= \bigwedge^0_{\mcO_{\underline{X}}} \mcE) \migiincl \bigwedge^\bullet_{\mcO_{\underline{X}}} \mcE$ defines  a morphism
\begin{align} \label{E14}
\langle \beta \rangle^\circledS_{\underline{X}, \mcE} : \langle \underline{X}, \mcE \rangle^\circledS \migi \underline{X}
\end{align}
of superschemes.

\vspace{3mm}
\bde \label{d0d99ff}
\leavevmode\\
\ \ \
We shall say that a superscheme $Z^\circledS$  is {\bf split} if $Z^\circledS \cong \langle \underline{X}, \mcE \rangle^\circledS$ for some scheme $\underline{X}$  and a coherent $\mcO_{\underline{X}}$-module $\mcE$.
 \ede
\vspace{3mm}

Next, we shall recall the notion of an affine superscheme.

\vspace{3mm}
\bde \label{d0d99}
\leavevmode\\
\ \ \
If $R := R_b \oplus R_f$ is a superring, then we shall write
\begin{align}
\mr{Spec}(R)^\circledS
\end{align}
 for the superspectrum of $R$ (i.e., ``$\underline{\mr{Spec}} R$" in the sense of ~\cite{CR}, Definition 10.1.1).
We shall say that a superscheme $X^\circledS$  is {\bf affine} if it is isomorphic to $\mr{Spec}(R)^\circledS$ for some superring $R$ (in particular,  both $X_b$  and $X_t$ are affine schemes).
 \ede
\vspace{3mm}

If both $m$ and  $n$ are nonnegative integers and $S^\circledS$ is a superscheme, then we shall write
\begin{align}
\mbA_{S^\circledS}^{m |n} & := S^\circledS \times _{S_0} \mr{Spec}(R_0 [t_n, \cdots, t_m] \otimes_{R_0}  {\bigwedge}_{R_0} (R_0  \psi_1 \oplus \cdots \oplus R_0 \psi_n))^\circledS \\
\big(& \  = S^\circledS \times_{S_0} \langle \mbA_{S_0}^{m|0}, \bigoplus_{l=1}^n \mcO_{\mbA_{S_0}^{m|0}} \psi_i\rangle^\circledS \big), \notag
\end{align}
where the $t_1, \cdots, t_m$ are ordinary  indeterminates and $\psi_1, \cdots, \psi_n$ are fermionic (i.e., anticommuting) indeterminates.

The following assertion is immediately verified from the definition of $\mbA^{m|n}_{(-)}$.
\vspace{3mm}
\bpr \label{P0} \leavevmode\\
 \ \ \
Let $f^\circledS : Y^\circledS \migi X^\circledS$ be a morphism of superschemes.
Then, the functorial (with respect to $Y^\circledS$) map of sets
\begin{align} \label{FF04}
\{ \widetilde{f}^\circledS \in \mr{Hom}_{\mfS \mfc \mfk_{/S_0}^\circledS} (Y^\circledS, \mbA_{X^\circledS}^{m|n}) \ | \ \mr{pr} \circ \widetilde{f}^\circledS = f^\circledS \} & \migi \Gamma (Y_b, \mcO_{Y_b})^{\oplus m} \oplus \Gamma (Y_b, \mcO_{Y_f})^{\oplus n} \\
\widetilde{f}^\circledS \hspace{10mm} & \mapsto ((\widetilde{f}^\flat (t_l))_{l=1}^m, (\widetilde{f}^\flat (\psi_{l'}))_{l'=1}^n) \notag
\end{align}
is bijective, where $\mr{pr}$ denotes the natural projection $\mbA_{X^\circledS}^{m|n} \migi X^\circledS$.
 \epr

\vspace{5mm}
\subsection{Supersmooth morphisms} \label{S15}
\leavevmode\\
\vspace{-4mm}

Let   $m$ and $n$ be  nonnegative integers.

\bde \label{d0819}
\leavevmode\\
\ \ \
Let $f^\circledS : Y^\circledS \migi X^\circledS$ be a morphism of superschemes.
We shall say that $f^\circledS$ is {\bf supersmooth of relative superdimension $m|n$}
if  there exists,  super\'{e}tale locally on $Y^\circledS$,  a super\'{e}tale morphism $Y^\circledS \migi X^\circledS \times_{S_0} \mbA_{S_0}^{m|n}$
over $X^\circledS$. 
 \ede

\begin{rema} \label{r488} \leavevmode\\
 \ \ \ 
Let $n$ be a nonnegative integer.
A morphism $f^\circledS :Y^\circledS \migi X^\circledS$ of superschemes
  is supersmooth of relative superdimension $n|0$ 
if and only if $f^\circledS$ is superflat and $f_b : Y_b \migi X_b$ is, in the classical sense,  smooth  of relative dimension $n$  (i.e., all nonempty fibers are  equidimensional of dimension $n$).
In particular, $f^\circledS$ is   supersmooth of relative superdimension $0|0$ if and only if it is super\'{e}tale.
 \end{rema}

\bpr \label{p0201} \leavevmode\\
 \ \ \
 Let $X^\circledS$ be
 a superscheme. 
 Then, the following two conditions (a) and  (b) are equivalent:
\vspace{1mm}
 \begin{itemize}
 \item[(a)]
 $X^\circledS$ is supersmooth over $S_0$ of relative superdimension $m|n$;
\vspace{1mm}
\item[(b)]
$X_t$ is smooth over $S_0$ of relative dimension $m$, the $\mcO_{X_t}$-module $\mr{gr}_{X^\circledS}^1$ is locally free of  rank $n$,  and there exists, super\'{e}tale  locally on $X^\circledS$, 
an isomorphism $X^\circledS \isom \langle X_t, \mr{gr}_{X^\circledS}^1 \rangle$ 
which makes the diagram
\begin{align}
\xymatrix{
 & X_t \ar[ld]_{\tau^\circledS_X} \ar[rd]^{\tau^\circledS_{ \langle X_t, \mr{gr}_{X^\circledS}^1 \rangle}}& 
 \\
 X^\circledS \ar[rr]^{\sim} & & \langle X_t, \mr{gr}_{X^\circledS}^1 \rangle
}
\end{align}
commute.
  \end{itemize} 
  In particular, $X^\circledS$ is split and supersmooth over $S_0$ of relative superdimension $m|n$ if and only if $X^\circledS \cong \langle \underline{X}, \mcE\rangle^\circledS$ for some smooth scheme $\underline{X}$ over $S_0$ of relative dimension $m$ and some vector bundle  $\mcE$ on $\underline{X}$ of rank $n$.
   \epr
\begin{proof}
Since the latter assertion follows directly  from the former assertion, it suffices to prove  only  the former assertion, i.e., the equivalence (a) $\Leftrightarrow$ (b).
The implication (b) $\Rightarrow$ (a) is clear.
We shall prove  (a) $\Rightarrow$ (b).
After possibly replacing $X^\circledS$ with its super\'{e}tale covering,  we may assume, without loss of generality, that $X_b$ is affine and  there exists a (globally defined) super\'{e}tale morphism $\pi^\circledS \ (:= (\pi_b, \pi^\flat)) : X^\circledS \migi \mbA_{S_0}^{m|n}$ over $S_0$.
Then, $\pi^\flat$ restricts to isomorphisms
 \begin{align} \label{E5}
 \pi_b^*(\mcO_{(\mbA_{S_0}^{m|n})_f}^i) \cong \mcO_{X_f}^i \hspace{3mm} (i =1,2)
 \end{align}
  (cf. Remark \ref{r48}).
  In particular, 
 the commutative square diagram 
\begin{align}
\xymatrix{
X_t  \ar[r]^{\gamma_X} \ar[d]_{ \pi_t} & X_b \ar[d]^{\pi_b}\\
  (\mbA_{S_0}^{m|n})_t \ar[r]_{\gamma_{\mbA_{S_0}^{m|n}}} & (\mbA_{S_0}^{m|n})_b 
}
\end{align}
is cartesian.
It follows that  $X_t$ is \'{e}tale over   $(\mbA_{S_0}^{m|n})_t$ ($=\mbA_{S_0}^{m|0}$), and hence,  smooth over $S_0$ of relative dimension $m$.
Since $X_b$ is affine, there exists a morphism $\iota_X  : X_b \migi X_t$ over $S_0$ satisfying that $\iota_X \circ \gamma_X = \mr{id}_{X_t}$.
The isomorphisms (\ref{E5}) yield an isomorphism
\begin{align} \label{e097}
\mr{gr}^1_{\pi^\flat} :  (\pi_t^* (\mr{gr}^1_{\mbA_{S_0}^{m|n}}) =) \ \pi_b^* (\mr{gr}^1_{\mbA_{S_0}^{m|n}}) \isom \mr{gr}^1_{X^\circledS}.
\end{align}
In particular, we have  $\mr{gr}^1_{X^\circledS} \cong \mcO^{\oplus n}_{X_t}$.
By Proposition \ref{P0}, one may find
  a morphism 
  \begin{align}
  \widetilde{\iota}_X^\circledS : X^\circledS \migi \langle X_t, \mr{gr}_{X^\circledS}^1 \rangle^\circledS
  \end{align} 
  whose composite with the projection $\langle \beta \rangle_{X_t, \mr{gr}_{X^\circledS}^1}^\circledS  : \langle X_t, \mr{gr}_{X^\circledS}^1 \rangle^\circledS \migi X_t$ coincides with $\iota_X \circ \beta_X^\circledS : X^\circledS \migi X_t$ and which makes
 the diagram 
\begin{align}
\xymatrix{X \ar[rr]^{\widetilde{\iota}_X^\circledS} \ar[rd]_{\pi^\circledS}& &  \langle X_t, \mr{gr}_{X^\circledS}^1 \rangle^\circledS  \ar[dl]^{{\pi'}^\circledS} \\
& \mbA^{m|n}_{S_0} &
}
\end{align}
commute, where 
${\pi'}^\circledS$ denotes  the  morphism determined (by means of the canoncial isomorphism  $\mbA^{m|n}_{S_0} \isom \langle \mbA^{m|0}_{S_0}, \pi_t^* (\mr{gr}^1_{\mbA_{S_0}^{m|n}}) \rangle^\circledS$) by  both $\pi_t : X_t \migi \mbA^{m|0}_{S_0}$ and the morphism 
$\mr{gr}^1_{\pi^\flat}$.
Observe that  both $\pi^\circledS$ and ${\pi'}^\circledS$ are superflat  (since $\langle X_t, \mr{gr}_{X^\circledS}^1 \rangle^\circledS \cong \mbA_{S_0}^{m|n} \times_{\mbA_{S_0}^{m|0}, \pi_t} X_t$ and $\pi_t$ is \'{e}tale).
On the other hand, 
$\widetilde{\iota}_X$ restricts, via base-change by $\tau^\circledS_{\mbA_{S_0}^{m|n}} : (\mbA_{S_0}^{m|n})_t \migi \mbA_{S_0}^{m|n}$, to 
 the identity morphism of $X_t$.
This implies that $\widetilde{\iota}_X$ is an isomorphism, and hence,  completes the proof of Proposition \ref{p0201}.
\end{proof}

\vspace{1mm}

\subsection{Superstacks} \label{S16}
\leavevmode\\
\vspace{-5mm}

\bde\label{D029}\leavevmode\\
\vspace{-5mm}
\begin{itemize}
\item[(i)]
A  {\bf   superstack  (over $S_0$)}  is a category fibered in groupoids 
$Z^\circledS \migi \mfS \mfc \mfh^{\circledS}_{/S_0}$
 over  $\mfS \mfc \mfh^{\circledS}_{/S_0}$  which is a stack with respect to the  super\'{e}tale pretopology.
\item[(ii)]
Let  $Z^\circledS_1 \migi \mfS \mfc \mfh^{\circledS}_{/S_0}$ and $Z^\circledS_2 \migi \mfS \mfc \mfh^{\circledS}_{/S_0}$ be  superstacks.
A {\bf morphism of  superstacks} from $Z^\circledS_1$ to $Z^\circledS_2$ is a functor $Z^\circledS_1 \migi Z^\circledS_2$ over $\mfS \mfc \mfh^{\circledS}_{/S_0}$.
\end{itemize}
\ede

  One verifies  immediately that the superstacks and the morphisms of superstacks form a $2$-category, in which the $2$-fiber products and finite coproducts exist.
The natural inclusion from $\mfS \mfc \mfh_{/S_0}^\circledS$ into this category preserves such limits.


\begin{rema} \label{r4a8} \leavevmode\\
 \ \ \ 
 For a superscheme  $X^\circledS$,  the set-valued contravariant  functor $\mr{Hom}_{\mfS \mfc \mfh^\circledS_{/S_0}}(-, X^\circledS)$ on $\mfS \mfc \mfh^\circledS_{/S_0}$ is verified (from a standard argument in descent theory) to be a sheaf on $\mfS \mfc \mfh^\circledS_{/S_0}$ with respect to super\'{e}tale pretopology (cf. ~\cite{G}, Lemma 2.8).
 We always identify any superscheme  $X^\circledS$  with the superstack corresponding to  the sheaf $\mr{Hom}_{\mfS \mfc \mfh^\circledS_{/S_0}}(-, X^\circledS)$.
  \end{rema}

\vspace{1mm}

\begin{rema} \label{r4b8} \leavevmode\\
 \ \ \ 
 If $\underline{Z}$ is a  stack over $S_0$, then,  in a natural manner,  one may consider it as a superstack.
More precisely, let us define the category $\underline{Z}^{\circledS \text{-}\mr{triv}} \migi \mfS \mfc \mfh_{/S_0}^\circledS$ fibered in groupoids over $\mfS \mfc \mfh^\circledS_{/S_0}$  as follows:
\begin{itemize}
\item[$\bullet$]
The {\it  objects} in $\underline{Z}^{\circledS \text{-}\mr{triv}}$ are  pairs $(X^\circledS, x)$ consisting of a superscheme  $X^\circledS$ and a morphisms $x :  X_b \migi \underline{Z}$ of stacks;
\item[$\bullet$]
The {\it  morphisms}  from an object $(Y^\circledS, y)$ to an object $(X^\circledS, x)$ are morphisms $f^\circledS : Y^\circledS \migi X^\circledS$ of superschemes  satisfying that $x \circ f_b \cong y$;
\item[$\bullet$]
The {\it functor} $\underline{Z}^{\circledS \text{-}\mr{triv}} \migi \mfS \mfc \mfh^\circledS_{/S_0}$ is given by assigning $(X^\circledS, x) \mapsto X^\circledS$ (for any object $(X^\circledS, x)$ in $\underline{Z}^{\circledS \text{-}\mr{triv}}$) and $f^\circledS \mapsto f^\circledS$ (for any morphism $f^\circledS$ in $\underline{Z}^{\circledS \text{-}\mr{triv}}$).
\end{itemize}
Then, $\underline{Z}^{\circledS \text{-}\mr{triv}}$ forms a  superstack.
The assignment $\underline{Z} \mapsto \underline{Z}^{\circledS \text{-}\mr{triv}}$
determines a fully faithful functor from the category of stacks over $S_0$ to the category of superstacks over $S_0$.
In this manner, we always consider any stack as a  superstack.
  \end{rema}

\vspace{1mm}
\begin{rema} \label{r4c8} \leavevmode\\
  \ \ \ 
   Let $Z^\circledS \migi \mfS \mfc \mfh^\circledS_{/S_0}$ be a  superstack.
The restriction of  this superstack  to 
the subcategory $\mfS \mfc \mfh_{/S_0} \subseteq \mfS \mfc \mfh_{/S_0}^\circledS$ forms
a stack 
\begin{align} \label{E6}
Z_t \migi \mfS \mfc \mfh_{/S_0}
\end{align}
 over $S_0$.
 If, moreover,  $Z^\circledS$ may  represented by a superscheme $X^\circledS$ (i.e., $Z^\circledS \cong  \mr{Hom}_{\mfS \mfc \mfh^\circledS_{/S_0}}(-, X^\circledS)$), then $Z_t$ (in the sense of (\ref{E6}))  may be represented by $X_t$  (in the sense of (\ref{e480})).
 Moreover, if $W^\circledS \migi Z^\circledS$ is a morphism of  superstacks, then it induces a morphism $W_t \migi Z_t$ of stacks.
  \end{rema}


\vspace{3mm}
\bde \label{d089}
\leavevmode\\
\ \ \
Let  $f^\circledS : Y^\circledS \migi X^\circledS$ be a morphism
of   superstacks.
\begin{itemize}
\item[(i)]
  We shall say that $f^\circledS$ is {\bf representable} if, for any morphism ${X'}^\circledS \migi X^\circledS$ of  superstacks  (where ${X'}^\circledS$ is  a superscheme), the fiber product $Y^\circledS \times_{f^\circledS, X^\circledS} {X'}^\circledS$ is a  superscheme.
\item[(ii)]
We shall say that $f^\circledS$ is {\bf superproper} if
the underlying morphism $f_t : Y_t \migi X_t$ of stacks is proper in the classical sense.
 \end{itemize}
 \ede

 \vspace{5mm}

\subsection{Deligne-Mumford superstacks} \label{S17}
\leavevmode\\
\vspace{-5mm}

\bde \label{d0299}
\leavevmode\\
\ \ \
We shall say that a superstack $Z^\circledS$  is  {\bf Deligne-Mumford}
 if it satisfies the following two conditions:
 \begin{itemize}
 \item[(i)]
 The  diagonal morphism  $Z^{\circledS} \migi Z^{\circledS}  \times_{S_0} Z^{\circledS}$ is  representable and the associated (representable) morphism $Z_t \migi Z_t \times_{S_0} Z_t$ of stacks is separated and quasi-compact in the classical sense;
 \item[(ii)]
 There exists a superscheme $U^\circledS$ over $S_0$ together with a representable
 morphism  $U^\circledS \migi Z^\circledS$ of superstacks over $S_0$ such that
 for each superscheme $V^\circledS$  over $Z^\circledS$, the morphism
 $U^\circledS \times_{Z^\circledS} V^\circledS \migi V^\circledS$  of superschemes (where $U^\circledS \times_{Z^\circledS} V^\circledS$ is necessarily a superscheme thank to  condition (i))  is surjective and super\'{e}tale.
 \end{itemize}
 We shall refer to such
 a superscheme $U^\circledS$ (together with $U^\circledS \migi Z^\circledS$)
   as a {\bf complete versal family for $Z^\circledS$}.
 \ede

\vspace{3mm}
\bde \label{d08gg9}
\leavevmode\\
\ \ \
Let  $f^\circledS : Y^\circledS \migi X^\circledS$ be a morphism
of Deligne-Mumford   superstacks over $S_0$,  and let $m$, $n$ be nonnegative integers.
  We shall say that $f^\circledS$ is {\bf super\'{e}tale} (resp., {\bf supersmooth of relative superdimension $m|n$}) if for any $2$-commutative diagram
  \begin{align}
  \xymatrix{
  V^\circledS \ar[r] \ar[rd]_{h^\circledS} & Y'^{\circledS} \ar[r]  \ar[d]  \ar@{}[rd]|{\Box} & Y^\circledS \ar[d]^{f^\circledS} 
  \\
  & U^\circledS \ar[r] & X^\circledS,
  }
  \end{align}
  where $U^\circledS$ and $V^\circledS$  are complete versal families for $X^\circledS$ and $Y'^{\circledS} := Y^\circledS \times_{X^\circledS} U^\circledS$ respectively, the  morphism $h^\circledS : V^\circledS \migi U^\circledS$ of superschemes is super\'{e}tale in the sense of Definition \ref{d2} (resp., supersmooth of relative superdimension $m |n$ in the sense of Definition \ref{d0819}).
 \ede

\begin{rema} \label{r7gg8} \leavevmode\\
 \ \ \ 
Let $Z^\circledS$ be 
a Deligne-Mumford    superstack.
Then, the structure sheaf $\mcO_{Z^\circledS}$ on $Z^\circledS$ is defined to be   a  super\'{e}tale  sheaf  on $Z^\circledS$ such  that $\Gamma (T^\circledS, \mcO_{Z^\circledS}) := \Gamma (T_b, \mcO_{T^\circledS})$
  for any   superscheme $T^\circledS$ together with a super\'{e}tale  morphism $T^\circledS \migi Z^\circledS$.
Moreover, one may define the notion of an $\mcO_{Z^\circledS}$-supermodule, as usual (cf. the discussion following Remark \ref{r486}).
  \end{rema}

\vspace{1mm}

\subsection{Groupoids in the category of superschemes} \label{S18}
\leavevmode\\
\vspace{-4mm}

Now, we recall that if we are given  a {\it groupoid} $\Gamma$, then it may be described as a  certain collection of data $(U_0, R_0, s_0, t_0, c_0)$, where $U_0$ and $R_0$  denote the sets of {\it objects}  and {\it arrows} of  $\Gamma$ respectively, $s_0$ and $t_0$ denote the {\it source} and {\it target} maps $R_0 \migi U_0$ respectively, and $c_0$ denotes the {\it composition} map $R_0 \times_{t_0, U_0, s_0} R_0 \migi R_0$. 
\vspace{3mm}
\bde\label{D019}\leavevmode\\
\ \ \
A {\bf groupoid in $\mfS \mfc \mfh_{/S_0}^\circledS$} is a collection of data
\begin{align}
R^\circledS \stackrel{}{\rightrightarrows} U^\circledS :=(U^\circledS, R^\circledS, s^\circledS, t^\circledS, c^\circledS),
\end{align}
where 
\begin{itemize}
\item[$\bullet$]
$U^\circledS$ and  $R^\circledS$ are superschemes;
\item[$\bullet$]
 $s^\circledS, t^\circledS : R^\circledS \migi U^\circledS$ and $c^\circledS : R^\circledS \times_{s^\circledS, U^\circledS, t^\circledS} R^\circledS \migi R^\circledS$ are morphisms of superschemes
  \end{itemize}
   such that for any  $T^\circledS \in \mr{Ob} (\mfS \mfc \mfh^{\circledS}_{/S_0})$ the quintuple 
\begin{align}
(R^\circledS \rightrightarrows U^\circledS) (T^\circledS) := (U^\circledS(T^\circledS), R^\circledS(T^\circledS), s^\circledS(T^\circledS), t^\circledS (T^\circledS), c^\circledS(T^\circledS))
\end{align}
 forms  a groupoid (in the above sense) which is functorial with respect to $T^\circledS$.
In a similar vein, one may obtain the definition of a {\bf groupoid in $\mfS \mfc \mfh_{/S_0}$}.
\ede
\vspace{3mm}

As in the usual case of stacks, one may associate,  to each groupoid  
$R^\circledS \rightrightarrows U^\circledS$
 in $\mfS \mfc \mfh^\circledS_{/S_0}$, a superstack 
 \begin{align}
 [R^\circledS \rightrightarrows U^\circledS]
 \end{align} 
 over $S_0$.
 More precisely,  $[R^\circledS \rightrightarrows U^\circledS]$ is  the {\it stackification}
   (with respect to the super\'{e}tale pretopology) of the category fibered in groupoids $ [R^\circledS \rightrightarrows U^\circledS]'$ determined by  $R^\circledS \rightrightarrows U^\circledS$ (i.e., the  fiber of $ [R^\circledS \rightrightarrows U^\circledS]'$ over $T^\circledS \in \mr{Ob} (\mfS \mfc \mfh^\circledS_{/S_0})$ is the groupoid $(R^\circledS \rightrightarrows U^\circledS) (T^\circledS)$ defined above).
Denote by 
\begin{align} \label{e4056}
\pi_{R^\circledS \rightrightarrows U^\circledS}^\circledS : U^\circledS \migi [R^\circledS \rightrightarrows U^\circledS]
\end{align}
 the natural projection.

\vspace{3mm}
\begin{rema} \label{r78} \leavevmode\\
 \ \ \ Let  $Z^\circledS$ be  a Deligne-Mumford superstack.
 \begin{itemize}
 \item[(i)]
One verifies that  there exists an isomorphism $z^\circledS : [R^\circledS \rightrightarrows U^\circledS] \isom Z^\circledS$ of superstacks, where   $R^\circledS \rightrightarrows U^\circledS :=(U^\circledS, R^\circledS, s^\circledS, t^\circledS, c^\circledS)$ is a groupoid   in $\mfS \mfc \mfh_{/S_0}^\circledS$,  such that  the three morphisms   $s^\circledS$,  $t^\circledS$, and $z^\circledS \circ \pi_{R^\circledS \rightrightarrows U^\circledS}^\circledS$ are super\'{e}tale.
Indeed, 
if  $U^\circledS$  is   
 an arbitrary  complete versal family 
  for our $Z^\circledS$,
then one may obtain, by starting with  the data $U^\circledS$, 
the  desired groupoid $[R^\circledS \rightrightarrows U^\circledS]$  in $\mfS \mfc \mfh_{/S_0}^\circledS$
 as follows: 
\begin{itemize}
\vspace{1mm}
\item[$\bullet$]
$R^\circledS := U^\circledS \times_{Z^\circledS} U^\circledS$;
\vspace{1mm}
\item[$\bullet$]
$s^\circledS$ and $t^\circledS$ are  the second  and first  projections $U^\circledS \times_{Z^\circledS} U^\circledS \ (=R^\circledS) \migi U^\circledS$  respectively;
\vspace{1mm}
\item[$\bullet$]
$c^\circledS$ is  the projection 
\begin{align}
& \  (U^\circledS \times_{Z^\circledS} U^\circledS) \times_{s^\circledS, U^\circledS, t^\circledS}    (U^\circledS \times_{Z^\circledS} U^\circledS)  \  (= R^\circledS \times_{s^\circledS, U^\circledS, t^\circledS} R^\circledS)\\
 \migi & \  U^\circledS \times_{Z^\circledS} U^\circledS \ (= R^\circledS) \notag
\end{align}
 into the $(1, 4)$-th factor.
\end{itemize}
We shall refer to such a groupoid $R^\circledS \rightrightarrows U^\circledS$ (together with such an  isomorphism $z^\circledS$) as a {\bf representation}  of $Z^\circledS$.

 \item[(ii)]
 Let $R^\circledS \rightrightarrows U^\circledS$ be as in (i) and denote by
 $R_t \rightrightarrows U_t$ the groupoid in $\mfS \mfc \mfh_{/S_0}$ defined to be $R_t \rightrightarrows U_t := (U_t, R_t, s_t, t_t, c_t)$.
 Then, the isomorphism $z^\circledS$ induces an isomorphism
 \begin{align} \label{E9}
z_t : [R_t \rightrightarrows U_t] \isom Z_t
 \end{align}
of stacks.

\item[(iii)]
Let  $R^\circledS \rightrightarrows U^\circledS$ be  as in (i) again.
Then,  $R_b \rightrightarrows U_b := (U_b, R_b, s_b, t_b, c_b)$ forms a groupoid in $\mfS \mfc \mfh_{/S_0}$ and we obtain a stack
\begin{align} \label{E17}
Z_b := [R_b \rightrightarrows U_b]
\end{align}
together with a morphism $\beta^\circledS_{Z} : Z^\circledS \migi Z_b$.
For any  super\'{e}tale morphism $T^\circledS \migi Z^\circledS$ (where $T^\circledS$ is a superscheme), there exists an \'{e}tale  morphism $T_b \migi Z_b$ which makes the square diagram
\begin{align}
\xymatrix{
T^\circledS \ar[r] \ar[d]_{\beta^\circledS_T}& Z^\circledS \ar[d]^{\beta^\circledS_Z} \\
T_b \ar[r] & Z_b
}
\end{align}
commute and cartesian.
In particular, the structure sheaf $\mcO_{Z_b}$ of $Z_b$ may be identified, via $\beta^\circledS_Z$, with the bosonic part of $\mcO_Z$.
If, moreover, $Z^\circledS$ may be represented by a superscheme $X^\circledS := (X_b, \mcO_{X^\circledS})$,
then $Z_b$ and $\beta_Z^\circledS$ (in the sense of (\ref{E17})) are isomorphic to  $X_b$ and $\beta^\circledS_X$ (in the sense of (\ref{e22})) respectively.
\end{itemize}
  \end{rema}


\vspace{3mm}
\bpr \label{p0207} \leavevmode\\
 \ \ \
 Let $Z^\circledS$ be a superstack  and $m$, $n$ are   nonnegative integers.
 Then, the following three conditions (a), (b), and (c) are equivalent:
 \begin{itemize}
\vspace{1mm}
 \item[(a)]
 $Z^\circledS$ is a supersmooth Deligne-Mumford superstack over $S_0$  of relative superdimension $m|n$;
\vspace{1mm}
\item[(b)]
 $Z^\circledS$ is a Deligne-Mumford superstack for which there exists   a complete versal family  of the form $\langle \underline{U}, \mcE_U \rangle^\circledS$, where $\underline{U}$ denotes  a smooth locally noetherian scheme over $S_0$ of relative dimension $m$ and $\mcE_{\underline{U}}$ denotes  a vector bundle on $\underline{U}$ of rank $n$;
\vspace{1mm}
 \item[(c)]
 $Z^\circledS$ admits a representation $R^\circledS \rightrightarrows U^\circledS := (U^\circledS, R^\circledS, s^\circledS, t^\circledS, c^\circledS)$ satisfying the following properties:
\begin{itemize}
\vspace{1mm}
\item[(c-1)]
$U^\circledS = \langle \underline{U}, \mcE_{\underline{U}} \rangle^\circledS$, where  $\underline{U}$  is a  smooth  scheme over $S_0$ of relative dimension $m$ and $\mcE_{\underline{U}}$ is a rank $n$ vector bundle  on $\underline{U}$;
\vspace{1mm}
\item[(c-2)]
Both $s^\circledS$ and $t^\circledS$ are super\'{e}tale and the morphism $(s_t, t_t) : R_t \migi \underline{U} \times_{S_0} \underline{U}$ ($= U_t \times_{S_0} U_t$) is separated and quasi-compact.
\end{itemize}

  \end{itemize} 
   \epr
\begin{proof}
The equivalence (a) $\Leftrightarrow$ (c)  follows immediately from Proposition \ref{p0201} and the definition of a Deligne-Mumford stack.
The implication (b) $\Rightarrow$ (c) is clear.
Let us consider  (c) $\Rightarrow$ (b).
First, we prove that the diagonal morphism $\varDelta_Z^\circledS : Z^\circledS \migi Z^\circledS \times_{S_0} Z^\circledS$ is representable.
Let $V^\circledS$ be an object in $\mfS \mfc \mfh_{/S_0}^\circledS$,  and let $x^\circledS$ and  $y^\circledS : V^\circledS \migi Z^\circledS$ be morphisms in $\mfS \mfc \mfh_{/S_0}^\circledS$.
To prove that $Z^\circledS \times_{\varDelta_Z^\circledS, Z^\circledS \times_{S_0} Z^\circledS, (x^\circledS, y^\circledS)} V^\circledS$ is in $\mfS \mfc \mfh_{/S_0}^\circledS$,
one may replace (thanks to descent property  in the  super\'{e}tale pretopology)  $V^\circledS$ with its super\'{e}tale covering.
Hence, we suppose, without loss of generality, that
both $x^\circledS$ and $y^\circledS$ may lift to morphisms $\widetilde{x}^\circledS$, $\widetilde{y}^\circledS : V^\circledS \migi U^\circledS$.
Then, we have
\begin{align}
Z^\circledS \times_{\varDelta_Z^\circledS, Z^\circledS \times_{S_0} Z^\circledS, (x^\circledS, y^\circledS)} V^\circledS
 \isom R^\circledS \times_{(s^\circledS, t^\circledS), U^\circledS \times_{S_0} U^\circledS, (\widetilde{x}^\circledS, \widetilde{y}^\circledS)} V^\circledS,
\end{align}
where the right-hand side is evidently an object in $\mfS \mfc \mfh_{/S_0}^\circledS$.
Thus, $\varDelta_Z^\circledS$ is representable.
Moreover,   this representability implies (since both $s^\circledS$ and $t^\circledS$ are  super\'{e}tale) immediately that the projection $\pi_{R^\circledS \rightrightarrows U^\circledS} $ is representable, surjective, and super\'{e}tale.
Finally, by  means of the isomorphism (\ref{E9}), the latter  condition  of (c-2) implies that the diagonal morphism $Z_t \migi Z_t \times_{S_0} Z_t$ is separated and quasi-compact.
This completes the proof of the  implication (c) $\Rightarrow$ (b), and consequently, the proof of Proposition \ref{p0207}.
\end{proof}


\vspace{10mm}
\section{Logarithmic  structures on superschemes} \vspace{3mm}

In this section, we shall   give  briefly  a general formulation of log  superschemes (or more generally,   log superstacks).
The notion of logarithmic  structure on a superscheme, as well as a superstack (cf. Definition \ref{d2a} (i)) 
 is a supersymmetric generalization of the classical notion of logarithmic structure in the sense of J. M. Fontaine and L. Illusie.
 (Basic references for the notion of  {\it logarithmic structure} on a scheme are, e.g.,   ~\cite{KATO} and ~\cite{ILL}.)
One of the most important concepts in log supergeometry is log supersmoothness (cf. Definition \ref{D0519} (ii)).
At the end of this section, we show (cf. Proposition \ref{p0404} (ii) and Corollary \ref{c0404} (i)-(iii)) how
log supersmooth deformations of a log superstack or  a morphism of log superstacks 
are controlled by the sheaf of logarithmic  superderivations (cf. Definition \ref{D0192} for the definition of a  logarithmic  superderivation).

\vspace{5mm}
\subsection{Logarithmic  structures} \label{S21}
\leavevmode\\
\vspace{-5mm}


\bde \label{d2a}\leavevmode\\
 \ \ \ 
\vspace{-5mm}
\begin{itemize}
\item[(i)]
 \ \ \ Let $X^\circledS := (X_b, \mcO_{X^\circledS})$ be a superscheme (resp., a superstack).
A {\bf  logarithmic structure} (or {\bf log structure} for short) on $X^\circledS$ is  a logarithmic structure 
  $\alpha_{X_b} : \mcM_{X_b} \migi \mcO_{X_b}$ on $X_b$ 
(where $\mcM_{X_b}$ denotes an \'{e}tale  sheaf of commutative monoids on $X_b$).

A {\bf log superscheme}  (resp.,  {\bf log superstack}) is a triple 
 \begin{equation}
 Y^{\circledS\mr{log}} := (Y_b, \mcO_{Y^\circledS},   \mcM_{Y_b} \stackrel{ \alpha_{Y_b}}{\migi} \mcO_{Y_b})
 \end{equation}
consisting of  a superscheme (resp., a superstack) $Y^\circledS := (Y_b, \mcO_{Y^\circledS})$  and  a log structure $\alpha_{Y_b}$ on $Y_b$ (hence, of $Y^\circledS$).
We shall refer to $Y_b^\mr{log} := (Y_b, \alpha_{Y_b})$
 as the {\bf underlying log scheme} 
 (resp.,  {\bf underlying log stack})  of $Y^{\circledS \mr{log}}$
 and refer to $Y^\circledS$ as the
  {\bf underlying superscheme} (resp.,   {\bf underlying superstack}) of $Y^{\circledS \mr{log}}$.
Denote by $\beta_Y^{\circledS \mr{log}} : Y^{\circledS \mr{log}} \migi Y_b^\mr{log}$ the  morphism of log superschemes extending $\beta_Y^\circledS$.
\item[(ii)]
 \ \ \ 
 Let $X^{\circledS \mr{log}} := (X_b, \mcO_{X^\circledS}, \alpha_{X_b})$  and $Y^{\circledS \mr{log}} := (Y_b, \mcO_{Y^\circledS}, \alpha_{Y_b})$ be two log superschemes (resp., log superstacks).
A {\bf morphism of log superschemes} (resp., {\bf morphism of log superstacks}) from $Y^{\circledS \mr{log}}$ to $X^{\circledS \mr{log}}$ is 
a triple
\begin{align}
f^{\circledS \mr{log}} := (f_b : Y_b \migi X_b, f^\flat : f_b^* (\mcO_{X^\circledS}) \migi \mcO_{Y^\circledS}, f^\sharp_b : f_b^{-1}(\mcM_{X_b}) \stackrel{}{\migi} \mcM_{Y_b}),
\end{align}
where $f^\circledS := (f_b, f^\flat)$ forms  
a morphism  $Y^\circledS \migi X^\circledS$ between  the underlying superschemes (resp., underlying superstacks) and  $f^\mr{log}_b := (f_b, f_b^\sharp)$ forms a morphism $Y_b^\mr{log} \migi X_b^\mr{log}$ between the underlying log schemes (resp., underlying log stacks).  
\end{itemize}
 \ede

\vspace{3mm}
\bde \label{d2b}\leavevmode\\
 \ \ \ 
An {\bf fs log superscheme} (resp., {\bf fs log superstack}) is a log superscheme (resp., log superstack) whose underlying log scheme (resp., underlying log stack) is fine and saturated.
 \ede
\vspace{3mm}

Let  $\underline{X}^\mr{log} := (\underline{X}, \alpha_{\underline{X}})$ be an fs  log scheme over $S_0$ and $\mcE$  a coherent $\mcO_{\underline{X}}$-module.
Then, we shall write
\begin{align} \label{E13}
\langle \underline{X}, \mcE \rangle^{\circledS\mr{log}}
\end{align}
 for the log superscheme defined to be $\langle \underline{X}, \mcE \rangle$ (cf.  (\ref{E12})) equipped with the log structure pulled-back from $\underline{X}^\mr{log}$ via
 $\langle \beta \rangle^\circledS_{\underline{X}, \mcE}$ (cf. (\ref{E14})).

\vspace{3mm}
\bde \label{d2c}\leavevmode\\
 \ \ \ 
Let $X^{\circledS \mr{log}}$ and $Y^{\circledS \mr{log}}$ be   log superschemes (resp., log superstacks) and $f^{\circledS \mr{log}} : Y^{\circledS \mr{log}} \migi X^{\circledS \mr{log}}$  a morphism of log superschemes (resp., a morphism of log superstacks).
We shall say that $ f^{\circledS \mr{log}}$   is  {\bf  strict   super\'{e}tale} (resp., {\bf a strict closed immersion}) if 
$f^\circledS$ is super\'{e}tale (resp., a closed immersion) and  $f^\mr{log}_b$ is strict, in the sense of ~\cite{ILL}, \S\,1.2.
 \ede
\vspace{3mm}

We shall write
\begin{align} \label{EE13}
 \mfS \mfc \mfh^{\circledS \mr{log}}_{/S_0}
\end{align}
for the category
whose {\it objects} are  fs log superschemes   and whose {\it morphisms} are morphisms of  log superschemes.
Also, write $\mfS \mfc \mfh_{/S_0}^\mr{log}$ for the full subcategory of $\mfS \mfc \mfh^{\circledS \mr{log}}_{/S_0}$ consisting of fs log schemes (i.e., fs log superschemes  $X^{\circledS \mr{log}}$ with $\mcO_{X_f} =0$).
The fiber products and finite coproducts in  $\mfS \mfc \mfh^{\circledS \mr{log}}_{/S_0}$ exist, and 
$\mfS \mfc \mfh^{\circledS \mr{log}}_{/S_0}$ admits the Grothendieck pretopology  given by strict  super\'{e}tale morphisms;  we shall refer to it as the  {\bf strict  super\'{e}tale pretopology}.
In a natural manner,  any log superstack may be thought of as a stack
  over $ \mfS \mfc \mfh^{\circledS \mr{log}}_{/S_0}$ with respect to the strict  super\'{e}tale pretopology.

\vspace{5mm}
\subsection{Logarithmic superdifferentials} \label{S22}
\leavevmode\\
\vspace{-4mm}


Let  $S^{\circledS \mr{log}} : = (S_b, \mcO_{S^\circledS}, \alpha_{S_b})$ and $X^{\circledS \mr{log}} := (X_b, \mcO_{X^\circledS}, \alpha_{X_b})$ be fs log superschemes and 
$f^{\circledS \mr{log}} \ (:= (f_b, f^\flat,  f_b^\sharp)) : X^{\circledS \mr{log}} \migi S^{\circledS \mr{log}}$ a morphism of log superschemes.
In the following, we shall  define (in a functorial manner)   a ``{\it log super}" analogue  of the  sheaf  of relative differential $1$-forms, i.e., an $\mcO_{X^\circledS}$-supermodule $\Omega_{X^{\circledS \mr{log}}/S^{\circledS \mr{log}}}$ 
  together with the universal derivation $d : \mcO_{X^\circledS} \migi \Omega_{X^{\circledS \mr{log}}/S^{\circledS\mr{log}}}$ defined  as follows.
Let us write $\varDelta_X^\circledS : X^\circledS \migi X^\circledS \times_{S^\circledS} X^\circledS$ for  the diagonal morphism and write   $\mcJ := \mr{Ker} (\mcO_{X^\circledS  \times_{S^\circledS} X^\circledS} \migi \varDelta_{X*}^\circledS(\mcO_{X^\circledS}))$.
Then, we shall define
\begin{align}
\Omega_{X^\circledS/S^\circledS} := \varDelta_X^{\circledS *}(\mcJ/\mcJ^2)
\end{align}
 and write
 $d : \mcO_{X^\circledS} \migi \Omega_{X^\circledS/S^\circledS}$ for the $f^{-1}_b(\mcO_{S^\circledS})$-linear  morphism  given by assigning $a \mapsto d(a) := \overline{(a \otimes 1 - 1\otimes a)}$  for any local section $a \in \mcO_{X^\circledS}$.
For example, if $X^{\circledS} = \mbA^{m |n}_{S^\circledS}$, then we have 
\begin{align} \label{E44}
\Omega_{X^\circledS/S^\circledS} \cong (\bigoplus_{i=1}^m \mcO_{X^\circledS} d (t_i)) \oplus  (\bigoplus_{i=1}^n \mcO_{X^\circledS} d (\psi_i)),
\end{align}
 where $d (t_i)$ ($i =1, \cdots, m$) are bosonic elements in $\Omega_{X^\circledS/S^\circledS}$ and $d (\psi_i)$ ($i = 1, \cdots, n$) are fermionic  elements.

Moreover, let us define  the $\mcO_{X^{\circledS}}$-supermodule $\Omega_{X^{\circledS \mr{log}}/S^{\circledS \mr{log}}}$ to be 
\begin{align}
\Omega_{X^{\circledS \mr{log}}/S^{\circledS \mr{log}}} := (\Omega_{X^\circledS/S^\circledS} \oplus (\mcO_{X^\circledS} \otimes_\mbZ \mcM^\mr{gr}_{X_b}))/\mcN,
\end{align}
where 
\begin{itemize}
\item[(i)]
$ \mcM^\mr{gr}_{X_b}$ denotes the groupification of $\mcM_{X_b}$ whose local sections are  bosonic (hence we obtain the  $\mcO_{X^\circledS}$-supermodule $\mcO_{X^\circledS} \otimes_\mbZ \mcM^\mr{gr}_{X_b}$);
\item[(ii)]
  $\mcN$ denotes the $\mcO_{X^\circledS}$-subsupermodule generated  locally by local sections of  the following forms:
\begin{itemize}
\vspace{1mm}
\item[$\bullet$]
$(d (\alpha_{X_b} (a)), 0) - (0, \alpha_{X_b} (a) \otimes a)$ with $a \in \mcM_{X_b}$;
\item[$\bullet$]
$(0, 1\otimes a)$ with $a \in \mr{Im} (f_b^{-1} (\mcM_{S_b}) \stackrel{f_b^\sharp}{\migi} \mcM_{X_b})$.
\end{itemize}
\vspace{1mm}
\end{itemize}
The class of $(0, 1\otimes a)$ for $a \in \mcM_{X_b}$ in $\Omega_{X^{\circledS \mr{log}}/S^{\circledS \mr{log}}}$ is denoted by $d \mr{log} (a)$.
Finally, we write
\begin{equation}
\mcT_{X^{\circledS \mr{log}}/S^{\circledS \mr{log}}} := \Omega_{X^{\circledS \mr{log}}/S^{\circledS \mr{log}}}^\vee.
\end{equation}
i.e., the dual $\mcO_{X^\circledS}$-supermodule  of $\Omega_{X^{\circledS \mr{log}}/S^{\circledS \mr{log}}}$.

The following Propositions \ref{p0404f} and \ref{p04304} may be verified  immediately.

\vspace{3mm}
\bpr \label{p0404f} \leavevmode\\
 \ \ \
Let us consider a cartesian square diagram
\begin{align}
\xymatrix{
Y^{\circledS \mr{log}} \ar[r]^{f^{\circledS \mr{log}}} \ar[d] \ar@{}[rd]|{\Box} & X^{\circledS \mr{log}} \ar[d] \\
T^{\circledS \mr{log}} \ar[r] & S^{\circledS \mr{log}}
}
\end{align}
in $\mfS \mfc \mfh_{/S_0}^{\circledS \mr{log}}$.
Then, the natural $\mcO_{Y^\circledS}$-linear morphism
\begin{align}
f^{\circledS *} (\Omega_{X^{\circledS \mr{log}}/S^{\circledS \mr{log}}}) \migi \Omega_{Y^{\circledS \mr{log}}/T^{\circledS \mr{log}}}
\end{align}
is an isomorphism.
     \epr

\vspace{3mm}
\bpr \label{p04304} \leavevmode\\
\vspace{-5mm}
\begin{itemize}
\item[(i)]
Let $X^{\circledS \mr{log}}$ and $Y^\circledS$ be fs log superschemes  and $f^{\circledS \mr{log}} : Y^\circledS \migi X^\circledS$
  a morphism of log  superschemes.
Then, there exists an exact sequence
\begin{align} \label{E20}
f^{\circledS *} (\Omega_{X^{\circledS \mr{log}}/S^{\circledS \mr{log}}}) \migi \Omega_{Y^{\circledS \mr{log}}/S^{\circledS \mr{log}}} \migi \Omega_{Y^{\circledS \mr{log}}/X^{\circledS \mr{log}}} \migi 0
\end{align}
of $\mcO_{Y^{\circledS}}$-supermodule.
\item[(ii)]
Suppose further  that  $f^{\circledS \mr{log}}$ is strict super\'{e}tale.
Then,
$\Omega_{Y^{\circledS \mr{log}}/X^{\circledS \mr{log}}} =0$ and the first arrow $f^{\circledS *} (\Omega_{X^{\circledS \mr{log}}/S^{\circledS \mr{log}}}) \migi \Omega_{Y^{\circledS \mr{log}}/S^{\circledS \mr{log}}}$ in (\ref{E20})
is an isomorphism.
 \end{itemize}
     \epr

\vspace{5mm}
\subsection{Logarithmic superderivations} \label{S23}
\leavevmode\\
\vspace{-4mm}

Let $S^{\circledS \mr{log}}$,  $X^{\circledS \mr{log}}$, and $f^{\circledS \mr{log}} : X^{\circledS \mr{log}} \migi S^{\circledS \mr{log}}$ be as at the beginning of the previous subsection.

\bde\label{D0192}\leavevmode\\
\ \ \
Let $\mcE$ be an  $\mcO_{X^\circledS}$-supermodule.
A {\bf logarithmic superderivation of $(\mcO_{X^\circledS}, \mcM_{X_b})$ (over $S^{\circledS\mr{log}}$) with value in $\mcE$} is a pair $\partial := (D, \delta)$, where
\begin{itemize}
\item[$\bullet$]
 $D$ is a superderivation $\mcO_{X^\circledS} \migi \mcE$ over $S^\circledS$, i.e., an $f_b^{-1} (\mcO_{S^\circledS})$-linear morphism satisfying that 
 \begin{align}
 D (a \cdot b) = D (a) \cdot b + (-1)^{|D| \cdot |a|} a \cdot D(b)
 \end{align}
  for any local sections $a, b  \in \mcO_X$ (where $|D|$ denotes the parity of $D$); 
 \item[$\bullet$]
 $\delta$ is a monoid homomorphism $\mcM_{X_b} \migi \mcE$ such that 
 \begin{align}
 D (\alpha_{X_b} (m)) = \alpha_{X_b} (m) \cdot \delta (m)
 \end{align}
  for any local section $m \in \mcM_{X_b}$;
\item[$\bullet$]
$D (f^{-1}_b (b)) = \delta (f^\sharp_b (n)) =0$ for any sections $b \in \mcO_{S^\circledS}$ and $n \in \mcM_{S_b}$. 
\end{itemize}
If $\partial := (D, \delta)$ is a logarithmic superderivation, then we usually just write $\partial (a)$ and $\partial(m)$ (where $a \in \mcO_{X^\circledS}$ and $m \in \mcM_{X_b}$) for $D (a)$ and $\delta (m)$ respectively.
\ede

\vspace{3mm}
\begin{rema} \label{r408} \leavevmode\\
 \ \ \ 
Let $\mcE$ be an $\mcO_{X^\circledS}$-supermodule.
Denote by
\begin{align} \label{E45}
\mr{Def}_{S^{\circledS}} (X^{\circledS \mr{log}}; \mcE)
\end{align}
 the set of logarithmic superderivations of $(\mcO_{X^\circledS}, \mcM_{X_b})$ over $S^{\circledS \mr{log}}$ with value in $\mcE$.
The structure of $\mcO_{X^\circledS}$-supermodule on $\mcE$ gives rise to a structure of $\Gamma (X_b, \mcO_{X^\circledS})$-supermodule on $\mr{Def}_{S^{\circledS}} (X^{\circledS \mr{log}}; \mcE)$.
In particular, $\mr{Def}_{S^{\circledS}} (X^{\circledS \mr{log}}; \mcE)$ decomposes as 
\begin{align}
\mr{Def}_{S^{\circledS}} (X^{\circledS \mr{log}}; \mcE)  = \mr{Def}_{S^{\circledS}} (X^{\circledS \mr{log}}; \mcE)_b \oplus \mr{Def}_{S^{\circledS}} (X^{\circledS \mr{log}}; \mcE)_f.
\end{align}
It is clear that there exists a universal logarithmic  superderivation
\begin{align} \label{E46}
d \in  \mr{Def}_{S^{\circledS}} (X^{\circledS \mr{log}}; \Omega_{X^{\circledS \mr{log}}/S^{\circledS \mr{log}}})_b.
\end{align}
That is to say,
 the morphism
 \begin{align} \label{E47}
 \mr{Hom}_{\mcO_{X^\circledS}} (\Omega_{X^{\circledS \mr{log}}/S^{\circledS \mr{log}}}, \mcE) \isom & \    \mr{Def}_{S^{\circledS}} (X^{\circledS \mr{log}}; \mcE) \\
 h \hspace{10mm} \mapsto  & \hspace{7mm} h \circ d \notag 
 \end{align}
 is an isomorphism 
 of  $\Gamma (X_b, \mcO_{X^\circledS})$-supermodules.
In particular, (since the isomorphism (\ref{E47}) is compatible with restriction to each open subscheme of $X_b$) the case of $\mcE = \mcO_{X^\circledS}$ implies that
the dual $\mcT_{X^{\circledS \mr{log}}/S^{\circledS \mr{log}}}$ of $\Omega_{X^{\circledS \mr{log}}/S^{\circledS \mr{log}}}$  is isomorphic to the sheaf 
given by assigning $U \mapsto \mr{Def}_{S^{\circledS}} (X^{\circledS \mr{log}} \times_{X_b} U; \mcE |_{U})$ (for any open subscheme $U$ of $X_b$).
 By taking account of  ~\cite{Og}, Proposition 1.1.7, one verifies that 
 $\mcT_{X^{\circledS \mr{log}}/S^{\circledS \mr{log}}}$  admits a structure of Lie superalgebra over $f_b^{-1}(\mcO_{S^\circledS})$ with bracket operation  given by
 \begin{align}
 [\partial_1, \partial_2] := (D_1\circ D_2 - (-)^{|D_1| \cdot |D_2|} D_2 \circ D_1, D_1 \circ \delta_2  - (-)^{|D_1| \cdot |D_2|} D_2 \circ \delta_1)
 \end{align}
for any homogenous  logarithmic superderivations $\partial_1 := (D_1, \delta_1)$ and  $\partial_2 := (D_2, \delta_2)$.
 \end{rema}
\vspace{3mm}

\begin{rema} \label{r4018} \leavevmode\\
 \ \ \ 
 The discussions in \S\,2.2 and \S\,2.3 (especially, 
  Propositions \ref{p0404f} and \ref{p04304}) generalize naturally to the case where $X^{\circledS \mr{log}}$ is a log superstack.
 In fact,  $\Omega_{X^{\circledS \mr{log}}/S^{\circledS \mr{log}}}$ is constructed in such a way that if $t^{\circledS \mr{log}}  :  T^{\circledS \mr{log}} \migi X^{\circledS \mr{log}}$ (where $T^{\circledS \mr{log}}$ is a superscheme) is
 a strict super\'{e}tale morphism, then we have a functorial (with respect to $T^{\circledS \mr{log}}$) isomorphism $t^{\circledS *}(\Omega_{X^{\circledS \mr{log}}/S^{\circledS \mr{log}}}) \cong \Omega_{T^{\circledS \mr{log}}/S^{\circledS \mr{log}}}$.
 \end{rema}
\vspace{5mm}

\subsection{Log supersmooth morphisms} \label{S24}
\leavevmode\\
\vspace{-4mm}

Let $S^{\circledS \mr{log}}$ be an fs log superscheme,   
$X^{\circledS \mr{log}}$  an   fs log superstack, and $f^{\circledS \mr{log}} : X^{\circledS \mr{log}}\migi S^{\circledS \mr{log}}$  a morphism of log superstacks.
\vspace{3mm}
\bde\label{D0519}\leavevmode\\
\ \ \ 
Let   $m$, $n$  be nonnegative integers.

\begin{itemize}
\item[(i)]
An  {\bf $(m|n)$-chart} on  $X^{\circledS \mr{log}}/S^{\circledS \mr{log}}$ is a triple 
\begin{align}
(Y^{\circledS \mr{log}}, U^\mr{log}, \eta^{\circledS \mr{log}}),
\end{align}
where
\begin{itemize}
\vspace{2mm}
\item[$\bullet$]
$Y^{\circledS \mr{log}}$ is 
an affine  fs log  superscheme 
together with a strict super\'{e}tale morphism $Y^{\circledS \mr{log}} \migi X^{\circledS \mr{log}}$ over $S_0$;
\vspace{2mm}
\item[$\bullet$]
$U^\mr{log}$ is an fs log affine scheme
 together with an integral log smooth  morphism  $U^\mr{log} \migi S_b^\mr{log}$ of relative dimension $m$; 
\vspace{2mm}
\item[$\bullet$]
$\eta^{\circledS \mr{log}}$ is an isomorphism $Y^{\circledS \mr{log}} \isom U^\mr{log} \times_{S_b} \mbA^{0|n}_{S^{\circledS}}$ over $S^{\circledS \mr{log}}$.
\vspace{2mm}
\end{itemize}
\item[(ii)]
We shall say that $X^{\circledS \mr{log}}$ is {\bf log supersmooth over $S^{\circledS \mr{log}}$  of relative superdimension $m|n$} if
there exists
a collection $\{ (Y^{\circledS \mr{log}}_\gamma, U^\mr{log}_\gamma, \eta^{\circledS \mr{log}}_\gamma)\}_\gamma$ of $(m|n)$-charts  on $X^{\circledS \mr{log}}/S^{\circledS \mr{log}}$ for which  the morphism $\coprod_\gamma Y^{\circledS \mr{log}}_\gamma \migi X^{\circledS \mr{log}}$ is a  strict super\'{e}tale covering of $X^{\circledS \mr{log}}$.
\end{itemize}
\ede
\vspace{3mm}

\begin{rema} \label{r4078} \leavevmode\\
 \ \ \ 
It is clear that  if both $S^{\circledS \mr{log}}$ and $X^{\circledS \mr{log}}$ has trivial structures, then   $X^{\circledS \mr{log}}$ is log supersmooth over $S^{\circledS \mr{log}}$ of relative superdimension $m|n$ if and only if $X^\circledS$ is supersmooth  of relative superdimension $m|n$, in the sense of Definition \ref{d08gg9}, (ii).
 \end{rema}

\vspace{2mm}
\bpr \label{p0404} \leavevmode\\
 \ \ \
Suppose that $X^{\circledS \mr{log}}$ is log supersmooth over $S^{\circledS \mr{log}}$ of relative superdimension $m |n$ for some nonnegative integers $m$, $n$. Then, the following assertions hold.
\begin{itemize}
\item[(i)]
The $\mcO_{X^\circledS}$-supermodule $\Omega_{X^{\circledS \mr{log}}/S^{\circledS \mr{log}}}$ is a supervector bundle  of superrank $m|n$.

\vspace{1mm}
\item[(ii)]
Let us consider a commutative square diagram
\begin{align}
\xymatrix{
T^{\circledS \mr{log}}  \ar[d]_{i^{\circledS \mr{log}}} \ar[r]^{t_X^{\circledS \mr{log}}}& X^{\circledS \mr{log}} \ar[d]^{f^{\circledS \mr{log}}} \\
\widetilde{T}^{\circledS \mr{log}} \ar[r]^{\widetilde{t}_S^{\circledS \mr{log}}} & S^{\circledS \mr{log}},}
\end{align}
where $T^\circledS$ is affine and $i^{\circledS \mr{log}}$ is a strict closed immersion defined by a square nilpotent  superideal $\mcJ$ of $\mcO_{\widetilde{T}^\circledS}$. 
(Hence, $\mcJ$ may be thought of as an $\mcO_{T^\circledS}$-supermodule.)
We shall write 
\begin{align} \label{e04889}
\mcF  := \mcH om_{\mcO_{T^\circledS}} (t_X^{\circledS *}(\Omega_{X^{\circledS \mr{log}}/S^{\circledS \mr{log}}}), \mcJ)_{b}  \ \big(\cong (t_X^{\circledS *}(\mcT_{X^{\circledS \mr{log}}/S^{\circledS \mr{log}}} )\otimes \mcJ)_b \big)
\end{align}
i.e., 
 the $\mcO_{T_b}$-submodule of $\mcH om_{\mcO_{T^\circledS}} (t_X^{\circledS *}(\Omega_{X^{\circledS \mr{log}}/S^{\circledS \mr{log}}}), \mcJ)$ consisting of $\mcO_{T^\circledS}$-linear homomorphisms   
$t_X^{\circledS *}(\Omega_{X^{\circledS \mr{log}}/S^{\circledS \mr{log}}}) \migi  \mcJ$
 of even parity.
Also, we shall write 
\begin{align}
\mcD e f_{\widetilde{T}^{\circledS \mr{log}}}(t_X^{\circledS \mr{log}})
\end{align}
 for 
 the strict super\'{e}tale sheaf on $\widetilde{T}^{\circledS \mr{log}}$ which, to any strict super\'{e}tale morphism $\alpha^{\circledS \mr{log}} : \widetilde{T}_1^{\circledS \mr{log}} \migi \widetilde{T}^{\circledS \mr{log}}$, assigns 
 the set of 
 morphisms $\widetilde{t}_{1, X}^{\circledS \mr{log}} : \widetilde{T}^{\circledS \mr{log}} \migi X^{\circledS \mr{log}}$ which makes the diagram
 \begin{align}
\xymatrix{
  T^{\circledS \mr{log}} \times_{\widetilde{T}^{\circledS \mr{log}}} \widetilde{T}_1^{\circledS \mr{log}}  \ar[d]_{i^{\circledS \mr{log}}} \ar[r]^{}& X^{\circledS \mr{log}} \ar[d]^{f^{\circledS \mr{log}}} \\
 \widetilde{T}_1^{\circledS \mr{log}} \ar[ru]_{\widetilde{t}_{1, X}^{\circledS \mr{log}}}\ar[r]_{\widetilde{t}_S^{\circledS \mr{log}}\circ \alpha^{\circledS \mr{log}} } & S^{\circledS \mr{log}},}
\end{align}
commute,
where the upper horizontal arrow denotes the composite of $t_X^{\circledS \mr{log}}$ and  the  projection $ T^{\circledS \mr{log}} \times_{\widetilde{T}^{\circledS \mr{log}}} \widetilde{T}_1^{\circledS \mr{log}} \migi T^{\circledS \mr{log}}$ to the first factor.

 
 Then, $\mcD e f_{\widetilde{T}^{\circledS \mr{log}}}(t_X^{\circledS \mr{log}})$ is nonempty (i.e., admits locally a section), and moreover,  admits canonically a structure of affine space
 \begin{align} \label{E49}
 \mcD e f_{\widetilde{T}^{\circledS \mr{log}}}(t_X^{\circledS \mr{log}})  \times i_{b*}(\mcF) & \migi \mcD e f_{\widetilde{T}^{\circledS \mr{log}}}(t_X^{\circledS \mr{log}}) \\
 (\widetilde{t}^{\circledS \mr{log}}_X, \partial) \hspace{5mm} & \mapsto \ \  \  \widetilde{t}^{\circledS \mr{log}}_X \boxplus^\dagger \partial \notag
 \end{align}
 
  modeled on $i_{b*}(\mcF)$.
\end{itemize}
     \epr
\begin{proof}
Assertion (i) follows from (\ref{E44}), Proposition \ref{p04304} (ii),  and ~\cite{KATO}, Proposition (3.10).
Next, we shall prove the former assertion of (ii), i.e., that $\mcD e f_{\widetilde{T}^{\circledS \mr{log}}}(t_X^{\circledS \mr{log}})$ is nonempty.
After possibly replacing $\widetilde{T}^{\circledS \mr{log}}$ with its strict super\'{e}tale covering, one may assume, without loss of generality,  that $X^{\circledS \mr{log}} = U^\mr{log} \times_{S_b} \mbA_{S^\circledS}^{0|n}$ for some fs log affine scheme $U^\mr{log}$ together with an integral  log smooth morphism  $f_U^\mr{log} : U^\mr{log} \migi S_b^\mr{log}$ of relative dimension $m$.
Consider  the commutative diagram
\begin{align}
\begin{CD}
T_b^{\mr{log}}
@> \underline{t}_X^\mr{log} >>  U^\mr{log} 
\\
@V i_b^{\mr{log}} VV @VV f_U^\mr{log} V
\\
\widetilde{T}_b^{\mr{log}} 
@>>  ({\widetilde{t}_S})_b^{\mr{log}} > S_b^\mr{log},
\end{CD}
\end{align}
where the upper horizontal arrow $\underline{t}_X^\mr{log}$ denotes 
the composite $(t_X)_b^\mr{log} : T_b^\mr{log} \migi X_b^\mr{log}$ and the natural projection $X_b^\mr{log} \migi U^\mr{log}$.
Since $f^\mr{log}_U$ is log smooth,
 there exists a morphism $\widetilde{t}_{b, U}^\mr{log} : \widetilde{T}_b^\mr{log} \migi U^\mr{log}$
  such that 
  $\widetilde{t}_{b, U}^\mr{log} \circ i_b^\mr{log} = \underline{t}_X^\mr{log}$
   and
    $f_U^\mr{log} \circ \widetilde{t}_{b, U}^\mr{log} =  (\widetilde{t}_S)_b^\mr{log}$.
On the other hand, let us consider  the functorial  bijection (\ref{FF04})  obtained in Proposition \ref{P0} (of the case where $(m, n) =(0, 1)$). 
Then,  the composite of  $t_X^{\circledS \mr{log}} : T^{\circledS \mr{log}} \stackrel{}{\migi} U^\mr{log} \times_{S_b} \mbA^{0|1}_{S^\circledS}$ 
and the projection 
$U^\mr{log} \times_{S_b} \mbA^{0|1}_{S^\circledS} \stackrel{}{\migi} \mbA^{0 |1}_{S^\circledS}$ to the second factor  extends, strict super\'{e}tale locally on $\widetilde{T}^{\circledS \mr{log}}$,
 to a morphism $\widetilde{t}_\mbA^{\circledS \mr{log}} : \widetilde{T}^{\circledS \mr{log}} \migi \mbA^{0|1}_{S^\circledS}$.
Thus, the morphism 
\begin{align}
(\widetilde{t}_{b, U}^\mr{log} \circ \beta^{\circledS \mr{log}}_{\widetilde{T}^{}}, \widetilde{t}_\mbA^{\circledS \mr{log}}) : \widetilde{T}^{\circledS \mr{log}} \migi U^\mr{log} \times_{S_b} \mbA_{S^\circledS}^{0|1}
\end{align}
 determines a section of $\mcD e f_{\widetilde{T}^{\circledS \mr{log}}}(t_X^{\circledS \mr{log}})$.

Next, we shall prove the latter assertion of (ii).
One may assume, without loss of generality, that  there exists an element $\widetilde{t}_{X, 0}^{\circledS \mr{log}} \in \Gamma (\widetilde{T}^{}_b, \mcD e f_{\widetilde{T}^{\circledS \mr{log}}}(t_X^{\circledS \mr{log}}))$.
Since we have an isomorphism
\begin{align}
\Gamma (\widetilde{T}^{}_b, i_{b*}(\mcF)) \ ( \isom \mr{Hom}_{\mcO_{X^{\circledS}}} (\Omega_{X^{\circledS \mr{log}}/S^{\circledS \mr{log}}}, t_{X*}^{\circledS} (\mcJ))_b) 
\isom \mr{Der}_{S^{\circledS \mr{log}}} (X^{\circledS \mr{log}}; t_{X*}^{\circledS} (\mcJ))_b
\end{align}
(cf. \ref{E47}), 
 it suffices to construct  a functorial (with respect to $\widetilde{T}^{\circledS \mr{log}}$) bijection 
\begin{align} \label{E48}
\mr{Der}_{S^{\circledS \mr{log}}} (X^{\circledS \mr{log}};  (\widetilde{t}_{X, 0}^\circledS)_* (\mcJ))_b \isom \Gamma (\widetilde{T}^{}_b, \mcD e f_{\widetilde{T}^{\circledS \mr{log}}}(t_X^{\circledS \mr{log}})).
\end{align}
Denote by $\widetilde{t}_{X, 0}^{\circledS \dagger} : \mcO_{X^\circledS} \migi t_{X*}^\circledS (\mcO_{\widetilde{T}^\circledS})$ and $\widetilde{t}_{X, 0}^{\circledS \ddagger}  : \mcM_{X_b} \migi t_{X*}^\circledS (\mcM_{\widetilde{T}_b})$ the morphism arising naturally from $\widetilde{t}_{X, 0}^{\circledS}$ (where we consider both $\mcO_{\widetilde{T}^\circledS}$ and $\mcM_{\widetilde{T}_b}$ as sheaves on $T_b$ via the underlying  homeomorphism between topological spaces  of $i^{\circledS \mr{log}}$).
Let us take  an element
\begin{align}
\partial := (\mcO_{X^\circledS} \stackrel{D}{\migi}  (\widetilde{t}_{X, 0}^\circledS)_* (\mcJ), \mcM_{X_b} \stackrel{\delta}{\migi} (\widetilde{t}_{X, 0}^\circledS)_* (\mcJ))
\end{align}
of $ \mr{Der}_{S^{\circledS \mr{log}}} (X^{\circledS \mr{log}};  (\widetilde{t}_{X, 0}^\circledS)_* (\mcJ))_b$.
By applying  the inclusions 
\begin{align}
(\widetilde{t}_{X, 0}^\circledS)_* (\mcJ) \migiincl (\widetilde{t}_{X, 0}^\circledS)_*(\mcO_{\widetilde{T}^\circledS})
 \ \  \text{and} \ \ 
(\widetilde{t}_{X, 0}^\circledS)_* (1 + \mcJ)  \migiincl  (\widetilde{t}_{X, 0}^\circledS)_* (\mcM_{\widetilde{T}_b}),
\end{align}
 one may obtain  two maps $\widetilde{t}_{X, 0}^{\circledS \dagger} + D$ and $\widetilde{t}_{X, 0}^{\circledS \ddagger} + \delta$ given by 
\begin{align}
\widetilde{t}_{X, 0}^{\circledS \dagger} + D : & \ \mcO_{X^\circledS} \migi (\widetilde{t}_{X, 0}^\circledS)_*(\mcO_{\widetilde{T}^\circledS}),  &  \widetilde{t}_{X, 0}^{\circledS \ddagger} + \delta : & \  \mcM_{X_b} \migi (\widetilde{t}_{X, 0}^\circledS)_* (\mcM_{\widetilde{T}_b})  \\
& \  \  a \ \  \ \mapsto \ \ \  a + D (a) &   &  \ \ \ b  \ \ \  \mapsto \ \ \  1 + \delta (b) \notag
\end{align}
for any  local sections $a \in \mcO_{X^\circledS}$, $b \in \mcM_{X_b}$.
By the definition of  a logarithmic  superderivation, the pair $(\widetilde{t}_{X, 0}^{\circledS \dagger} + D, \widetilde{t}_{X, 0}^{\circledS \ddagger} + \delta)$ determines 
a new morphism  $\widetilde{t}_{X, 0}^{\circledS} \boxplus^\dagger \partial : \widetilde{T}^{\circledS \mr{log}} \migi X^{\circledS \mr{log}}$ in 
$\mcD e f_{\widetilde{T}^{\circledS \mr{log}}}(t_X^{\circledS \mr{log}})$.
One verifies immediately that this assignment $\partial \mapsto  \widetilde{t}_{X, 0}^{\circledS} \boxplus^\dagger \partial $ determines  the desired  bijection  (\ref{E48}).
This completes the proof of Proposition \ref{p0404}.
\end{proof}

\vspace{5mm}

\subsection{Log supersmooth liftings}\label{S25}
\leavevmode\\
\vspace{-5mm}

\bde\label{D0719}\leavevmode\\
\ \ \ 
Let 
$X^{\circledS \mr{log}}$ and  $S^{\circledS \mr{log}}$
 be as in Proposition \ref{p0404}.
Also, let $i_S^{\circledS \mr{log}} : S^{\circledS \mr{log}} \migi \widetilde{S}^{\circledS \mr{log}}$ be a strict closed immersion  
determined  by a   nilpotent superideal $\mcJ$ on $\mcO_{\widetilde{S}^\circledS}$.
\begin{itemize}
\item[(i)]
By a {\bf log supersmooth lifting} of $X^{\circledS \mr{log}}$ over $\widetilde{S}^{\circledS \mr{log}}$, we mean a triple 
\begin{align}
\widetilde{\mbX} :=  (\widetilde{X}^{\circledS \mr{log}}, \widetilde{f}^{\circledS \mr{log}}, i_X^{\circledS \mr{log}})
\end{align}
 consisting of  a log   superstack $\widetilde{X}^{\circledS \mr{log}}$,  a  log  supersmooth morphism $\widetilde{f}^{\circledS \mr{log}} : \widetilde{X}^{\circledS \mr{log}} \migi \widetilde{S}^{\circledS \mr{log}}$, and a strict closed immersion $i^{\circledS \mr{log}}_X : X^{\circledS \mr{log}} \migi \widetilde{X}^{\circledS \mr{log}}$ which make   the square  diagram
\begin{align}
\xymatrix{
 X^{\circledS \mr{log}} \ar[d]_{f^{\circledS \mr{log}}} \ar[r]^{i_X^{\circledS \mr{log}}}
  & \widetilde{X}^{\circledS \mr{log}} \ar[d]^{\widetilde{f}^{\circledS \mr{log}}}\\
 S^{\circledS \mr{log}} \ar[r]_{i_S^{\circledS \mr{log}}} & \widetilde{S}^{\circledS \mr{log}}
}
\end{align}
commute and cartesian.
\item[(ii)]
Let $\widetilde{\mbX}_l := (\widetilde{X}_l^{\circledS \mr{log}}, \widetilde{f}_l^{\circledS \mr{log}}, i_{X_l}^{\circledS \mr{log}})$  ($l =1,2$)
be log supersmooth liftings of $X^{\circledS \mr{log}}$ over $\widetilde{S}^{\circledS \mr{log}}$.
An {\bf isomorphism of log supersmooth liftings} from $\widetilde{\mbX}_1$ to  $\widetilde{\mbX}_2$ 
is an isomorphism $j^{\circledS \mr{log}} : \widetilde{X}^{\circledS \mr{log}}_1 \isom \widetilde{X}^{\circledS \mr{log}}_2$ such that $\widetilde{f}^{\circledS \mr{log}}_2 \circ j^{\circledS \mr{log}} = \widetilde{f}^{\circledS\mr{log}}_1$ and $\widetilde{i}^{\circledS \mr{log}}_{X_2} \circ j^{\circledS \mr{log}} = \widetilde{i}^{\circledS \mr{log}}_{X_1}$.
\end{itemize}
\ede
\vspace{3mm}

\begin{rema} \label{r4} \leavevmode\\
 \ \ \ 
 Suppose that we are given a log supersmooth lifting $\widetilde{\mbX} := (\widetilde{X}^{\circledS \mr{log}}, \widetilde{f}^{\circledS \mr{log}}, i_{X}^{\circledS \mr{log}})$ of $X^{\circledS \mr{log}}$ over $\widetilde{S}^{\circledS \mr{log}}$ and a strict super\'{e}tale morphism  $Y^{\circledS \mr{log}} \migi X^{\circledS \mr{log}}$  over $S^{\circledS \mr{log}}$.
 Then, by Proposition \ref{p0607}, there exists a  log supersmooth lifting of $Y^{\circledS \mr{log}}$ over $\widetilde{S}^{\circledS \mr{log}}$ which is uniquely determined up to isomorphism.
 We denote this log supersmooth lifting by 
 \begin{align} \label{E64}
\widetilde{\mbX} |_{Y^{\circledS \mr{log}}} :=  (\widetilde{X}^{\circledS \mr{log}} |_{Y^{\circledS \mr{log}}}, \widetilde{f}^{\circledS \mr{log}} |_{Y^{\circledS \mr{log}}}, i_{X}^{\circledS \mr{log}} |_{Y^{\circledS \mr{log}}}).
 \end{align}
  \end{rema}

\vspace{3mm}
\bco \label{c0404} \leavevmode\\
 \ \ \
Let us keep the notation in Definition \ref{D0719}. Suppose further that   $\mcJ$ is  square  nilpotent 
(hence $\mcJ$ may be thought of as an $\mcO_{S^\circledS}$-supermodule.)
Also, write 
\begin{align}
\mcF := \mcH om_{\mcO_{X^\circledS}} (\Omega_{X^{\circledS \mr{log}}/S^{\circledS \mr{log}}}, \mcJ \mcO_{X^\circledS})_{b}.
\end{align}
(Note that, since $X^{\circledS}$ is superflat over $S^\circledS$,  we have $\mcF \cong (\mcT_{X^{\circledS \mr{log}}/S^{\circledS \mr{log}}}\otimes f^{\circledS *}(\mcJ))_b$.)

\begin{itemize}
\item[(i)]
Suppose that we are given a log supersmooth lifting  $\widetilde{\mbX} :=(\widetilde{X}^{\circledS \mr{log}}, \widetilde{f}^{\circledS \mr{log}}, i_X^{\circledS \mr{log}})$  of $X^{\circledS \mr{log}}$ over $\widetilde{S}^{\circledS \mr{log}}$.
Then, the group of automorphisms of $\widetilde{\mbX}$
  is canonically isomorphic to  $\Gamma (X_b, \mcF)$.
\item[(ii)]
Suppose that we are given  two log smooth liftings $\widetilde{\mbX}_l := (\widetilde{X}_l^{\circledS \mr{log}}, \widetilde{f}_l^{\circledS \mr{log}}, i_{\mbX_l}^{\circledS \mr{log}})$ ($l =1, 2$)  of $X^{\circledS \mr{log}}$ over $\widetilde{S}^{\circledS \mr{log}}$.
Then, there exists a strict super\'{e}tale covering $Y^{\circledS \mr{log}} \migi X^{\circledS \mr{log}}$
such that $\widetilde{\mbX}_1 |_{Y^{\circledS \mr{log}}} \isom  \widetilde{\mbX}_2  |_{Y^{\circledS \mr{log}}}$.
In particular, if 
there exists a log supersmooth lifting of $X^{\circledS\mr{log}}$ over $\widetilde{S}^{\circledS \mr{log}}$, 
then  the set of isomorphism classes of log supersmooth liftings of $X^{\circledS \mr{log}}$  over $\widetilde{S}^{\circledS \mr{log}}$
forms canonically an affine space modeled on 
$H^1 (X_b, \mcF)$.
\item[(iii)]
A log supersmooth lifting of $X^{\circledS \mr{log}}$ over $\widetilde{S}^{\circledS \mr{log}}$ exists if  $H^2 (X_b, \mcF) =0$.
\end{itemize}

     \eco
\begin{proof}
Assertions (i),   (ii),  and (iii) follow from  Proposition \ref{p0404}, (ii) together with a routine argument in the theory of the classical (log) smoothness.
\end{proof}

\vspace{10mm}
\section{Stable log twisted  $\text{SUSY}_1$ curves} \vspace{3mm}

In this section, we shall 
consider, by means of various notions defined previously, supersymmetric analogues of a pointed  log twisted curve (with a canonical logarithmic structure).
We first recall (in \S\,\ref{S31}) the definition of a  twisted curve  and prove the Riemann-Roch theorem  for twisted curves  (cf. Theorem \ref{p04033}), which will be used in, e.g., computing the superdimension of the relevant  moduli introduced later.
Then,   log twisted $(1|1)$-curves (cf. Definition \ref{D0219f} (i))
are defined and characterized by local models, which are the fiber products of a (locally defined) log twisted curve and the affine superspace of superdimension $0|1$.
Moreover, by introducing  a logarithmic and twisted analogue of  superconformal structure,   we obtain the notion of a (pointed) log twisted $\text{SUSY}_1$ curve (cf. Definition \ref{D02}) which are central objects of the present paper.
 As shown in Corollary \ref{c01033}, it is  a basic property that suitable (with respect to log supersmooth  deformation) local models of  a (pointed) log twisted $\text{SUSY}_1$ curve may be chosen.
Finally, we introduce  the fibered category ${^{\S_1} \overline{\mfM}}_{g,r, \lambda}^{\circledS \mr{log}}$   (cf. (\ref{EE12})) classifying stable log twisted $\text{SUSY}_1$ curve of prescribed type $(g,r, \lambda)$ (cf. Definition \ref{D03}).
\vspace{5mm}

\subsection{The Riemann-Roch theorem  for twisted curves} \label{S31}
\leavevmode\\
\vspace{-5mm}

In this section, let us review the notion of a twisted curve and consider the Riemann-Roch theorem for twisted curves.
Here, recall that the {\it tameness} condition on a Deligne-Mumford stack $Z$ means that for every geometric point $q : \mr{Spec} (k) \migi Z$ the group $\mr{Aut} (q)$ of its stabilizers  has order prime to the characteristic of the algebraically closed field $k$.)

\vspace{3mm}
\bde \label{d0112}
\leavevmode\\
  \ \ \ Let $\underline{S}$ be a  scheme.
\begin{itemize}
\item[(i)]
A {\bf local twisted curve} over $\underline{S}$ is a flat morphism $\underline{f} : \underline{U} \migi \underline{S}$ of tame Deligne-Mumford stacks satisfying the following three conditions (i-1)-(i-3):
\begin{itemize}
\vspace{1mm}
\item[(i-1)]
The geometric  fibers of $\underline{f}$ are purely $1$-dimensional  and, \'{e}tale locally on $\underline{U}$,  isomorphic to  nodal curves; 
\vspace{1mm}
\item[(i-2)]
The smooth locus $\underline{U}^\mr{sm}$ of $\underline{U}$ (over $\underline{S}$) is an algebraic space;
\vspace{1mm}
\item[(i-3)]
For each node $q$ of a geometric  fiber of $\underline{f}$, there exists a commutative diagram:
\begin{align} \label{e012}
\xymatrix{  V \ar[r]^{\hspace{-3mm}c} \ar[dr]_d& [T/\mu_{l'}] \ar[r]^b& R \ar[d]^a \\
& \underline{U} \ar[r]^{\underline{f}} & \underline{S},
}
\end{align}
where
\begin{itemize}
\item[$\bullet$]
$R = \mr{Spec} (A)$ for some commutative ring $A$ and $a$ denotes  an \'{e}tale  neighborhood of $\underline{f} (q) \in \underline{S}$;
\item[$\bullet$]
$T = \mr{Spec}(A[z, w]/(zw -t))$ for some $t \in A$; 
\item[$\bullet$]
$[T/\mu_{l'}]$ denotes the quotient stack of $T$ by $\mu_{l'}$, where  $l'$ is a positive integer,  $\mu_{l'}$ denotes the group scheme over $R$ of $l'$-th roots of unity, 
and the action of $\mu_{l'}$ is given by $(z,w) \mapsto (\xi \cdot z, \xi^{-1} \cdot  w)$ for any $\xi \in \mu_{l'}$.
\item[$\bullet$]
 $b$ denotes  the natural projection, $c$ denotes  an \'{e}tale morphism 
and $d$ is an  \'{e}tale neighborhood of $q$.
\end{itemize}
\end{itemize}
\item[(ii)]
Let $g$ be a nonnegative integer.
A {\bf twisted curve (of genus $g$)} over $\underline{S}$ is a local twisted curve over $\underline{S}$ which is proper and whose coarse moduli space becomes a semistable  curve (of genus $g$) over $\underline{S}$ (cf.  ~\cite{Chi1},  Definition 2.4.1). 
\end{itemize}
\ede
\vspace{3mm}

We prove the  following assertion, which   is   the Riemann-Roch theorem for line bundles on a twisted curve.

\vspace{3mm}
\bt \label{p04033} \leavevmode\\
 \ \ \
 Let $X$ be a twisted curve of genus $g$ ($\geq 0$) over an algebraically closed field $k$,  and let $\mcL$ be a line bundle on $X$ of total degree $m$.   (Here, it follows from ~\cite{Chi1}, Proposition 2.5.6,  that $m$ is necessarily an integer.)
 We shall write 
 \begin{align}
  \chi (X, \mcL) := \mr{dim}_k(H^0(X, \mcL)) - \mr{dim}_k (H^1 (X, \mcL))).
 \end{align}
 Then, we have that   $H^2 (X, \mcL) = 0$ and  $ \chi (X, \mcL)  = \mr{deg} (\mcL) - g +1$.
   \et
\begin{proof}
Write  $|X|$ for  the coarse moduli space of $X$ and  $\pi : X \migi |X |$ for the projection.
First,  we shall compare $\chi (X, \mcL)$ with $\chi (|X|, \pi_*(\mcL))$ (i.e., the Eular characteristic of $\pi_*(\mcL)$ in the classical sense).
Denote by $e_1, \cdots, e_s : \mr{Spec} (k) \migi |X|$   the set of  nodes of $|X|$.
For each $i \in \{ 1, \cdots, s \}$, 
the fiber product $X \times_{|X|, e_i} \mr{Spec} (k)$ is isomorphic to the classifying stack $B (\mr{Aut} (\widetilde{e}_i))$ of the group $\mr{Aut} (\widetilde{e}_i)$ of stabilizers of the (unique) point $\widetilde{e}_i \in X(k)$  over $e_i$.
Since $\mr{Aut} (\widetilde{e}_i)$ is, by the definition of a twisted curve, isomorphic to $\mu_{l_i}$ for some $l_i \geq 1$, 
 the pull-back $\widetilde{e}_i^* (\mcL)$ of $\mcL$ may be thought of as the trivial  $k$-module $k$ (of rank one) equipped with a $\mu_{l_i}$-action.
Write 
\begin{align}
E^{\mr{nt}} := \{ i \ | \ 1 \leq i \leq s \ \text{and the $\mu_{l_i}$-action on  $\widetilde{e}_i^* (\mcL)$ is nontrivial} \}.
\end{align}
If $i \in E^\mr{nt}$, then one verifies immediately that
\begin{align} \label{E0010}
H^q (\mu_{l_i}, \widetilde{e}_i^* (\mcL)) =0 \  (q =0, 2) \ \text{and} \ H^1 (\mu_{l_i}, \widetilde{e}_i^* (\mcL)) =k.
\end{align}
On the other hand, if $i \notin E^\mr{nt}$, then
\begin{align} \label{E0011}
H^q (\mu_{l_i}, \widetilde{e}_i^* (\mcL)) =0 \  (q =1, 2) \ \text{and} \ H^0 (\mu_{l_i}, \widetilde{e}_i^* (\mcL)) =k.
\end{align}
Now, let us consider the Leray spectral sequence
\begin{align} \label{S001}
E_2^{p, q} := H^p (|X|, \mbR^q \pi_* (\mcL)) \Longrightarrow H^{p+q} (X, \mcL)
\end{align}
associated with $\mcL$.
Since $\mbR^q \pi_*(-)$ ($q \geq 1$)  vanishes  on the smooth locus $|X|^{\mr{sm}}$ of $|X|$,
we have  that 
\begin{align} \label{E0005}
H^j (|X|, \mbR^i \pi_*(-)) = 0 \ \ \text{unless $j =0$ or $(i, j) = (0, 1)$}.
\end{align}
Moreover, $\mbR^q \pi_* (\mcL)$  is isomorphic to $\bigoplus_{i=1}^s e_{i*} (H^q (\mu_i, \widetilde{e}_i^* (\mcL)))$.
Hence, it follows  from (\ref{E0010}), (\ref{E0011}), and  (\ref{E0005}) that    $H^2 (X, \mcL) = 0$  (which completes the proof of the former equality) and 
 \begin{align} \label{E0006}
& \  \mr{dim}_k (H^1 (X, \mcL)) - \mr{dim}_k (H^1 (|X|, \pi_*(\mcL)))  \\
  =& \   \Sigma_{i\in E^\mr{nt}} \mr{dim}_k (H^1 (\mu_i, \widetilde{e}_i^* (\mcL)))\notag  \\
   =&  \ \sharp E^\mr{nt}.  \notag
 \end{align}
Since  the equality $H^0 (X, \mcL) = H^0 (|X|, \pi_* (\mcL))$ is evidently verified, 
 we have 
\begin{align}
\label{E0004}
\chi (X, \mcL) - \chi (|X|, \pi_*(\mcL)) =  \sharp E^\mr{nt}.
\end{align}

Next, we shall compare the total  degree of $\mcL$ and $\pi^*(\pi_*(\mcL))$.
For each $i \in \{ 1, \cdots, s \}$, 
the formal neighborhood $\widetilde{T}_i$ of $X$ at  $\widetilde{e}_i$ is isomorphic to 
the quotient  stack $[\mr{Spec}(R)/\mu_{l_i}]$, where  $R :=  k[[z, w]]/(zw)$ and  the $\mu_{l_i}$-action on $\mr{Spec}(R)$ is given by  $(z, w) \mapsto (\xi \cdot z, \xi^{-1} \cdot w)$ for any $\xi \in \mu_{l_i}$.
In particular, if $T_i$ denotes  the formal neighborhood  of $|X|$ at $e_i$, 
then we have $T_i \cong \mr{Spec} (R^{l_i})$, where  $R^{l_i} := k[[z^{l_i}, w^{l_i}]]/(z^{l_i}y^{l_i})$, and the  morphism $\widetilde{T}_i \migi T_i$ induced by $\pi$ is given by the natural  inclusion $R^{l_i} \migiincl R$.
A choice of   trivialization $\mcL |_{\mr{Spec}(R)} \isom \mcO_{\mr{Spec}(R)}$ allows us to identify the total space of the line bundle $\mcL |_{\widetilde{T}_i}$ with the quotient stack  $[(\widetilde{T}_i \times \mr{Spec} (k[t]))/\mu_{l_i}]$, where the $\mu_{l_i}$-action is given by $(z, w, t) \mapsto (\xi \cdot z, \xi^{-1} \cdot w, \xi^{m_{i}}\cdot t)$ for some integer $m_i$ with $0 \leq m_i \leq l$. 
Then, $\mcO_{T_i}$-module  $\pi_*(\mcL) |_{T_i}$ corresponds to  the  ideal $(z^{l - m_{i}}, w^{m_i}) \subseteq R^{l_i}$.
The restriction to $\widetilde{T}_i$ of the natural morphism $\pi^*(\pi_*(\mcL)) \migi \mcL$
may be identified with the natural  inclusion  of the ideal $(z^{l - m_{i}}, w^{m_i}) \subseteq R$.
If $i \notin E^\mr{nt}$ (resp., $i \in E^\mr{nt}$), then the length of $\mcL/ \pi^*(\pi_*(\mcL))$ at $\widetilde{e}_i$ is $0$ (resp.,  $\frac{1}{l} \cdot \mr{length} (R/(z^{l - m_{i}}, w^{m_i})) = 1$).
Since  $\pi^*(\pi_*(\mcL)) \migi \mcL$  is injective  and its cokernel is only supported at $\bigcup_{i=1}^s \mr{Im} (\widetilde{e}_i)$, 
we have
\begin{align} \label{E0003}
\mr{deg} (\mcL) = \mr{deg} (\pi^* (\pi_*(\mcL))) + \Sigma_{i \in E^\mr{nt}} 1 = \mr{deg} (\pi_*(\mcL)) + \sharp E^\mr{nt}.
\end{align}
By combining (\ref{E0004}) and (\ref{E0003}), we have the  equality $\chi (X, \mcL) = \mr{deg} (\mcL) -g+1$, as desired.
\end{proof}

\vspace{5mm}

\subsection{Log  twisted $(1|1)$-curves} \label{S32}
\leavevmode\\
\vspace{-4mm}

Let  $\underline{S}$ be a scheme and 
$\underline{f} : \underline{U} \migi \underline{S}$ a local  twisted curve over $\underline{S}$.
According to (the proof of) ~\cite{O1}, Theorem 3.6,  
there exist canoincally log structures 
\begin{align}
\alpha_{\underline{U}}^{\underline{f}} : \mcM_{\underline{U}} \migi \mcO_{\underline{U}} \ \ \text{and} \ \ \alpha^{\underline{f}}_{\underline{S}} : \mcM_{\underline{S}}\migi \mcO_{\underline{S}}
\end{align}
  on $\underline{U}$ and $\underline{S}$ respectively (where we denote the resulting log stacks by 
  $\underline{U}^{\underline{f}\text{-}\mr{log}}$ and $\underline{S}^{\underline{f}\text{-}\mr{log}}$ respectively), and moreover, a special morphism 
  \begin{align}
  \underline{f}^{\underline{f}\text{-}\mr{log}} := (\underline{f}, \underline{f}^\flat : \underline{f}^{-1}(\mcM_{\underline{S}})\migi \mcM_{\underline{U}}) :\underline{U}^{\underline{f}\text{-}\mr{log}}  \migi \underline{S}^{\underline{f}\text{-}\mr{log}}
  \end{align}
   (cf. ~\cite{O1}, Theorem 3.5 for the definition of ``{\it special}") extending $\underline{f}$.
The data $(\alpha_{\underline{U}}^{\underline{f}}, \alpha^{\underline{f}}_{\underline{S}},  \underline{f}^{\underline{f}\text{-}\mr{log}})$ is uniquely determined  up to unique isomorphism.

\vspace{3mm}
\bde\label{D0219h}\leavevmode\\
\ \ \
Let $\underline{S}$ be a scheme and $\alpha_{\underline{S}}^1 : \mcM^1_{\underline{S}} \migi \mcO_{\underline{S}}$, $\alpha^2_{\underline{S}} : \mcM_{\underline{S}}^2 \migi \mcO_{\underline{S}}$ two log structures on $\underline{S}$.
We shall say that a morphism $(\underline{S}, \alpha^1_{\underline{S}}) \migi  (\underline{S}, \alpha^2_{\underline{S}})$ of log schemes is {\bf log-like} (over $\underline{S}$) if its  underlying endomorphism of $\underline{S}$ coincides with the identity morphism. 
\ede

\vspace{3mm}
\bde\label{D0219g}\leavevmode\\
\ \ \
Let $\underline{S}^\mr{log}$ be an fs  log scheme.
A {\bf log local twisted curve over $\underline{S}^\mr{log}$} is a morphism $\underline{f}^\mr{log} : \underline{U}^\mr{log} \migi \underline{S}^\mr{log}$ of log  stacks 
satisfying  the following two  conditions:
\begin{itemize}
\item[(i)]
The underlying morphism $\underline{f} : \underline{U} \migi \underline{S}$ is a local twisted curve over $\underline{S}$;
\item[(ii)]
There exist a log-like morphisms $\underline{S}^\mr{log} \migi \underline{S}^{\underline{f}\text{-}\mr{log}}$ and $\underline{U}^\mr{log} \migi \underline{U}^{\underline{f}\text{-}\mr{log}}$ over $\underline{S}$ and $\underline{U}$ respectively which make
the square diagram
\begin{align}
\xymatrix{
\underline{U}^\mr{log} \ar[r] \ar[d]_{\underline{f}^\mr{log}} & \underline{U}^{\underline{f}\text{-}\mr{log}} \ar[d]^{  \underline{f}^{\underline{f}\text{-}\mr{log}}} \\
\underline{S}^\mr{log} \ar[r] & \underline{S}^{\underline{f}\text{-}\mr{log}}
}
\end{align} 
commute and cartesian.
\end{itemize}
\ede
\vspace{3mm}
Now, let us fix an fs log  superscheme  $S^{\circledS \mr{log}}$.

\vspace{3mm}
\bde\label{D0219f}\leavevmode\\
\vspace{-5mm}
\begin{itemize}
\item[(i)]
A {\bf log twisted $(1|1)$-curve} over $S^{\circledS \mr{log}}$ is a log superstack $X^{\circledS \mr{log}}$ 
over $S^{\circledS \mr{log}}$  such that $X^\circledS/S^\circledS$ is proper and, 
 for each geometric point $q$ of $X_b$, 
there exists a $(1|1)$-chart $(Y^{\circledS \mr{log}}, U^\mr{log}, \eta^{\circledS \mr{log}})$ on $X^{\circledS \mr{log}}/S^{\circledS \mr{log}}$ around $q$ such that $U^\mr{log}$ is a log local twisted curve over $S_b^\mr{log}$.
We shall refer to such a $(1|1)$-chart $(Y^{\circledS \mr{log}}, U^\mr{log}, \eta^{\circledS \mr{log}})$ as a {\bf log twisted $(1|1)$-chart} (around $q$) on $X^{\circledS \mr{log}}/S^{\circledS \mr{log}}$.
\item[(ii)]
Let $X^{\circledS \mr{log}}/S^{\circledS \mr{log}}$ be a log twisted $(1|1)$-curve over $S^{\circledS \mr{log}}$.
Then, the induced   stack $X_t$  is  a twisted curve over $S_t$.
We shall say that 
 $X^{\circledS \mr{log}}/S^{\circledS \mr{log}}$ is {\bf of genus $g$} if $X_t/S_t$ is of genus $g$ in the sense of ~\cite{Chi1},  Definition 2.4.1. 
\end{itemize}
\ede

\vspace{5mm}

\subsection{Pointed log twisted $(1|1)$-curves} \label{S33}
\leavevmode\\
\vspace{-5mm}


In the rest of the present paper, let us fix a pair  of nonnegative integers $(g,r)$ satisfying that $2g-2+r >0$.

\vspace{3mm}
\bde
\label{D010}\leavevmode\\
\ \ \
An {\bf $r$-pointed log   twisted $(1|1)$-curve of genus $g$} over $S^{\circledS \mr{log}}$ is a collection of data
\begin{equation} \label{e015}
\mfX^{\circledS \bigstar} :=\big(X^{\circledS \mr{log}}/S^{\circledS \mr{log}}, \{ [\sigma^\circledS_i]
 \}_{i =1}^r \big),
\end{equation}
where
\begin{itemize}
\item[$\bullet$]
$X^{\circledS \mr{log}}$ denotes  a log twisted $(1|1)$-curve  of genus $g$ over $S^{\circledS \mr{log}}$;
\vspace{1mm}
\item[$\bullet$]
$ [\sigma^\circledS_i]$ (for  each $i =1, \cdots, r $) denotes  a closed subsuperscheme of $X^\circledS$ over $S^\circledS$ represented by a closed immersion $\sigma_i^\circledS : \mbA^{0|1}_{S^\circledS}  \migi X^\circledS$ over $S^\circledS$,
 \end{itemize}
 satisfying the following conditions:
 \begin{itemize}
 \vspace{2mm}
 \item[(i)]
$\mbA^{0|1}_{S^\circledS} \times_{\sigma^\circledS_i, X^\circledS, \sigma^\circledS_j} \mbA^{0|1}_{S^\circledS} = \emptyset$
  for any pair $(i, j)$ with $i \neq j$;
\vspace{2mm}
\item[(ii)]
The smooth locus $X_t^{\mr{sm}}$ of $X_t$ (over $S_t$) may be represented  by a  scheme over $S_t$ and 
the image $\mr{Im} ( (\sigma_i)_t)$ of each $(\sigma_i)_t$ lies  in $X_t^{\mr{sm}}$.
\vspace{1mm}
  \end{itemize}
Let 
$\mfX^{\circledS \bigstar} := (X^{\circledS \mr{log}}/S^{\circledS \mr{log}}, \{[\sigma_i^\circledS] \}_{i=1}^r)$  
be an $r$-pointed log twisted $(1|1)$-curve  of genus $g$  over $S^{\circledS \mr{log}}$.
Then, 
the collection of data 
\begin{align} \label{e016}
\mfX_t^\bigstar := (X_t / S_t, \{ [ (\sigma_i)_t ] \}_{i=1}^r)
\end{align}
forms an $r$-pointed twisted curve of genus $g$ over $S_t$ (in the sense of ~\cite{AV1}, Definition 4.1.2);  we shall refer to it as  the {\bf underlying (pointed) twisted curve} of $\mfX^{\circledS \bigstar}$.
 \ede

\vspace{3mm}
\bpr \label{p0468} \leavevmode\\
 \ \ \ 
 Let $\mfX^{\circledS \bigstar} := (X^{\circledS \mr{log}}/S^{\circledS \mr{log}}, \{ [\sigma_i^\circledS ]\}_{i=1}^r)$ be an $r$-pointed twisted $(1|1)$-curve of genus $g$ over $S^{\circledS \mr{log}}$.
(We shall fix a representative $\sigma_i^\circledS$ of $[\sigma^\circledS_i]$ for each $i$.)
Also, we shall fix    $i \in \{ 1, \cdots, r \}$ and a geometric point $q$ of $\mr{Im}([\sigma_i])$.
Then, there exists 
   a collection of data
   \begin{align} \label{E001}
\mbU^\bigstar :=  (Y^{\circledS \mr{log}} \stackrel{\pi^{\circledS \mr{log}}}{\migi} X^{\circledS \mr{log}}, U^\mr{log}, \eta^{\circledS \mr{log}}, \Sigma^U, \sigma^U), 
   \end{align}
   where 
   \begin{itemize}
   \item[$\bullet$]
   the triple $(Y^{\circledS \mr{log}}, U^\mr{log}, \eta^{\circledS \mr{log}})$ is  a log twisted $(1|1)$-chart around $q$ of $X^{\circledS \mr{log}}/S^{\circledS \mr{log}}$;
\item[$\bullet$]
$\pi^{\circledS \mr{log}}$ denotes  the structure morphism of $Y^{\circledS \mr{log}}$ over $X^{\circledS \mr{log}}$;
\item[$\bullet$]
$\Sigma^U$ is an \'{e}tale  scheme  over $S_b$;
\item[$\bullet$]
$\sigma^U$ is  a closed immersion $\Sigma^U \migi U$ over  $S_b$
\end{itemize}
such that  the square  diagram
\begin{align}
\begin{CD}
  \Sigma^U \times_{S_b} \mbA_{S^\circledS}^{0|1} @> \mr{pr}_2^\circledS >>  \mbA^{0|1}_{S^\circledS}
  \\
  @V \sigma^U \times \mr{id} VV @VV \sigma_i^\circledS V
  \\
  U \times_{S_b} \mbA_{S^\circledS}^{0 |1}   @>> \pi^\circledS \circ (\eta^\circledS)^{-1} >  X^\circledS. 
\end{CD}
\end{align}
is commutative and cartesian, where  $\mr{pr}_2^\circledS$ denotes the projection to the second factor.
We shall refer to such a collection of data $\mbU^\bigstar$ as  a {\bf pointed log twisted $(1|1)$-chart}  (around $q$) on $X^{\circledS \mr{log}}/S^{\circledS \mr{log}}$.
\epr
\begin{proof}
We may suppose, without loss of generality, that $S^\circledS$ is affine.
Let us take a log twisted $(1|1)$-chart $(Y^{\circledS \mr{log}}, U^\mr{log}, \eta^{\circledS \mr{log}})$ around $q$ on  $X^{\circledS \mr{log}}/S^{\circledS \mr{log}}$ such that there is no nodal point in $U$.
It follows from Proposition \ref{p0607} that  there exist
an \'{e}tale  scheme $\Sigma^U$ over $S_b$
  and  a closed immersion
 $\sigma^{Y \circledS} :  \Sigma^U \times_{S_b} \mbA^{0|1}_{S^\circledS} \migi Y^\circledS$  
 which make  the diagram 
\begin{align}
\xymatrix{
 \Sigma^U \times_{S_b} \mbA^{0|1}_{S^\circledS} 
\ar[r]^{\hspace{5mm}\mr{pr}_2^\circledS} \ar[d]_{\sigma^{Y \circledS}} & \mbA^{0|1}_{S^\circledS} \ar[d]^{\sigma_i^{\circledS}} \\
Y^\circledS \ar[r]_{\pi^{\circledS \mr{log}}}& X^{\circledS}
}
\end{align}
commutate and cartesian.
Consider the composite
 \begin{align}
 \sigma_{\eta}^{Y\circledS} : \Sigma^U \times_{S_b} \mbA_{S^\circledS}^{0|1}   \stackrel{\sigma^{Y \circledS}}{\migi}   Y^\circledS  \stackrel{\eta^\circledS}{\migi} U \times_{S_b} \mbA^{0|1}_{S^\circledS}.
 \end{align}
By applying Lemma \ref{L0102} below, there exists an automorphism $u^\circledS$ of $U \times_{S_b} \mbA^{0|1}_{S^\circledS}$ over $S^\circledS$ such that $u^\circledS \circ \sigma_{\eta}^{Y\circledS} = \sigma^U \times \mr{id}_{\mbA^{0|1}_{S^\circledS}}$ for some closed immersion $\sigma^U : \Sigma^U \migi U$.
Notice that (since 
there is no nodal point in $U$) 
the log structure of $U^\mr{log}$ coincides with the  pull-back from $S_b^\mr{log}$.
Hence, $u^\circledS$ extends to an automorphism $u^{\circledS \mr{log}}$ of $U^\mr{log} \times_{S_b} \mbA^{0|1}_{S^\circledS}$ over $S^{\circledS \mr{log}}$.
Thus, the collection of data
\begin{align}
\mbU^\bigstar := (Y^{\circledS \mr{log}}, U^\mr{log}, u^{\circledS \mr{log}} \circ \eta^{\circledS \mr{log}}, \Sigma^U, \sigma^U)
\end{align}
obtained in this manner forms the desired  collection.
This completes the proof of Proposition  \ref{p0468}.
\end{proof}
\vspace{3mm}

The following lemma was used in the proof of Proposition \ref{p0468}.

\vspace{3mm}
\ble
 \label{L0102} \leavevmode\\
 \ \ \
 Suppose that $S^\circledS$ is affine.
 Let $\Sigma$ be an affine scheme over $S_b$  and $U$  a smooth affine  scheme over $S_b$.
Then, for any  closed immersion $\sigma^\circledS : \Sigma \times_{S_b} \mbA^{0|1}_{S^\circledS} \migi U \times_{S_b} \mbA^{0|1}_{S^\circledS}$ over $S^{\circledS}$,
  there exists an automorphism $u^\circledS$ of $U \times_{S_b} \mbA^{0|1}_{S^\circledS}$ over $S^\circledS$ such that  the composite $u^\circledS  \circ \sigma^\circledS : \Sigma \times_{S_b} \mbA^{0|1}_{S^\circledS} \migi U \times_{S_b} \mbA^{0|1}_{S^\circledS}$  is of the form $\sigma_0 \times \mr{id}_{\mbA^{0 |1}_{S^\circledS}}$ for some closed immersion $\sigma_0 : \Sigma \migi U$ over $S_b$.
   \ele
\begin{proof}

In the following, let us construct two morphisms  $\delta_1^\circledS : U \times_{S_b} \mbA^{0|1}_{S^\circledS} \migi  U$ and  $\delta_2^\circledS : U \times_{S_b} \mbA^{0|1}_{S^\circledS} \migi  \mbA^{0|1}_{S^\circledS}$.

First, we shall consider $\delta_2^\circledS$.
 The map $\Gamma (U, \mcO_{(U \times_{S_b} \mbA^{0|1}_{S^\circledS})_f}) \migi \Gamma (\Sigma, \mcO_{(\Sigma \times_{S_b} \mbA^{0|1}_{S^\circledS})_f})$ induced by the closed immersion $\sigma^\circledS$ is surjective.
Hence, by Proposition \ref{P0}, 
the map 
\begin{align} \label{EE44}
\mr{Hom}_{S^\circledS}(U \times_{S_b} \mbA^{0|1}_{S^\circledS},  \mbA^{0|1}_{S^\circledS}) & \migi \mr{Hom}_{S^\circledS} (\Sigma \times_{S_b} \mbA^{0|1}_{S^\circledS},  \mbA^{0|1}_{S^\circledS}) \\
h^\circledS \hspace{15mm} & \mapsto \hspace{10mm} h^\circledS \circ \sigma^\circledS \notag
\end{align}
obtained by composing with $\sigma^\circledS$  is surjective.
Then, let us   take $\delta_2^\circledS$ to be  an inverse image 
 of the projection $\Sigma \times_{S_b} \mbA^{0|1}_{S^\circledS}  \migi \mbA^{0|1}_{S^\circledS}$ via the surjection (\ref{EE44}).

Next, we shall consider $\delta_1^\circledS$.
 We shall write
\begin{align}
 \sigma_1^\circledS : \Sigma \times_{S_b} \mbA^{0|1}_{S^\circledS} \stackrel{\sigma^\circledS}{\migi} U \times_{S_b} \mbA^{0|1}_{S^\circledS} \stackrel{\mr{pr}^\circledS}{\migi} U, 
\end{align}
where the second arrow $\mr{pr}^\circledS$ denotes the projection to the first factor.
Also, write
 \begin{align}
 \widetilde{\sigma}_t^\circledS : \Sigma \times_{S_b} \mbA^{0|1}_{S^\circledS} \migi \Sigma \stackrel{\sigma_t}{\migi} U,
 \end{align}
 where the first arrow denotes  the projection $\Sigma \times_{S_b} \mbA^{0|1}_{S^\circledS} \migi \Sigma$ to the first factor.
 Moreover, denote  by $\mcJ_{\Sigma}$ and $\mcJ_U$ 
 the (square nilpotent) ideal  of $\mcO_{(\Sigma \times_{S_b} \mbA^{0|1}_{S^\circledS})_b}$ and $\mcO_{(U \times_{S_b} \mbA^{0|1}_{S^\circledS})_b}$ corresponding to the closed   immersions 
 \begin{align}
 \gamma_{\Sigma \times_{S_b} \mbA^{0|1}_{S^\circledS}} :  \Sigma \migi (\Sigma \times_{S_b} \mbA^{0|1}_{S^\circledS})_b \ \ \text{and} \  \ 
 \gamma_{U \times_{S_b} \mbA^{0|1}_{S^\circledS}} :  U \migi (U \times_{S_b} \mbA^{0|1}_{S^\circledS})_b
 \end{align}
  respectively.
Then, there exists an element $\partial \in \Gamma (\Sigma, \sigma_t^*(\mcT_{U/S_b})\otimes \mcJ_\Sigma)$ such that $(\sigma_1)_b \boxplus^\dagger \partial = (\widetilde{\sigma}_t^{\circledS})_b$
 (cf. Proposition \ref{p0404} (ii) for the definition of ``$\boxplus^\dagger$").
Since the morphism 
\begin{align} \label{FF07}
\Gamma (U, \mcT_{U/S_b} \otimes \mcJ_U) \migi \Gamma (\Sigma, \sigma_t^*(\mcT_{U/S_b})\otimes \mcJ_\Sigma)
\end{align}
 induced by $\sigma_t$
 is surjective, we obtain an inverse image $\widetilde{\partial} \in \Gamma (U, \mcT_{U/S_b}\otimes \mcJ_U) $ of $\partial$ via (\ref{FF07}).
Thus, we obtain a morphism 
\begin{align}
\delta^\circledS_1 :=  (\mr{pr}_b \boxplus^\dagger \widetilde{\partial}) \circ \beta^\circledS_{U \times_{S_b} \mbA^{0|1}_{S^\circledS}} :  U \times_{S_b} \mbA^{0|1}_{S^\circledS} \migi U
\end{align}
over $S_b$.

It follows from the definitions of $\delta_1^\circledS$ and $\delta_2^\circledS$ that
the endomorphism $u^\circledS := (\delta_1^\circledS, \delta_2^\circledS)$ of $U \times_{S_b} \mbA^{0|1}_{S^\circledS}$ turns out to be the desired  automorphism.
This completes the proof of  Lemma \ref{L0102}.
\end{proof}
\vspace{3mm}

\begin{rema} \label{r4f828} \leavevmode\\
 \ \ \ 
Let $\mfX^{\circledS \bigstar} := (X^{\circledS \mr{log}}/S^{\circledS \mr{log}}, \{ [ \sigma_i^\circledS ] \}_{i=1}^r)$ be an $r$-pointed log twisted $(1|1)$-curve of genus $g$ over $S^{\circledS \mr{log}}$, $q$ a geometric point of $X_b$, and  $\mbU := (Y^{\circledS \mr{log}}, U^\mr{log}, \eta^{\circledS \mr{log}})$  a log twisted $(1|1)$-chart around $q$ on $X^{\circledS \mr{log}}/S^{\circledS \mr{log}}$.
Suppose that we are given a strict super\'{e}tale morphism  $\pi_Y^{\circledS \mr{log}} : Y'^{\circledS \mr{log}} \migi Y^{\circledS \mr{log}}$ such that $Y'^\circledS$ is affine and  the image of   the composite  $Y'^{\circledS \mr{log}} \stackrel{\pi_Y^{\circledS \mr{log}}}{\migi} Y^{\circledS \mr{log}} \migi X^{\circledS \mr{log}}$  contains $q$.
Then, by Proposition \ref{p0607}, there exist
a strict  \'{e}tale morphism $\pi_U^{\mr{log}} : U'^{\mr{log}} \migi U^\mr{log}$   and an isomorphism
$\eta'^{\circledS \mr{log}} : Y'^{\circledS \mr{log}} \isom U'^{\mr{log}} \times_{S_b} \mbA^{0|1}_{S^\circledS}$ over $S^{\circledS \mr{log}}$  (hence $U'$ is affine) which make the square diagram
\begin{align}
\xymatrix{
Y'^{\circledS \mr{log}} \ar[r]^{\eta'^{\circledS \mr{log}}}  \ar[d]_{\pi_Y^{\circledS \mr{log}}} & U'^{\mr{log}} \times_{S_b} \mbA^{0|1}_{S^\circledS} \ar[d]^{\pi_U^{\mr{log}} \times \mr{id}}  \\
Y^{\circledS \mr{log}} \ar[r]^{\eta^{\circledS \mr{log}}} & U^\mr{log} \times_{S_b} \mbA^{0|1}_{S^\circledS} 
}
\end{align}
commute and cartesian.
(Such a pair $(U'^{\mr{log}}, \eta^{\circledS \mr{log}})$ is uniquely determined up to isomorphism).
We shall write
\begin{align}
\mbU |_{Y'^{\circledS \mr{log}}} := (Y'^{\circledS \mr{log}}, U'^{\mr{log}}, \eta'^{\circledS \mr{log}}),
\end{align}
which forms  a log twisted $(1|1)$-chart around $q$ on $X^{\circledS \mr{log}}/S^{\circledS \mr{log}}$.

Suppose  further that $q \in \mr{Im} ((\sigma_i)_b)$ and we are given 
 $\Sigma^U$ and $\sigma^U$ as in (\ref{E001}) for which the collection of data  $\mbU^\bigstar := (Y^{\circledS \mr{log}}, U^\mr{log}, \eta^{\circledS \mr{log}}, \Sigma^U, \sigma^U)$ forms 
a pointed log twisted $(1|1)$-chart   around $q$ on $X^{\circledS \mr{log}}/S^{\circledS \mr{log}}$.
Let us  write 
\begin{align}
\Sigma^{U'} := \Sigma^U \times_{\sigma^U, U, \pi_U} U' \ \ \text{and} \ \ \sigma^{U'} := \sigma^U \times \mr{id}_{U'} : \Sigma^{U'} \migi U'.
\end{align}
Then,
\begin{align} \label{E1001}
\mbU^\bigstar |_{Y'^{\circledS \mr{log}}} := (Y'^{\circledS \mr{log}}, U'^{\mr{log}}, \eta'^{\circledS \mr{log}}, \Sigma^{U'}, \pi^{U'}, \sigma^{U'})
\end{align} 
forms a  pointed log twisted $(1|1)$-chart $\mbU^\bigstar$  around $q$ of $X^{\circledS \mr{log}}/S^{\circledS \mr{log}}$.
 \end{rema}

\vspace{5mm}

\subsection{Superconformal structures} \label{S34}
\leavevmode\\
\vspace{-4mm}

Let  us fix   an $r$-pointed log  twisted $(1|1)$-curve  $\mfX^{\circledS \bigstar} :=(X^{\circledS \mr{log}}/S^{\circledS \mr{log}}, \{ [\sigma_i^{\circledS}]\}_{i =1}^r)$ of genus $g$ over $S^{\circledS \mr{log}}$.
We shall construct  a new log structure  on $X^\circledS$ as follows.
The ideal sheaf $\mcI_i \subseteq \mcO_{X_b}$ ($i =1, \cdots, r$) defining the closed immersion $(\sigma_i)_b$  is, by Proposition \ref{p0468},  an invertible sheaf.
 As explained in ~\cite{KATO}, Complement 1, it corresponds to a log structure $\alpha_{X_b}^{\sigma_i} : \mcM_{X_b}^{\sigma_i} \migi \mcO_{X_b}$.
We shall write
\begin{align}
\mcM^\bigstar_{X_b} := \mcM_{X_b} \oplus_{\mcO_{X_b}^\times} (\mcM^{\sigma_1}_{X_b} \oplus_{\mcO^\times_{X_b}} \cdots \oplus_{\mcO^\times_{X_b}} \mcM^{\sigma_r}_{X_b})
\end{align}
and 
define a log structure $\alpha^{\bigstar}_{X_b}$
  to be the amalgam
\begin{align}
\alpha^\bigstar_{X_b} := (\alpha_{X_b}, (\alpha^{\sigma_1}_{X_b}, \cdots, \alpha_{X_b}^{\sigma_r})) : \mcM^\bigstar_{X_b}  \migi \mcO_{X_b}.
\end{align}
We shall denote by
\begin{align}
X^{\circledS \bigstar\text{-}\mr{log}} := (X^{\circledS }, \alpha^\bigstar_{X_b})
\end{align}
the resulting log superstack over $S^{\circledS \mr{log}}$, which admits   a natural morphism $X^{\circledS \bigstar\text{-}\mr{log}}  \migi X^{\circledS \mr{log}}$.
If $f^{\circledS \mr{log}} : X^{\circledS \mr{log}} \migi S^{\circledS \mr{log}}$ denotes the structure morphism of $X^{\circledS \mr{log}}/S^{\circledS \mr{log}}$, then we shall write 
\begin{align}
f^{\circledS \bigstar\text{-}\mr{log}} : X^{\circledS \bigstar \text{-}\mr{log}} \migi S^{\circledS \mr{log}}
\end{align}
for the composite of  $f^{\circledS \mr{log}}$ with $X^{\circledS \bigstar \text{-}\mr{log}} \migi X^{\circledS \mr{log}}$.
One verifies that $ X^{\circledS \bigstar \text{-}\mr{log}}/S^{\circledS \mr{log}}$ is log supersmooth of relative superdimension $1|1$.
In particular, the $\mcO_{X^\circledS }$-supermodule  $\mcT_{X^{\circledS \bigstar \text{-}\mr{log}}/S^{\circledS \mr{log}}}$ (as well as $\Omega_{X^{\circledS \bigstar \text{-}\mr{log}}/S^{\circledS \mr{log}}}$) is a  supervector bundle of superrank $1|1$  (cf. Proposition \ref{p0404} (i)).


\vspace{3mm}
\bde
\label{D02}\leavevmode\\
\vspace{-5mm}
\begin{itemize}
\item[(i)]
Let $X^{\circledS \mr{log}}$ be a log supersmooth superscheme over $S^{\circledS \mr{log}}$ of relative superdimension $1|1$.
A {\bf superconformal structure} on $X^{\circledS \mr{log}}/S^{\circledS \mr{log}}$ is a  subsupervector bundle $\mcD$ of superrank $0|1$  of $\mcT_{X^{\circledS \mr{log}}/S^{\circledS \mr{log}}}$ (i.e., $\mcT_{X^{\circledS \mr{log}}/S^{\circledS \mr{log}}}/\mcD$ is a supervector bundle of superrank $1|0$) such that the {\it $\mcO_{X^\circledS}$-linear}  morphism
\begin{align} \label{e032}
(\mcD^{\otimes 2} :=) \ \mcD \otimes_{\mcO_{X^\circledS}} \mcD &\migi \mcT_{X^{\circledS \mr{log}}/S^{\circledS \mr{log}}}/\mcD \\
\partial_1 \otimes \partial_2  \hspace{5mm} &\mapsto \frac{1}{2} \cdot\overline{ [\partial_1, \partial_2]} \notag
\end{align}
(where  $\partial_1$ and $\partial_2$ are local sections of $\mcD$) is an isomorphism.
\item[(ii)]
An {\bf $r$-pointed log  twisted  $\text{SUSY}_1$ curve of genus $g$}  over $S^{\circledS \mr{log}}$ is a collection of data
\begin{equation}
{^{\S_1} \mfY}_{}^{\circledS \bigstar} :=(Y^{\circledS \mr{log}}/S^{\circledS \mr{log}}, \{[\sigma_i^\circledS] \}_{i=1}^r,  \mcD)
\end{equation}
 consisting of an   $r$-pointed  log twisted $(1|1)$-curve $(Y^{\circledS \mr{log}}/S^{\circledS \mr{log}},  \{[\sigma_i^\circledS ] \}_{i=1}^r)$  of genus $g$ over $S^{\circledS \mr{log}}$ and a superconformal structure $\mcD$ on $Y^{\circledS \bigstar \text{-}\mr{log}}/S^{\circledS \mr{log}}$.
\end{itemize}
 \ede

\vspace{3mm}
\bde
\label{D021}\leavevmode\\
\ \ \
For $j \in \{1, 2\}$,  let $S_j^{\circledS \mr{log}}$  be  an fs log superscheme
 and 
 ${^{\S_1} \mfX}^{\circledS \bigstar}_j := (f_j^{\circledS \mr{log}} : X_j^{\circledS \mr{log}} \migi S_j^{\circledS \mr{log}}, \{ [\sigma^\circledS_{j, i}]  \}_{i=1}^r, \mcD_j)$
   an  $r$-pointed log  twisted $\text{SUSY}_1$ curve of genus $g$ over $S_j^{\circledS \mr{log}}$.
\begin{itemize}
\item[(i)]
A {\bf superconformal morphism} from ${^{\S_1} \mfX}^{\circledS \bigstar}_1$ to ${^{\S_1} \mfX}^{\circledS \bigstar}_2$ is a pair 
\begin{align}
{^{\S_1} \Phi}^{\circledS \bigstar} : = (\Phi^{\circledS \mr{log}},  \phi^{\circledS \mr{log}})
\end{align}
consisting of two morphisms $\Phi^{\circledS \mr{log}}  : X_1^{\circledS \mr{log}} \migi  X_2^{\circledS \mr{log}}$, $\phi^{\circledS \mr{log}} : S_1^{\circledS \mr{log}} \migi  S_2^{\circledS \mr{log}}$ such that
\vspace{1mm}
\begin{itemize}
\item[$\bullet$]
the square   diagram
\begin{align}
\xymatrix{
X_1^{\circledS \mr{log}} \ar[r]^{\Phi^{\circledS \mr{log}}} \ar[d]_{f_1^{\circledS \mr{log}}} & X_2^{\circledS \mr{log}} \ar[d]^{f_2^{\circledS \mr{log}}}  \\
S_1^{\circledS \mr{log}} \ar[r]_{\phi^{\circledS \mr{log}}} & S_1^{\circledS \mr{log}} 
}
\end{align}
($i = 1, \cdots, r$)
 is  commutative and cartesian;
\item[$\bullet$]
$[\sigma^\circledS_{1, i}] = \Phi^{\circledS *} ([\sigma^\circledS_{2, i}])$ (for any $i \in \{ 1, \cdots, r \}$)  and 
  $\mcD_1 = \Phi^{\circledS *} (\mcD_2)$  via the  isomorphism $\mcT_{X_1^{\circledS \bigstar\text{-}\mr{log}}/S_1^{\circledS \mr{log}}}\isom  \Phi^{\circledS*}(\mcT_{X_2^{\circledS \bigstar\text{-}\mr{log}}/S_2^{\circledS \mr{log}}})$ induced by $\Phi^{\circledS \mr{log}}$.
\end{itemize}
\item[(ii)]
Suppose further that $S_1^{\circledS \mr{log}} = S_2^{\circledS \mr{log}}$ ($=: S^{\circledS \mr{log}}$).
A {\bf superconformal isomorphism over $S^{\circledS \mr{log}}$}  from 
${^{\S_1} \mfX}^{\circledS \bigstar}_1$ to ${^{\S_1} \mfX}^{\circledS \bigstar}_2$
 is a superconformal morphism ${^{\S_1} \Phi}^{\circledS \bigstar} :=(\Phi^{\circledS \mr{log}}, \phi^{\circledS \mr{log}}) :{^{\S_1} \mfX}^{\circledS \bigstar}_1 \migi {^{\S_1} \mfX}^{\circledS \bigstar}_2$ 
 such that $\phi^{\circledS \mr{log}} = \mr{id}_{S^{\circledS \mr{log}}}$ and $\Phi^{\circledS \mr{log}}$  is an  isomorphism.
\end{itemize}
 \ede
\vspace{3mm}

In the following Proposition \ref{p01033}, we discuss an explicit description of a superconformal structure (cf. e.g., ~\cite{Witten2}, Lemma 3.1,  for the case of smooth $\text{SUSY}_1$ curves over $\mbC$, i.e., super Riemann surfaces.)

\vspace{3mm}
\bpr \label{p01033} \leavevmode\\
 \ \ \
 Let $U^\mr{log}$ be a log smooth scheme over $S_b^\mr{log}$ of relative dimension $1$.
 (In particular, $Z^{\circledS \mr{log}} := U^\mr{log} \times_{S_b} \mbA^{0|1}_{S^\circledS}$ is a log supersmooth superscheme over $S^{\circledS \mr{log}}$ of relative superdimension $1|1$).
 Suppose that
 we are given an element $z \in \Gamma (U, \mcM_{U})$ such that $\Omega_{U^\mr{log}/S_b^\mr{log}} \cong \mcO_U \cdot d \mr{log} (z)$.
 Let us  regard $d \mr{log} (z)$ and  $d (\psi)$ as sections of $\Omega_{Z^{\circledS \mr{log}}/S^{\circledS \mr{log}}}$ via the projections $Z^{\circledS \mr{log}} \migi U^\mr{log}$ and $Z^{\circledS \mr{log}}  \migi \mbA^{0|1}_{S^\circledS}$ respectively;  these sections give a decomposition
 $\Omega_{Z^{\circledS \mr{log}}/S^{\circledS \mr{log}}} \cong \mcO_{Z^\circledS} \cdot d (\psi) \oplus \mcO_{Z^{\circledS}} \cdot d \mr{log} (z)$.
 In particular, we have 
 \begin{align}
 \mcT_{Z^{\circledS \mr{log}}/S^{\circledS \mr{log}}} \cong   \mcO_{Z^\circledS} \cdot \partial_\psi \oplus \mcO_{Z^{\circledS}} \cdot \partial_z,
 \end{align}
  where $\{ \partial_\psi, \partial_z \}$ is  the dual basis of $\{ d (\psi), d \mr{log} (z) \}$.
 Then, the following assertions are satisfied.
 \begin{itemize}
 \item[(i)]
 For each $a \in \Gamma (Z_b, \mcO_{Z_b}^\times)$, the subsupermodule 
 \begin{align}
 \mcD_a := \mcO_{Z^\circledS} \cdot (\partial_\psi + a \psi \cdot \partial_z)
 \end{align}
  of $\mcT_{Z^{\circledS \mr{log}}/S^{\circledS \mr{log}}}$ forms a superconformal structure on $Z^{\circledS \mr{log}}/S^{\circledS \mr{log}}$.
 Moreover, the assignment $a \mapsto \mcD_a$ determines a bijection between the set  $\Gamma (Z_b, \mcO_{Z_b}^\times)$ and the set of superconformal structres on $Z^{\circledS \mr{log}}/S^{\circledS \mr{log}}$.
\vspace{1mm}
 \item[(ii)]
Let us take two superconformal structures on $Z^{\circledS \mr{log}}/S^{\circledS \mr{log}}$ of the form $\mcD_a$,  $\mcD_b$ for some $a$, $b \in \Gamma (Z_b, \mcO_{Z_b}^\times)$.
Suppose that there exists an element $c \in \Gamma (Z_b, \mcO_{Z_b}^\times)$ such that $c^2 \cdot a = b$.
(According to Proposition \ref{p0607}, such an element $c$ exists after possibly  replacing $U$ with
its  \'{e}tale covering.)
If we write $\iota_{c}$ for the automorphism of $Z^{\circledS \mr{log}}$ over $U^\mr{log} \times_{S_b} S^\circledS$ given by assigning $\psi \mapsto c \cdot \psi$, then the isomorphism $\mcT_{Z^{\circledS \mr{log}}/S^{\circledS \mr{log}}} \isom \iota_{c}^* (\mcT_{Z^{\circledS \mr{log}}/S^{\circledS \mr{log}}})$ induces an isomorphism $\mcD_a \isom \iota_{c}^*(\mcD_b)$. 
  \end{itemize}

   \epr
\begin{proof}
First, we consider assertion (i).
For each $a \in \Gamma (Z_b, \mcO_{Z_b}^\times)$, we have
\begin{align}
& \ \ \ \  \frac{1}{2} \cdot [\partial_\psi + a \psi \cdot \partial_z, \partial_\psi + a \psi \cdot \partial_z] \\
& = \partial_\psi^2 + ((a \psi \cdot  \partial_z) \circ \partial_\psi + \partial_\psi \circ (a \psi \cdot \partial_z)) + (a \psi \cdot \partial_z)^2 \notag \\
& = 0  + a \cdot \partial_z + 0  \notag \\
& = a \cdot \partial_z. \notag 
\end{align}
Hence, (since  the sections  $\{ \partial_\psi + a \psi \cdot \partial_z, a \cdot  \partial_z \}$ generate $\mcT_{Z^{\circledS \mr{log}}/S^{\circledS \mr{log}}}$) $\mcD_a$ forms a superconformal structure on $Z^{\circledS \mr{log}}/S^{\circledS \mr{log}}$. 

Next, we shall consider the bijectivity of  the assignment $a \mapsto \mcD_a$.
Let  $\mcD$ be a superconformal structure on $Z^{\circledS \mr{log}}/S^{\circledS \mr{log}}$.
There exists an open covering $\{ Z_\gamma \}_\gamma$ of $Z_b$ such that
each restriction $\mcD |_{Z_\gamma}$ may be generated by  some $\partial_\gamma \in \Gamma (Z_\gamma, \mcD)$.
The section  $\partial_\gamma$ may be described  as  $\partial_\gamma := a_\gamma \psi \cdot \partial_z + b_\gamma \cdot \partial_\psi$ (where $a_\gamma, b_\gamma \in \Gamma (Z_\gamma, \mcO_{X_b})$).
 Then, 
 \begin{align}
 \partial^{2}_\gamma := \frac{1}{2} \cdot [ \partial_\gamma, \partial_\gamma ] = a_\gamma b_\gamma \cdot \partial_z + a_\gamma (\partial_z (b_\gamma)) \psi  \cdot \partial_\psi.
 \end{align}
 Since $\{\partial_\gamma, \partial^{2}_\gamma \}$ generates $\mcT_{Z^{\circledS \mr{log}}/S^{\circledS \mr{log}}}|_{Z_\gamma}$,   both $a_\gamma$ and $b_\gamma$  lie in $\Gamma (Z_\gamma, \mcO_{Z_b}^\times)$.
 Thus, there exists {\it uniquely} an element  $\partial'_\gamma$ (i.e., $\partial'_\gamma := b_\gamma^{-1} \cdot \partial_\gamma$)  in $\Gamma (Z_\gamma, \mcD)$ of the form $\partial_\psi + a'_\gamma \psi \cdot \partial_z$ (for some $a'_\gamma \in \Gamma (Z_\gamma, \mcO_{Z_b}^\times)$).
 In particular, $\{ \partial'_\gamma \}_\gamma$ may be glued together to an element  of $\Gamma (Z_b, \mcD)$ of the form $\partial_\phi + a' \psi \cdot \partial_z$ (for a unique $a' \in \Gamma (Z_b, \mcO_{Z_b}^\times)$).
This assignment $\mcD \mapsto a'$ determines an inverse to the assignment $a \mapsto \mcD_a$.
Consequently, $a \mapsto \mcD_a$ is bijective, as desired.
 
 Finally,  assertion (ii) follows immediately  from the definition of $\iota_c$.
\end{proof}
\vspace{3mm}

In particular, we have  the following assertion.

\vspace{3mm}
\bco \label{c01033} \leavevmode\\
 \ \ \
  Let ${^{\S_1} \mfX}^{\circledS \bigstar} := (X^{\circledS \mr{log}}/S^{\circledS \mr{log}}, \{ [\sigma_i^\circledS] \}_{I-1}^r, \mcD)$ be an $r$-pointed  log twisted  $\text{SUSY}_1$ curve of genus  $g$ over $S^{\circledS \mr{log}}$.
  Then, there exists a 
collection of data
\begin{align} \label{E50}
\{  (Y_\gamma^{\circledS \mr{log}} \stackrel{\pi_\gamma^{\circledS \mr{log}}}{\migi}X^{\circledS \mr{log}}, U_\gamma^\mr{log}, \eta_\gamma^{\circledS \mr{log}}, z_\gamma)\}_\gamma, 
\end{align}
where
\begin{itemize}
\item[$\bullet$]
$\{ (Y_\gamma^{\circledS \mr{log}} \stackrel{\pi_\gamma^{\circledS \mr{log}}}{\migi}X^{\circledS \mr{log}}, U_\gamma^\mr{log}, \eta_\gamma^{\circledS \mr{log}}) \}_\gamma$ is a collection of log twisted $(1|1)$-chart on $X^{\circledS \mr{log}}/S^{\circledS \mr{log}}$ such that $\coprod_\gamma Y^{\circledS \mr{log}}_\gamma \migi X^{\circledS \mr{log}}$ is a strict super\'{e}tale covering of $X^{\circledS \mr{log}}$;
\item[$\bullet$]
Each $z_\gamma$ is an element of $\Gamma (U_\gamma, \mcM_{U_\gamma})$ such that $d \mr{log} (z_\gamma)$ generates  $\Omega_{U_\gamma^\mr{log}/S_b^\mr{log}}$ and 
the superconformal structure $\mcD |_{Y_\gamma^{\circledS \mr{log}}}$  on $Y^{\circledS \mr{log}}_\gamma/S^{\circledS \mr{log}}$  obtained by restricting  $\mcD$ to $Y_\gamma^{\circledS \mr{log}}$ coincides with 
\begin{align}
\mcO_{U_\gamma^\mr{log} \times_{S_b}\mbA_{S^\circledS}} \cdot (\partial_\psi + \psi \cdot \partial_{z_\gamma}) \subseteq \mcT_{U_\gamma^\mr{log} \times_{S_b}\mbA_{S^\circledS}^{0|1}/S^{\circledS \mr{log}}}
\end{align}
 (where $\{ \partial_\psi, \partial_{z_\gamma} \}$ is the dual basis of $\{ d (\psi), d \mr{log} (z_\gamma) \}$) via the isomorphism $\mcT_{Z^{\circledS \mr{log}}/S^{\circledS \mr{log}}} \isom (\eta_\gamma^{\circledS})^*(\mcT_{U_\gamma^\mr{log} \times_{S_b}\mbA_{S^\circledS}^{0|1}/S^{\circledS \mr{log}}})$ induced by $\eta_\gamma^\circledS$.
\end{itemize}
   \eco

\vspace{5mm}

\subsection{Kodaira-Spencer morphisms} \label{S35}
\leavevmode\\
\vspace{-4mm}


Let ${^{\S_1} \mfX}^{\circledS \bigstar} := (f^{\circledS \mr{log}} : X^{\circledS \mr{log}} \migi S^{\circledS \mr{log}}, \{ [\sigma_i^\circledS] \}_{i=1}^r, \mcD)$ be an $r$-pointed log twisted $\text{SUSY}_1$ curve of genus $g$ over $S^{\circledS \mr{log}}$.
Let us  define an $f_b^{-1}(\mcO_{S^\circledS})$-subsupermodule  
$\mcT_{X^{\circledS \bigstar\text{-}\mr{log}}/S_0}^\mcD$
 (resp.,  $\mcT_{X^{\circledS \bigstar \text{-}\mr{log}}/S^{\circledS \mr{log}}}^\mcD$) 
  of $\mcT_{X^{\circledS \bigstar\text{-}\mr{log}}/S_0}$ (resp.,  $\mcT_{X^{\circledS \bigstar \text{-}\mr{log}}/S^{\circledS \mr{log}}}$)  to be
\begin{align}
\mcT_{X^{\circledS \bigstar\text{-}\mr{log}}/S_0}^\mcD  & := \{ \partial \in \mcT_{X^{\bigstar\text{-}\mr{log}}/S_0} \ | \ [\partial, \mcD] \subseteq \mcD \} \\
(\text{resp.,} \  \mcT_{X^{\circledS \bigstar \text{-}\mr{log}}/S^{\circledS \mr{log}}}^\mcD & := \mcT_{X^{\circledS \bigstar \text{-}\mr{log}}/S^{\circledS \mr{log}}}  \cap \mcT_{X^{\circledS \bigstar\text{-}\mr{log}}/S_0}^\mcD). \notag
\end{align}  
Since $X^{\circledS \bigstar \text{-}\mr{log}}$ is log supersmooth over $S^{\circledS \mr{log}}$,
the dual of the sequence (\ref{E20}) gives rise to  a short exact 
sequence of $\mcO_{X^\circledS}$-supermodules:
\begin{align} \label{S0011}
0 \migi  \mcT_{X^{\circledS \bigstar \text{-}\mr{log}}/S^{\circledS \mr{log}}}^\mcD  \migi \mcT_{X^{\circledS \bigstar\text{-}\mr{log}}/S_0}^\mcD   \migi f^{\circledS *}(\mcT_{S^{\circledS \mr{log}}/S_0}) \migi 0.
\end{align}
(Here, the pulled-back $\mcO_{S^\circledS}$-supermodule  $f^{\circledS *}(-)$ via $f^\circledS$ defined preceding Definition \ref{d32}  may be also defined in our situation, i.e.,  $X^\circledS$ is a superstack.)
The higher direct image $\mbR^1 f_{b*}(\mcT_{X^{\circledS \bigstar \text{-}\mr{log}}/S^{\circledS \mr{log}}}^\mcD)$ admits naturally a structure of $\mcO_{S^\circledS}$-supermodule.
Denote by $\mbR^1 f_{*}^\circledS(\mcT_{X^{\circledS \bigstar \text{-}\mr{log}}/S^{\circledS \mr{log}}}^\mcD)$ the resulting $\mcO_{S^\circledS}$-supermodule.
The connecting homomorphism of (\ref{S0011}) yields an $\mcO_{S^\circledS}$-linear morphism
\begin{align}
\mcK \mcS({^{\S_1} \mfX}^{\circledS \bigstar}) : \mcT_{S^{\circledS \mr{log}}/S_0}  \migi \mbR^1 f_{*}^\circledS (\mcT_{X^{\circledS \bigstar \text{-}\mr{log}}/S^{\circledS \mr{log}}}^\mcD),
\end{align}
which is referred to as the {\bf Kodaira-Spencer morphism} of ${^{\S_1} \mfX}^{\circledS \bigstar}$.

The following proposition will be used in the discussion in Remark  \ref{r4} and Proposition \ref{P001}.

\vspace{3mm}
\bpr
\label{p010} \leavevmode\\
 \ \ \
There exists a canonical  $f_b^{-1}(\mcO_{S^\circledS})$-linear isomorphism $\mcT_{X^{\circledS \bigstar \text{-}\mr{log}}/S^{\circledS \mr{log}}}^\mcD  \isom \mcD^{\otimes 2}$.
   \epr
\begin{proof}
The assertion follows from an argument similar to the argument in the proof of ~\cite{LR}, Lemma 2.1,  together with Corollary \ref{c01033} of the present paper.
\end{proof}

\vspace{5mm}

\subsection{Stable log twisted $\text{SUSY}_1$ curves} \label{S36}
\leavevmode\\
\vspace{-4mm}

Let $\lambda$ be 
 an even positive integer  invertible in $S_0$.
Let us recall  from ~\cite{Chi1},  Definition 4.1.3 and Remark 4.2.6,  the notion of a $\lambda$-stable twisted curve.
We shall write 
\begin{align}
{^\text{tw} \overline{\mfM}}_{g,r, \lambda}
\end{align}
for the  moduli stack classifying  $r$-pointed  $\lambda$-stable twisted curves over $S_0$ of genus $g$.
It is a geometrically connected,  proper,  and smooth Deligne-Mumford stack over $S_0$ of relative dimension $3g-3+r$  (cf. ~\cite{Chi1}, Corollary 4.2.8).
Denote by $(\mfC, \{ [\sigma_{\mfC, i}] \}_{i=1}^r)$ the tautological $r$-pointed $\lambda$-stable twisted  curve over ${^\text{tw} \overline{\mfM}}_{g,r, \lambda}$.
Both ${^\text{tw} \overline{\mfM}}_{g,r, \lambda}$ and $\mfC$
 admit canonically  log structures (cf. ~\cite{O1}, Theorem 1.9).
If  ${^\text{tw} \overline{\mfM}}_{g,r, \lambda}^{\mr{log}}$ and $\mfC^\mr{log}$ denote  the resulting log stacks, then the structure morphism of $\mfC$ over ${^\text{tw} \overline{\mfM}}_{g,r, \lambda}$ extends to
a log smooth morphism $\mfC^\mr{log} \migi {^\text{tw} \overline{\mfM}}_{g,r, \lambda}^{\mr{log}}$. 

Let $s^\mr{log} : \underline{S}^\mr{log} \migi {^\text{tw} \overline{\mfM}}_{g,r, \lambda}^{\mr{log}}$ be a morphism whose underlying morphism of stacks classifies an $r$-pointed twisted curve $\underline{\mfX}^\bigstar :=  (\underline{X}/\underline{S}, \{ [\underline{\sigma}_i]\}_{i=1}^r)$ of genus $g$.
Then, by equipping $\underline{X}$ with   the log structure pulled-back from $\mfC^\mr{log} \times_{{^\text{tw} \overline{\mfM}}_{g,r, \lambda}^{\mr{log}}, s^\mr{log}}  \underline{S}^\mr{log} $ via the isomorphism $\underline{X} \isom  \mfC \times_{{^\text{tw} \overline{\mfM}}_{g,r, \lambda}, s} \underline{S}$ induced by $s$, we have a log stack 
\begin{align} \label{E1111}
\underline{X}^{\bigstar \text{-} \mr{log}}
\end{align}
 together with  a log smooth morphism $\underline{X}^{\bigstar \text{-} \mr{log}} \migi \underline{S}^\mr{log}$.

Moreover, let us  write $\mfM_{g,r}$ for  the moduli stack classifying $r$-pointed  proper {\it smooth} curves over $S_0$ of genus $g$.
 By the natural inclusion $\mfM_{g,r} \migiincl {^\text{tw} \overline{\mfM}}_{g,r, \lambda}^{\mr{log}}$, we may regard  $\mfM_{g,r}$  as a dense open substack of ${^\text{tw} \overline{\mfM}}_{g,r, \lambda}^{\mr{log}}$.
Also, this  open locus of ${^\text{tw} \overline{\mfM}}_{g,r, \lambda}^{\mr{log}}$
coincides with 
  the locus  in which the log structure of  ${^\text{tw} \overline{\mfM}}_{g,r, \lambda}^{\mr{log}}$ becomes  trivial.

\vspace{3mm}
\bde
\label{D03}\leavevmode\\
\ \ \
A {\bf stable log twisted  $\text{SUSY}_1$  curve of type $(g,r, \lambda)$}
  over $S^{\circledS \mr{log}}$ is 
an $r$-pointed  log twisted $\text{SUSY}_1$  curve of genus $g$ over $S^{\circledS \mr{log}}$ whose underlying pointed   twisted  curve is $\lambda$-stable.
 \ede
\vspace{3mm}


Let  $(g, r, \lambda)$ be 
 a  triple of nonnegative integers  satisfying that $2g-2+r >0$ and $\lambda$ is even.
 Then,  the stable log  twisted $\text{SUSY}_1$ curves of type $(g, r, \lambda)$
  over  log  superschemes
  and superconformal morphisms between them form  a  category 
fibered in groupoids over $\mfS \mfc \mfh_{/S_0}^{\circledS \mr{log}}$:
\begin{align} \label{EE12}
{^{\S_1}  \overline{\mfM}}_{g,r, \lambda}^{\circledS \mr{log}} \hspace{2mm}  &\migi \mfS \mfc \mfh_{/S_0}^{\circledS \mr{log}}  \\
{^{\S_1}  \mfX}^{\circledS \bigstar} \ \text{(over $S^{\circledS \mr{log}}$)} \hspace{2mm} &\mapsto  \hspace{5mm} S^{\circledS \mr{log}}. \notag
\end{align}
One verifies from a standard argument in descent theory  that ${^{\S_1}  \overline{\mfM}}_{g,r, \lambda}^{\circledS \mr{log}} $ forms  a stack with respect to the strict super\'{e}tale pretopology in $\mfS \mfc \mfh_{/S_0}^{\circledS \mr{log}}$. 
We shall denote by
\begin{align} \label{EE11}
({^{\S_1}  \overline{\mfM}}_{g,r, \lambda})_t^{\mr{log}}
\end{align}
the restriction  of ${^{\S_1}  \overline{\mfM}}_{g,r, \lambda}^{\circledS \mr{log}}$ to the full subcategory $\mfS \mfc\mfh_{/S_0}^{\mr{log}} \subseteq \mfS \mfc \mfh_{/S_0}^{\circledS \mr{log}}$.
The assignment 
 from each
stable log twisted $\text{SUSY}_1$ curve  over an fs log scheme  to its underlying pointed twisted curve
  determines a  morphism
  $({^{\S_1}  \overline{\mfM}}_{g,r, \lambda})_t^{\mr{log}} \migi  {^\mr{tw} \overline{\mfM}}_{g,r,\lambda}$; it extends to a morphism 
\begin{align} \label{e010}
({^{\S_1}  \overline{\mfM}}_{g,r, \lambda})_t^{\mr{log}} \migi  {^\mr{tw} \overline{\mfM}}^\mr{log}_{g,r,\lambda}
\end{align}
   of  log stacks.


\vspace{10mm}
\section{Superconformal structure v.s. spin structure}

This section is devoted  to understand the structure of  the reduced stack $({^{\S_1}  \overline{\mfM}}_{g,r, \lambda})_t^{\mr{log}}$ of ${^{\S_1}  \overline{\mfM}}_{g,r, \lambda}^{\circledS \mr{log}}$.
The point is that to giving a  pointed log twisted $\text{SUSY}_1$ curve over  a log scheme   is, via a natural procedure,  equivalent to giving a pointed log twisted curve equipped with an additional data called a pointed    spin structure (cf. Definition \ref{De2}).
Thus,  if  ${^{\mr{tw}} \overline{\mfM}}_{g,r, \lambda, \mr{spin}}^\mr{log}$ (cf. (\ref{EE33})) denotes the moduli stack classifying  $\lambda$-stable log  twisted  curves of type $(g,r)$  equipped with a parabolic spin structure, then it is canonically isomorphic to   $({^{\S_1} \overline{\mfM}}_{g,r, \lambda})_t^\mr{log}$ of ${^\S \overline{\mfM}}_{g,r, \lambda}^{\circledS \mr{log}}$, as shown in  Proposition \ref{P66}.

In the following, {\it we suppose that $r$ is even}.

\vspace{5mm}

\subsection{Parabolic spin structures} \label{S41}
\leavevmode\\
\vspace{-4mm}

Let $\underline{S}^\mr{log}$ be an fs log  scheme and $\underline{\mfX}^\bigstar  := (\underline{X}/\underline{S}, \{[ \underline{\sigma}_i]\}_{i=1}^r)$ be an $r$-pointed  twisted curve of genus $g$ over the underlying scheme $\underline{S}$ of $\underline{S}^\mr{log}$.
Hence, by the discussion preceding Definition \ref{D03}, we have a log smooth morphism $\underline{X}^{\bigstar \text{-} \mr{log}} \migi \underline{S}^{\mr{log}}$ and 
$\Omega_{\underline{X}^{\bigstar \text{-} \mr{log}}/\underline{S}^{\mr{log}}}$ is a line bundle of total degree $2g-2+r$.
Note that for each $i \in \{1, \cdots, r \}$,  there exists a canonical isomorphism
\begin{align}
\Lambda_i : \underline{\sigma}_i^* (\Omega_{\underline{X}^{\bigstar \text{-} \mr{log}}/\underline{S}^{\mr{log}}}) \isom \mcO_{\underline{S}}
\end{align}
 which maps any local section of the form
$\underline{\sigma}_i^* (d \mr{log} (x))$ to $1 \in \mcO_{\underline{S}}$, where $x$ is a local  
function defining the closed substack $[ \underline{\sigma}_i]$ of $\underline{X}$.
We shall write 
\begin{equation}
\mfS \mfp \mfi \mfn_{\underline{\mfX}^\bigstar}
\end{equation}
for the groupoid 
defined as follows:
\begin{itemize}
\item[$\bullet$]
 The {\it  objects} in $\mfS \mfp \mfi \mfn_{\underline{\mfX}^\bigstar}$ are pairs
$(\mcL, \eta)$,
 where $\mcL$ denotes a line bundle on $\underline{X}$  such that $\underline{\sigma}_i^*(\mcL) \cong  \mcO_{\underline{S}}$ for any $i \in \{ 1, \cdots, r \}$  and  $\eta$ denotes  an isomorphism $\mcL^{\otimes 2} \isom \Omega_{\underline{X}^{\bigstar \text{-} \mr{log}}/\underline{S}^{\mr{log}}}$.
\item[$\bullet$]
The {\it morphisms} from   $(\mcL, \eta)$ to  $(\mcL', \eta')$ (where both $(\mcL, \eta)$ and $(\mcL', \eta')$ are objects in $\mfS \mfp \mfi \mfn_{\underline{\mfX}^\bigstar}$) are 
 isomorphisms $\iota : \mcL_1 \isom \mcL_2$ satisfying the equality $\eta_2 \circ \iota^{\otimes 2} = \eta_1$.
\end{itemize}

\vspace{3mm}
\bde
\label{De2}\leavevmode\\
\ \ \
We shall refer to such a pair   $(\mcL, \eta)$ (i.e., an object of $\mfS \mfp \mfi \mfn_{\underline{\mfX}^\bigstar}$) as a {\bf pointed  spin structure} on $\underline{\mfX}^\bigstar$.
 \ede

\vspace{3mm}
\begin{rema} \label{r9009} \leavevmode\\
 \ \ \ 
 Suppose that we are given a line bundle $\mcL_0$ on $\underline{X}$ together with an isomorphism $\eta_0 : \mcL_0^{\otimes 2} \isom \Omega_{\underline{X}^{\bigstar \text{-} \mr{log}}/\underline{S}^{\mr{log}}}$.
Since  the composite   $\Lambda_i \circ  \underline{\sigma}_i^*(\eta_0)$ (for  each $i \in \{ 1, \cdots, r \}$)  is  an isomorphism $\underline{\sigma}_i^*(\mcL)^{\otimes 2} \isom \mcO_{\underline{S}}$,  the line bundle  $\underline{\sigma}_i^*(\mcL)$ defines a $\mu_2$-torsor over $\underline{S}$.
Hence, 
after possibly base-changing $\underline{\mfX}^\bigstar$ via an \'{e}tale covering  $\underline{S}'  \migi \underline{S}$ of  $\underline{S}$,
the pair $(\mcL_0, \eta_0)$ becomes a pointed spin structure on $\underline{\mfX}^\bigstar$.
Indeed, if $\underline{S}_i$ denotes the total space of the $\mu_2$-torsor corresponding to $\underline{\sigma}_i^*(\mcL)$,
then it suffices to choose the \'{e}tale covering  $\underline{S}' = \underline{S}_1 \times_{\underline{S}} \underline{S}_2 \times_{\underline{S}} \cdots \times_{\underline{S}} \underline{S}_r$ of $\underline{S}$.
 \end{rema}
\vspace{3mm}

Denote by
\begin{equation} \label{EE33}
{^{\mr{tw}} \overline{\mfM}}_{g,r, \lambda, \mr{spin}}
\end{equation}
the  category fibered in groupoids over ${^\mr{tw} \overline{\mfM}}_{g,r}$ whose fiber over
$\underline{S} \migi {^\mr{tw} \overline{\mfM}}_{g,r}$ (where $\underline{S}$ is a scheme) classifying an $r$-pointed  $\lambda$-stable  twisted curve $\underline{\mfX}^\bigstar$  is the groupoid $\mfS \mfp \mfi \mfn_{\underline{\mfX}^\bigstar}$.
One verifies from  ~\cite{Chi1}, Corollary 4.2.8 (and the fact that $\mr{deg} (\Omega_{\underline{X}^{\bigstar \text{-} \mr{log}}/\underline{S}^{\underline{f}\text{-}\mr{log}}}) = 2g-2+r$ is even)  that ${^{\mr{tw}} \overline{\mfM}}_{g,r, \lambda, \mr{spin}}$  may be represented by
 a smooth  proper Deligne-Mumford stack over $S_0$ of relative dimension $3g-3+r$ and the forgetting   morphism
\begin{align} \label{e040}
 {^{\mr{tw}} \overline{\mfM}}_{g,r, \lambda, \mr{spin}} \migi {^\mr{tw} \overline{\mfM}}_{g,r, \lambda}
\end{align}
 is finite and \'{e}tale.
 Indeed,  according to the discussion in Remark \ref{r9009}, 
 ${^{\mr{tw}} \overline{\mfM}}_{g,r, \lambda, \mr{spin}}$ turns out to be  finite and \'{e}tale over the moduli stack classifying  $r$-pointed  $\lambda$-stable twisted curves  $(\underline{X}/\underline{S}, \{ [\underline{\sigma}_i] \}_{i=1}^r)$ of genus $g$ equipped with a square root of $\Omega_{\underline{X}^{\bigstar \text{-} \mr{log}}/\underline{S}^{\underline{f}\text{-}\mr{log}}}$.
We equip ${^{\mr{tw}} \overline{\mfM}}_{g,r, \lambda, \mr{spin}}$ with 
 the  log structure pulled-back from ${^\mr{tw} \overline{\mfM}}_{g,r, \lambda}^\mr{log}$.
 Write 
 \begin{align}
 {^{\mr{tw}} \overline{\mfM}}^\mr{log}_{g,r, \lambda, \mr{spin}}
 \end{align}
for the resulting fs log stack (hence, (\ref{e040}) extends to ${^{\mr{tw}} \overline{\mfM}}^\mr{log}_{g,r, \lambda, \mr{spin}} \migi {^\mr{tw} \overline{\mfM}}_{g,r, \lambda}^\mr{log}$).

\vspace{5mm}

\subsection{From $({^{\S_1} \overline{\mfM}}_{g,r, \lambda})_t^\mr{log}$ to ${^\mr{tw} \overline{\mfM}}_{g,r, \lambda, \mr{spin}}^\mr{log}$} \label{S42}
\leavevmode\\
\vspace{-4mm}

The main goal of this section is to prove Proposition \ref{P66} described at the end of this section, i.e., to construct an equivalence of categories $({^{\S_1} \overline{\mfM}}_{g,r, \lambda})_t^\mr{log} \isom {^\mr{tw} \overline{\mfM}}_{g,r, \lambda, \mr{spin}}^\mr{log}$.
To this end, we construct first 
a morphism
$({^{\S_1} \overline{\mfM}}_{g,r, \lambda})_t^\mr{log} \migi {^\mr{tw} \overline{\mfM}}_{g,r, \lambda, \mr{spin}}^\mr{log}$ over ${^\mr{tw} \overline{\mfM}}_{g,r,\lambda}^\mr{log}$.

Let $\underline{S}^\mr{log}$ be  an fs log scheme and  $\underline{S}^\mr{log} \migi ({^{\S_1} \overline{\mfM}}_{g,r, \lambda})_t^\mr{log}$  a morphism classifying 
  a stable log twisted  $\text{SUSY}_1$ curve 
  \begin{align}
  {^{\S_1} \mfX}^{\circledS \bigstar} := (f^{\circledS \mr{log}} : X^{\circledS \mr{log}} \migi \underline{S}^\mr{log}, \{ [\sigma^{\circledS}_i] \}_{i=1}^r, \mcD)
  \end{align}
   of type $(g, r, \lambda)$  over $\underline{S}^\mr{log}$.
The morphism  $\gamma_X : X_t \migi X_b$ is an isomorphism,  
and  allows us to  identify
   the   $r$-pointed $\lambda$-stable  twisted curve
\begin{align}
\mfX_b^\bigstar := (f_b : X_b \migi \underline{S}, \{ [(\sigma_i)_b]\}_{i=1}^r)
\end{align}
     with    the underlying pointed twisted curve  of $(X^{\circledS \mr{log}}/\underline{S}^\mr{log}, \{ [\sigma^{\circledS}_i]\}_{i=1}^r)$.
One verifies that there exist
a line bundle $\mcL^\bigstar$ on $X_b$ and an isomorphism $\Upsilon : X^{\circledS} \isom \langle X_b,  \mcL^\bigstar \rangle^\circledS$ over $\underline{S}$
which sends $[\sigma_i^\circledS]$ (for each $i \in \{1, \cdots, r\}$)  to the closed subsuperscheme of $\langle X_b,  \mcL^\bigstar \rangle^\circledS$ represented  by the closed immersion $\langle \underline{S}, (\sigma_i)_b^*(\mcL^\bigstar)\rangle \migi \langle X_b, \mcL^\bigstar \rangle^\circledS$ extending $(\sigma_i)_b$.
In particular, we have  
\begin{align} \label{FF09}
(\sigma_i)_b^*(\mcL^\bigstar) \cong \mcO_{\underline{S}}.
\end{align}
The isomorphism  $\Upsilon$   gives rise to an isomorphism
\begin{align} \label{e06}
\Omega_{X^{\circledS \bigstar \text{-}\mr{log}}/\underline{S}^\mr{log}} \isom ( \mcO_{X^\circledS} \otimes_{\mcO_{X_b}}\Omega_{X_b^{\bigstar\text{-}\mr{log}}/\underline{S}^\mr{log}} ) \oplus ( \mcO_{X^\circledS} \otimes_{\mcO_{X_b}}  \mcL^\bigstar )
\end{align}
of $\mcO_{X^\circledS}$-supermodules, where, in the right-hand side,   the sections
of the forms $(1 \otimes a, 0)$ and $(0, 1\otimes  b)$ (for some $a \in \Omega_{X_b^{\bigstar\text{-}\mr{log}}/\underline{S}^\mr{log}}$ and $b \in \mcL^\bigstar$)
  are defined to be  bosonic and fermionic  sections respectively.
Consider its dual
\begin{align} \label{e0998}
\mcT_{X^{\circledS \bigstar \text{-}\mr{log}}/\underline{S}^\mr{log}} \isom  (\mcO_{X^\circledS} \otimes_{\mcO_{X_b}}\mcT_{X_b^{\bigstar\text{-}\mr{log}}/\underline{S}^\mr{log}} ) \oplus ( \mcO_{X^\circledS} \otimes_{\mcO_{X_b}}  \mcL^{\bigstar \vee}).
\end{align}
It follows from Proposition \ref{p01033} (i) that 
the composite  morphism 
\begin{align} \label{e077}
\mcD &  \  \migiincl \mcT_{X^{\circledS \bigstar \text{-}\mr{log}}/\underline{S}^\mr{log}}  \\
& \stackrel{(\ref{e0998})}{\isom}  (\mcO_{X^\circledS} \otimes_{\mcO_{X_b}}\mcT_{X_b^{\bigstar\text{-}\mr{log}}/\underline{S}^\mr{log}} ) \oplus ( \mcO_{X^\circledS} \otimes_{\mcO_{X_b}}  \mcL^{\bigstar \vee}) \notag \\
&  \  \migisurj   \mcO_{X^\circledS} \otimes_{\mcO_{X_b}}  \mcL^{\bigstar \vee} \notag 
\end{align}
(where the third morphism denotes the projection to the second factor)
between  supervector bundles of superrank $0|1$
 is surjective, and hence,  an isomorphism. 
Moreover, we have a composite isomorphism
\begin{align} \label{e088}
\mcO_{X^\circledS} \otimes_{\mcO_{X_b}}\mcT_{X_b^{\bigstar\text{-}\mr{log}}/\underline{S}^\mr{log}} &    \migi   (\mcO_{X^\circledS} \otimes_{\mcO_{X_b}}\mcT_{X_b^{\bigstar\text{-}\mr{log}}/\underline{S}^\mr{log}} ) \oplus ( \mcO_{X^\circledS} \otimes_{\mcO_{X_b}}  \mcL^{\bigstar \vee}) \\
& \stackrel{(\ref{e0998})^{-1}}{\isom} \mcT_{X^{\circledS \bigstar \text{-}\mr{log}}/\underline{S}^\mr{log}} \notag \\
& \migisurj \mcT_{X^{\circledS \bigstar \text{-}\mr{log}}/\underline{S}^\mr{log}} /\mcD \notag
\end{align}
(where the first morphism denotes the inclusion into the first factor).
The isomorphism (\ref{e032}) in our situation  may be described  
 as an isomorphism
\begin{align} \label{e036}
(\mcO_{X^\circledS} \otimes_{\mcO_{X_b}}  (\mcL^{\bigstar \vee})^{\otimes 2}  =) \ (\mcO_{X^\circledS} \otimes_{\mcO_{X_b}}  \mcL^{\bigstar \vee})^{\otimes 2} \isom \mcO_{X^\circledS} \otimes_{\mcO_{X_b}}\mcT_{X_b^{\bigstar\text{-}\mr{log}}/\underline{S}^\mr{log}}
\end{align}
 via the composite isomorphisms (\ref{e077}) and (\ref{e088}).
The restriction of  this isomorphism to the bosonic part (and taking its dual) becomes  an isomorphism $\eta^\bigstar : (\mcL^\bigstar)^{\otimes 2} \isom \Omega_{X_b^{\bigstar\text{-}\mr{log}}/\underline{S}^\mr{log}}$.
Thus, the pair $(\mcL^\bigstar, \eta^\bigstar)$  forms  (thanks to (\ref{FF09}))
 a pointed  spin structure on $\mfX_b^\bigstar$.
 If $\underline{S} \migi {^\mr{tw} \overline{\mfM}}_{g,r, \lambda, \mr{spin}}$ denotes the  classifying morphism of $(\mcL^\bigstar, \eta^\bigstar)$, then it
  extends uniquely to a morphism $\underline{S}^\mr{log} \migi {^\mr{tw} \overline{\mfM}}^\mr{log}_{g,r, \lambda, \mr{spin}}$  over ${^\mr{tw} \overline{\mfM}}_{g,r,\lambda}^\mr{log}$.
The assignment ${^{\S_1} \mfX}^{\circledS \bigstar} \mapsto (\mcL^\bigstar, \eta^\bigstar)$ is functorial with respect to $\underline{S}^\mr{log}$ and hence, determines  a morphism 
\begin{align} \label{e445}
({^{\S_1} \overline{\mfM}}_{g,r, \lambda})_t^\mr{log} \migi {^\mr{tw} \overline{\mfM}}_{g,r, \lambda, \mr{spin}}^\mr{log}
\end{align}
 over ${^\mr{tw} \overline{\mfM}}_{g,r,\lambda}^\mr{log}$.

\vspace{5mm}

\subsection{From ${^\mr{tw} \overline{\mfM}}_{g,r, \lambda, \mr{spin}}^\mr{log}$  to $({^{\S_1} \overline{\mfM}}_{g,r, \lambda})_t^\mr{log}$} \label{S43}
\leavevmode\\
\vspace{-4mm}

Conversely, we shall construct a morphism ${^\mr{tw} \overline{\mfM}}_{g,r, \lambda, \mr{spin}}^\mr{log} \migi ({^{\S_1} \overline{\mfM}}_{g,r,\lambda})^\mr{log}_t$.
Let $\underline{S}^\mr{log}$ be an fs log scheme and 
$\underline{S}^\mr{log} \migi {^\mr{tw} \overline{\mfM}}_{g,r,\lambda, \mr{spin}}^\mr{log}$  a morphism whose underlying morphism $\underline{S} \migi {^\mr{tw} \overline{\mfM}}_{g,r,\lambda, \mr{spin}}$ classifies 
 a spin structure $(\mcL, \eta)$ on an $r$-pointed $\lambda$-stable  twisted curve $\underline{\mfX}^\bigstar := (\underline{f} : \underline{X} \migi \underline{S}, \{ [\underline{\sigma}_i]\}_{i=1}^r)$. (In particular, we have a morphism $\underline{S}^\mr{log} \migi \underline{S}^{\underline{f} \text{-} \mr{log}}$.)
By fixing an isomorphism $\underline{\sigma}_i^*(\mcL) \isom \mcO_{\underline{S}}$ (for each $i = 1, \cdots, r$),
we obtain a composite closed immersion
\begin{align}
\underline{\sigma}^\circledS_{(\mcL, \eta), i} : 
 \mbA^{0|1}_{\underline{S}} \isom \langle \underline{S}, \underline{\sigma}_i^*(\mcL) \rangle^\circledS \migi 
\langle \underline{X}, \mcL \rangle^\circledS
\end{align}
 extending $\underline{\sigma}_i$.
We shall write  $\underline{X}^\mr{log} := \underline{X}^{\underline{f} \text{-} \mr{log}} \times_{\underline{X}^{\underline{f} \text{-} \mr{log}}} \underline{S}^\mr{log}$, and hence,  obtain 
a log superstack $\langle \underline{X}, \mcL \rangle^{\circledS \mr{log}}$ (cf. (\ref{E13})) over $\underline{S}^\mr{log}$.
The collection of data 
\begin{align}
\underline{\mfX}^{\circledS \bigstar}_{(\mcL, \eta)} := (\langle \underline{X}, \mcL \rangle^{\circledS \mr{log}}/\underline{S}^\mr{log}, \{ [\underline{\sigma}^\circledS_{(\mcL, \eta), i}]\}_{i=1}^r)
\end{align}
 forms an $r$-pointed  log  twisted $(1|1)$-curve of genus $g$  over $\underline{S}^\mr{log}$.
Since $\mcO_{\langle \underline{X}, \mcL \rangle^\circledS} \cong \mcO_{\underline{X}}\oplus \mcL$,  we obtain (cf. (\ref{e06}))  a composite  isomorphism
\begin{align} \label{e057}
 & \ \ \  \ \ \mcT_{\langle \underline{X}, \mcL \rangle^{\circledS \bigstar \text{-}\mr{log}} /\underline{S}^{\mr{log}}} \\
  & \isom  (\mcO_{\langle \underline{X}, \mcL \rangle^\circledS}  \otimes_{\mcO_X} \mcT_{\underline{X}^{\bigstar \text{-}\mr{log}}/\underline{S}^{\underline{f}\text{-}\mr{log}}})   \oplus 
 (\mcO_{\langle X, \mcL \rangle^\circledS}  \otimes_{\mcO_X} \mcL^\vee)\notag \\
& \isom \mcT_{\underline{X}^{\bigstar \text{-}\mr{log}}/\underline{S}^{\underline{f}\text{-}\mr{log}}}  \oplus \mcL^{\vee}\oplus \mcL^{\vee } \oplus \mcO_{\underline{X}} \notag
\end{align}
of $\mcO_{\underline{X}}$-modules.
Consider  the $\mcO_{\langle \underline{X}, \mcL \rangle^\circledS}$-linear injection
\begin{align} \label{E30}
 (\mcO_{\langle \underline{X}, \mcL \rangle^\circledS}  \otimes_{\mcO_{\underline{X}}} \mcL^\vee =:) \ \mcL^{\vee } \oplus \mcO_{\underline{X}} & \migiincl  \mcT_{\underline{X}^{\bigstar \text{-}\mr{log}}/\underline{S}^{\underline{f}\text{-}\mr{log}}}   \oplus \mcL^{\vee}\oplus \mcL^{\vee } \oplus \mcO_{\underline{X}} \\
(\,a, \  b\,) \hspace{3mm} &\mapsto \hspace{15mm} (\,0, \  \ a, \  \ a, \ \ b\,).\notag
\end{align}
Write $\mcD_{(\mcL, \eta)}$ for the  subsupervector bundle  (of superrank $0|1$) of $\mcT_{\langle \underline{X}, \mcL \rangle^{\circledS \bigstar \text{-}\mr{log}} /\underline{S}^{\mr{log}}}$ corresponding, via the composite isomorphism (\ref{e057}), to the image of  (\ref{E30}).
Then, one verifies immediately that the collection of data 
\begin{align}
{^\S \underline{\mfX}}^{\circledS  \bigstar}_{(\mcL, \eta)} := (\langle \underline{X}, \mcL \rangle^{\circledS \mr{log}}/\underline{S}^\mr{log}, \{ [\underline{\sigma}^\circledS_{(\mcL, \eta), i}]\}_{i=1}^r, \mcD_{(\mcL, \eta)})
\end{align}
forms a stable log  twisted $\text{SUSY}_1$  curve  of type $(g, r, \lambda)$  over $\underline{S}^\mr{log}$ whose underlying pointed  twisted curve is isomorphic to $\underline{\mfX}^\bigstar$.
It determines a morphism $\underline{S}^\mr{log} \migi ({^{\S_1} \overline{\mfM}}_{g,r,\lambda})^\mr{log}_t$.
By varying $\underline{S}^\mr{log}$ with the various fs log schemes, we obtain a morphism
\begin{align} \label{e345}
{^\mr{tw} \overline{\mfM}}_{g,r, \lambda, \mr{spin}}^\mr{log} \migi ({^{\S_1} \overline{\mfM}}_{g,r,\lambda})^\mr{log}_t  
\end{align}
 over ${^\mr{tw} \overline{\mfM}}_{g,r,\lambda}^\mr{log}$.

\vspace{3mm}
\begin{rema} \label{r4} \leavevmode\\
 \ \ \ 
 We keep the above notation.
 By applying Proposition \ref{p010}, we have a composite  isomorphism of $\underline{f}^{-1} (\mcO_{\underline{S}})$-modules: 
 \begin{align} \label{e059}
  \mcT^{\mcD_{(\mcL, \eta)}}_{\langle \underline{X}, \mcL \rangle^{\circledS \bigstar \text{-}\mr{log}}/\underline{S}^\mr{log}} & \isom (\mcD_{(\mcL, \eta)})^{\otimes 2} \\
&  \isom \mcO_{\langle \underline{X}, \mcL \rangle^\circledS} \otimes_{\mcO_{\underline{X}}} (\mcL^{\vee})^{\otimes 2} \notag  \\
&  \isom  \mcT_{\underline{X}^{\bigstar \text{-}\mr{log}}/\underline{S}^{\underline{f} \text{-}\mr{log}}}  \oplus \mcL^\vee,  \notag 
 \end{align}
where the second isomorphism follows from the definition of $\mcD_{(\mcL, \eta)}$ (cf. (\ref{E30})).
Let us  equip $ \mcT^{\mcD_{(\mcL, \eta)}}_{\langle \underline{X}, \mcL \rangle^{\circledS \bigstar \text{-}\mr{log}}/\underline{S}^\mr{log}}$ with a structure of  $\mcO_{\langle \underline{X}, \mcL \rangle^\circledS}$-supermodule via the first isomorphism in (\ref{e059}).
Then, 
(\ref{e059}) induces  two isomorphisms of $\mcO_{\underline{X}}$-modules
\begin{align} \label{E006}
(  \mcT^{\mcD_{(\mcL, \eta)}}_{\langle \underline{X}, \mcL \rangle^{\circledS \bigstar \text{-}\mr{log}}/\underline{S}^\mr{log}})_b \cong  \mcT_{\underline{X}^{\bigstar \text{-}\mr{log}}/\underline{S}^{\underline{f} \text{-}\mr{log}}}, \hspace{5mm}   (  \mcT^{\mcD_{(\mcL, \eta)}}_{\langle \underline{X}, \mcL \rangle^{\circledS \bigstar \text{-}\mr{log}}/\underline{S}^\mr{log}})_f \cong \mcL^\vee.
\end{align}
In particular, we have
\begin{align} \label{E066}
\mr{deg} ((  \mcT^{\mcD_{(\mcL, \eta)}}_{\langle \underline{X}, \mcL \rangle^{\circledS \bigstar \text{-}\mr{log}}/\underline{S}^\mr{log}})_b) = -2g+2-r, \hspace{5mm}  \mr{deg} ((  \mcT^{\mcD_{(\mcL, \eta)}}_{\langle \underline{X}, \mcL \rangle^{\circledS \bigstar \text{-}\mr{log}}/\underline{S}^\mr{log}})_f) = -g+1-\frac{r}{2}.
\end{align}
 The inclusion $ \mcT^{\mcD_{(\mcL, \eta)}}_{\langle \underline{X}, \mcL \rangle^{\circledS \bigstar \text{-}\mr{log}}/\underline{S}^\mr{log}}  \migiincl    \mcT_{\langle \underline{X}, \mcL \rangle^{\circledS \bigstar \text{-}\mr{log}}/\underline{S}^\mr{log}}$ corresponds, via the isomorphisms (\ref{e057}) and (\ref{e059}), to
 the  inclusion
 \begin{align}
\mcT_{\underline{X}^{\bigstar \text{-}\mr{log}}/\underline{S}^{\underline{f} \text{-}\mr{log}}}   \oplus \mcL^\vee & \migi  \mcT_{\underline{X}^{\bigstar \text{-}\mr{log}}/\underline{S}^{\underline{f} \text{-}\mr{log}}}   \oplus \mcL^{\vee}\oplus \mcL^{\vee } \oplus \mcO_{\underline{X}} \\
(\,a, \  b\,) \hspace{3mm} &\mapsto \hspace{22mm} (\,a, \ b, \ b, \ 0\,).\notag
 \end{align}
  \end{rema}

\vspace{3mm}
\bpr 
\label{Pt66}
\leavevmode\\
 \ \ \ 
Let ${^{\S_1} \mfX}^{\circledS \bigstar} := (f^{\circledS \mr{log}} : X^{\circledS \mr{log}} \migi S^{\circledS \mr{log}}, \{ [\sigma_i^\circledS] \}_{i=1}^r, \mcD)$ be a stable  log twisted $\text{SUSY}_1$ curve of type  $(g, r, \lambda)$ over $S^{\circledS \mr{log}}$.
Then, $\mbR^2f_{*}^\circledS (\mcT^\mcD_{X^{\circledS \bigstar \text{-}\mr{log}}/S^{\circledS \mr{log}}}) =0$.
Also, the $\mcO_{S^\circledS}$-supermodule $\mbR^1f_{*}^\circledS (\mcT^\mcD_{X^{\circledS \bigstar \text{-}\mr{log}}/S^{\circledS \mr{log}}})$ is locally free of rank  $3g-3+r | 2g-2+\frac{r}{2}$ and the formulation of $\mbR^1f_{*}^\circledS (\mcT^\mcD_{X^{\circledS \bigstar \text{-}\mr{log}}/S^{\circledS \mr{log}}})$ commutes with base-change with respect to $S^\circledS$.
Moreover,  if $S^\circledS$ is affine (i.e., $S^{\circledS} = S_b$) and $\mcF$ is a quasi-coherent  $\mcO_{S_b}$-module, then 
the natural morphism
\begin{align}
\mbR^1f_{b*} ((\mcT^\mcD_{X^{\circledS \bigstar \text{-}\mr{log}}/S_b^{\mr{log}}})_{(-)}) \otimes \mcF \migi \mbR^1f_{b*} ((\mcT^\mcD_{X^{\circledS \bigstar \text{-}\mr{log}}/S_b^{\mr{log}}} \otimes_{f^{-1}_b (\mcO_{S_b})} f^{-1}_b (\mcF))_{(-)}),
\end{align}
where $(-)$ denotes either ``$b$" or ``$f$",  is an  isomorphism.
 \epr
\begin{proof}
In the following, we consider  $\mcT^\mcD_{X^{\circledS \bigstar \text{-}\mr{log}}/S^{\circledS \mr{log}}}$ as being equipped with a structure of $\mcO_{X^{\circledS}}$-supermodule by transposing the structure of $\mcO_{X^{\circledS}}$-supermodule on $\mcD^{\otimes 2}$ via the isomorphism $\mcT^\mcD_{X^{\circledS \bigstar \text{-}\mr{log}}/S^{\circledS \mr{log}}} \isom \mcD^{\otimes 2}$ obtained in Proposition \ref{p010}.
Let us take  an algebraically closed field $k$ and a morphism $v^{\circledS} : \mr{Spec} (k) \migi S^\circledS$, where we equip $\mr{Spec} (k)$ with a log structure pulled-back from $S^{\circledS \mr{log}}$ via $v^\circledS$. Write
$\mr{Spec}(k)^\mr{log}$ for the resulting log scheme and $v^{\circledS \mr{log}} : \mr{Spec}(k)^\mr{log} \migi S^{\circledS \mr{log}}$ the morphism extending $v^\circledS$.
Also, write $(X_v^{\circledS \mr{log}}/k^{\mr{log}}, \{ [\sigma_{i, v}^\circledS ] \}_{i=1}^r, \mcD_v)$
for the stable log twisted $\text{SUSY}_1$ curve over $\mr{Spec}(k)^\mr{log}$ defined as  the base-change of ${^{\S_1} \mfX}^{\circledS \bigstar}$ via $v^{\circledS \mr{log}}$.
It follows from 
Theorem \ref{p04033} that
$H^2 ((X_v)_b, \mcT^{\mcD_v}_{X_v^{\circledS \bigstar \text{-}\mr{log}}/k^{\mr{log}}}) =0$.
By replacing $v^\circledS$ with the  various points of $S^\circledS$, one verifies   from ~\cite{Har}, Ch.\,III, Theorem 12.11, (a) (or ~\cite{QL}, \S\,5, Remark 3.21, (c)) that 
$\mbR^2f_{b*} (\mcT^\mcD_{X^{\circledS \bigstar \text{-}\mr{log}}/S^{\circledS \mr{log}}}) =0$.
(Here, note that although the result of {\it loc.\,cit.} deals only  with the case of schemes, one may prove, by a similar argument, its  analogous assertion   for  the case of superstacks.
 Thus, in the proof of this proposition, we  apply the result of {\it loc.\,cit.} as the  result  corresponding to  the case of superstacks.)
Hence, the natural morphism 
 \begin{align}
 v^{\circledS *}(\mbR^1f_{b*} (\mcT^\mcD_{X^{\circledS \bigstar \text{-}\mr{log}}/S^{\circledS \mr{log}}})) \migi H^1 ((X_v)_b, \mcT^{\mcD_v}_{X_v^{\circledS \bigstar \text{-}\mr{log}}/k^{\mr{log}}})
 \end{align}
  is surjective for all $v^\circledS$ (cf. ~\cite{Har}, Ch.\,III, Theorem 12.11, (b)) and hence, is an isomorphism.
Moreover, the last reference or  ~\cite{QL}, \S\,5, Remark 3.21, (c) (applied to the case $i=1$) shows that $\mbR^1f_{b*} (\mcT^\mcD_{X^{\circledS \bigstar \text{-}\mr{log}}/S^{\circledS \mr{log}}})$ is locally free.
The rank of this $\mcO_{S^\circledS}$-supermodule may be calculated by the dimension of $H^1 ((X_v)_b, \mcT^{\mcD_v}_{X_v^{\circledS \bigstar \text{-}\mr{log}}/k^{\mr{log}}})$.
But, it follows from (\ref{E066}) and Theorem \ref{p04033} that
\begin{align}
\mr{dim}_k (H^1 ((X_v)_b, (\mcT^{\mcD_v}_{X_v^{\circledS \bigstar \text{-}\mr{log}}/k^{\mr{log}}})_b)) = 3g-3+r
\end{align}
and
\begin{align}
\mr{dim}_k (H^1 ((X_v)_b, (\mcT^{\mcD_v}_{X_v^{\circledS \bigstar \text{-}\mr{log}}/k^{\mr{log}}})_f)) = 2g-2+\frac{r}{2},
\end{align}
as desired. The commutativity (with respect to base-change over superschemes over $S^\circledS$) of the formulation of  $\mbR^1f_{*}^\circledS (\mcT^\mcD_{X^{\circledS \bigstar \text{-}\mr{log}}/S^{\circledS \mr{log}}})$ follows  from the above discussion and the discussion in {\it loc.\,cit.}.
 This completes the proof of the former assertion.

The latter assertion follows from ~\cite{QL}, \S\,5, Remark 3.21, (c).
\end{proof}

\vspace{5mm}

\subsection{ ${^\mr{tw} \overline{\mfM}}_{g,r, \lambda, \mr{spin}}^\mr{log}$  is isomorphic to   $({^{\S_1} \overline{\mfM}}_{g,r, \lambda})_t^\mr{log}$} \label{S44}
\leavevmode\\
\vspace{-4mm}

One verifies that the morphisms (\ref{e445}) and (\ref{e345})  obtained previously
   are the  inverse  morphisms of each other.
Thus, we have the following Proposition \ref{P66}.
In particular, for each $r$-pointed $\lambda$-stable twisted curve $\underline{\mfX}^\bigstar : = (\underline{f} : \underline{X} \migi \underline{S}, \{ [\underline{\sigma}_{i} ]\}_{i=1}^r)$ of genus $g$ over a scheme $\underline{S}$, there exists canonically an equivalence of categories between   $\mfS \mfp \mfi \mfn_{\underline{\mfX}^\bigstar}$ and the groupoid of stable  log twisted $\text{SUSY}_1$ curves  ${^{\S_1} \mfX}^{\circledS \bigstar}$  over $S^{\underline{f} \text{-} \mr{log}}$ having $\underline{\mfX}^\bigstar$ as the underlying  pointed twisted curve.

\vspace{3mm}
\bpr 
\label{P66}
\leavevmode\\
 \ \ \ 
There exists a canonical isomorphism of fibered categories 
\begin{align} \label{E37}
({^{\S_1} \overline{\mfM}}_{g,r,\lambda})^\mr{log}_t  \isom {^\mr{tw} \overline{\mfM}}_{g,r, \lambda, \mr{spin}}^\mr{log}
\end{align}
over ${^\mr{tw} \overline{\mfM}}_{g,r, \lambda}^\mr{log}$.
In particular, $({^{\S_1} \overline{\mfM}}_{g,r,\lambda})^\mr{log}_t$ may be represented by a proper smooth Deligne-Mumford stack over $S_0$ of relative dimension $3g -3 +r$.
 \epr

\vspace{10mm}
\section{Deformations of  stable log twisted $\text{SUSY}_1$ curves} \vspace{3mm}

In this final section, we prove the  main assertion, i.e., Theorem A.
As discussed in \S\,\ref{S54}, the main step of  the proof  is  to construct  a canonical fermionic deformation of $ {^{\mr{tw}} \overline{\mfM}}_{g,r, \lambda, \mr{spin}}^\mr{log}$ ($\cong ({^{\S_1} \overline{\mfM}}_{g,r, \lambda})^\mr{log}_t$ by Proposition \ref{P66})  in a way that a universal stable log twisted $\text{SUSY}_1$ curve  exists (uniquely).
To this end, 
we develop (in \S\S\,5.1-5.2)  log smooth deformation theory concerning  log twisted $\text{SUSY}_1$  curves.
By applying the results obtained in these  discussions,
one may construct  (cf. Corollary \ref{p01}) a universal family of stable  log twisted $\text{SUSY}_1$ curves over  a  fermionic deformation of
a representation (in the sense of Remark \ref{r78} (i)) of $ {^{\mr{tw}} \overline{\mfM}}_{g,r, \lambda, \mr{spin}}^\mr{log}$.
It gives rise to a representation (by a groupoid in $\mfS \mfc \mfh_{/S_0}^\circledS$) of ${^{\S_1} \overline{\mfM}}_{g,r, \lambda}^{\circledS \mr{log}}$ itself
 (equipped  with a natural  log structure).
 This  implies immediately the proof of Theorem A, as desired (cf. \S\,\ref{S54} for the detailed discussion).


\vspace{5mm}

\subsection{Deformation spaces of  stable log twisted $\text{SUSY}_1$ curves} \label{S51}
\leavevmode\\
\vspace{-4mm}

Let $\widetilde{S}^{\circledS \mr{log}}$ be an fs log scheme and  
 $S^{\circledS \mr{log}}$ be a strict closed subsuperscheme of $\widetilde{S}^{\circledS \mr{log}}$ determined by  a nilpotent superideal $\mcI \subseteq \mcO_{\widetilde{S}^\circledS}$ contained in $\mcN_{\widetilde{S}^\circledS}$.
Let 
${^{\S_1} \mfX}^{\circledS \bigstar} :=(f^{\circledS \mr{log}} : X^{\circledS  \mr{log}} \migi  S^{\circledS  \mr{log}}, \{ [\sigma^\circledS_{i}]\}_{i=1}^r, \mcD)$
 be a stable log twisted  $\text{SUSY}_1$ curve of type $(g,r, \lambda)$ over $S^{\circledS \mr{log}}$.
Write ${^{\S_1} \underline{\mfX}}^{\circledS \bigstar} := (\underline{f}^{\circledS \mr{log}} : \underline{X}^{\circledS \mr{log}} \migi S^\mr{log}_t, \{ [\underline{\sigma}_i]\}_{i=1}^r , \underline{\mcD})$ for
the base-change of ${^{\S_1} \mfX}^{\circledS \bigstar}$ via the strict closed immersion $\tau_{S}^{\circledS \mr{log}} : S_t^\mr{log} \migi S^{\circledS \mr{log}}$ extending $\tau_{S}^\circledS$.
 Also, write
\begin{align}
\mr{Def}_{\widetilde{S}^{\circledS \mr{log}}} ({^{\S_1} \mfX}^{\circledS \bigstar})
\end{align}
for the set of superconformal isomorphism classes of  stable log twisted $\text{SUSY}_1$ curves of type $(g,r, \lambda)$
 over $\widetilde{S}^{\circledS \mr{log}}$ extending  ${^{\S_1} \mfX}^{\circledS \bigstar}$.

\vspace{3mm}
\bpr  [cf. ~\cite{LR}, Lemma 2.4] 
\label{P001}
\leavevmode\\
 \ \ \ 
 Suppose that $\widetilde{S}^\circledS$ is affine and  that   $\mcN_{\widetilde{S}^\circledS} \mcI =0$ (which implies that $\mcI$ is square nilpotent and  may be thought of as an $\mcO_{S_t}$-module).
Then, $\mr{Def}_{\widetilde{S}^{\circledS \mr{log}}} ({^{\S_1} \mfX}^{\circledS \bigstar})$ is nonempty and has a canonical structure 
\begin{align} \label{E003}
\mr{Def}_{\widetilde{S}^{\circledS \mr{log}}} ({^{\S_1} \mfX}^{\circledS \bigstar}) \times H^1 (\underline{X}_b, (\mcT^{\underline{\mcD}}_{\underline{X}^{\circledS \bigstar \text{-}\mr{log}}/S_t^{\mr{log}}} \otimes_{\underline{f}_b^{-1}(\mcO_{S_t})} \underline{f}^{-1}_b(\mcI))_b) & \migi \mr{Def}_{\widetilde{S}^{\circledS \mr{log}}} ({^{\S_1} \mfX}^{\circledS \bigstar}) \\
( \ a, \hspace{5mm} b \ ) \hspace{65mm}& \mapsto \hspace{5mm} a \boxplus^\ddagger b
 \notag
\end{align}
 of 
affine space modeled on $H^1 (\underline{X}_b, (\mcT^{\underline{\mcD}}_{\underline{X}^{\circledS \bigstar \text{-}\mr{log}}/S_t^{\mr{log}}} \otimes_{\underline{f}_b^{-1}(\mcO_{S^\circledS})} \underline{f}^{-1}_b(\mcI))_b)$.
Also, if 
${^{\S_1} \widetilde{\mfX}}^{\circledS \bigstar}$ is a  stable log twisted $\text{SUSY}_1$ curve in $\mr{Def}_{\widetilde{S}^{\circledS \mr{log}}} ({^{\S_1} \mfX}^{\circledS \bigstar})$, then 
 there is no nontrivial superconformal automorphism of ${^{\S_1} \widetilde{\mfX}}^{\circledS \bigstar}$ over $\widetilde{S}^{\circledS \mr{log}}$ which restricts to the identity morphism of ${^{\S_1} \mfX}^{\circledS \bigstar}$.
\epr
\begin{proof}
First, we shall prove that $\mr{Def}_{\widetilde{S}^{\circledS \mr{log}}} ({^{\S_1}\mfX}^{\circledS \bigstar})$ is nonempty.
Let us take a collection of  data
\begin{align}
\mbU_I := \{ (Y_{\gamma}^{\circledS \mr{log}} \stackrel{\pi_\gamma^{\circledS \mr{log}}}{\migi} X^{\circledS \mr{log}}, U_{\gamma}^\mr{log}, \eta^{\circledS\mr{log}}_{\gamma}, z_{\gamma}) \}_{\gamma \in I}
\end{align}
(indexed by a set $I$) obtained by applying Corollary \ref{c01033} to our  ${^{\S_1} \mfX}^{\circledS \bigstar}$. 
For each $i \in \{ 1, \cdots, r\}$, we shall write $I_i$ for the subset of $I$ consisting of elements $\gamma$ satisfying that $Y_{\gamma}^\circledS \times_{X^\circledS, \sigma^\circledS_{ i}}  \mbA^{0|1}_{S^\circledS}    \neq \emptyset$.
By Proposition  \ref{p0468},  we may assume, without loss of generality, that
for each $\gamma \in I_i$ there exists  
a pair $(\Sigma^{U_\gamma}, \sigma^{U_\gamma})$ for which 
the collection of data
\begin{align}
\mbU^\bigstar_\gamma := (Y_\gamma^{\circledS \mr{log}}, U_\gamma^{\circledS \mr{log}}, \eta_\gamma^{\circledS \mr{log}}, \Sigma^{U_\gamma}, \sigma^{U_\gamma})
\end{align}
 satisfies the condition described in Proposition \ref{p0468}.
For each $\gamma \in I$ there exists (since $Y_{\gamma}^\circledS$ is affine) a log supersmooth lifting  $\widetilde{\mbY}_\gamma :=  (\widetilde{Y}^{\circledS \mr{log}}_\gamma, \widetilde{f}_\gamma^{\circledS \mr{log}}, \widetilde{i}^{\circledS \mr{log}}_{Y_\gamma})$ of $Y^{\circledS \mr{log}}_{\gamma}$ over $\widetilde{S}^{\circledS \mr{log}}$ (cf. Corollary \ref{c0404} (iii)) and a log smooth lifting $\widetilde{U}_\gamma^\mr{log}$ of $U_\gamma^\mr{log}$ over $\widetilde{S}_b^\mr{log}$ together with an isomorphism $\widetilde{\eta}_\gamma^{\circledS \mr{log}} : \widetilde{Y}_\gamma^{\circledS \mr{log}} \isom \widetilde{U}^{\mr{log}}_\gamma \times_{\widetilde{S}_b} \mbA_{\widetilde{S}^\circledS}^{0|1}$ over $\widetilde{S}^{\circledS \mr{log}}$ lifting $\eta^{\circledS \mr{log}}_\gamma$.
Also, $z_\gamma$ lifts to an element  $\widetilde{z}_\gamma \in \mcM_{\widetilde{U}_\gamma}$.
The $\mcO_{\widetilde{U}_\gamma^\mr{log} \times_{\widetilde{S}_b}\mbA_{\widetilde{S}^\circledS}}$-subsupermodule $\mcO_{\widetilde{U}_\gamma^\mr{log} \times_{\widetilde{S}_b}\mbA_{\widetilde{S}^\circledS}} \cdot (\partial_\psi + \psi \cdot \partial_{\widetilde{z}_\gamma}) \subseteq \mcT_{\widetilde{U}_\gamma^\mr{log} \times_{\widetilde{S}_b}\mbA_{\widetilde{S}^\circledS}/\widetilde{S}^{\circledS \mr{log}}}$ defines, via $\widetilde{\eta}^{\circledS \mr{log}}_\gamma$, a superconformal structure $\widetilde{\mcD}_{\gamma}$  on $\widetilde{Y}^{\circledS \mr{log}}_\gamma/\widetilde{S}^{\circledS \mr{log}}$ extending  $\mcD |_{Y^{\circledS \mr{log}}_\gamma}$.
After possibly replacing $\mbU_I$ with its refinement in an evident sense (cf. Remark \ref{r4f828}),
we may suppose that the following three properties (i)-(iii) are satisfied:
\begin{itemize}
 \item[(i)]
For each pair $(\gamma_1, \gamma_2)$ with
$Y_{\gamma_1, \gamma_2}^{\circledS \mr{log}} := Y_{\gamma_1}^{\circledS \mr{log}} \times_{X^{\circledS \mr{log}}} Y_{\gamma_2}^{\circledS \mr{log}} \neq \emptyset$
 (hence  $Y_{\gamma_1, \gamma_2}^{\circledS \mr{log}}$ is affine),  
   there exists an isomorphism
\begin{align}
\nu^{\circledS \mr{log}}_{\gamma_1,  \gamma_2} : \widetilde{\mbY}_{\gamma_1}^{\circledS \mr{log}} |_{Y_{\gamma_1, \gamma_2}^{\circledS \mr{log}}}   \isom \widetilde{\mbY}_{\gamma_2}^{\circledS \mr{log}}  |_{Y_{\gamma_1, \gamma_2}^{\circledS \mr{log}}}
\end{align}
of log supersmooth liftings which sends 
   $\widetilde{\mcD}_{\gamma_1}  |_{Y_{\gamma_1, \gamma_2}^{\circledS \mr{log}}}$ to $\widetilde{\mcD}_{\gamma_2}  |_{Y_{\gamma_1, \gamma_2}^{\circledS \mr{log}}}$ (cf. Corollary \ref{c0404} (ii) and Proposition \ref{p01033}).
\item[(ii)]
If $\gamma \in I_i$  ($i =1, \cdots, r$), then 
 there exists a pair $(\Sigma_i^{\widetilde{U}_\gamma}, \sigma_i^{\widetilde{U}_\gamma})$ consisting of     a scheme $\Sigma_i^{\widetilde{U}_\gamma}$ over $\widetilde{S}_b$ with $\Sigma_i^{\widetilde{U}_\gamma} \times_{\widetilde{S}_b} S_b \cong \Sigma_i^{U_\gamma}$   and a closed immersion $\sigma_i^{\widetilde{U}_\gamma} : \Sigma_i^{\widetilde{U}_\gamma} \migi  \widetilde{U}_\gamma$ over $\widetilde{S}_b$ extending  $\sigma_i^{U_\gamma}$;

\item[(iii)]
 For each  pair $(\gamma_1, \gamma_2) \in I_i \times I_i$ ($i =1, \cdots, r$)  with $Y^{\circledS \mr{log}}_{\gamma_1, \gamma_2} \neq \emptyset$, then 
 the restrictions  to $Y_{\gamma_1, \gamma_2}^{\circledS \mr{log}}$ of  two composites
 \begin{align}
 \Sigma_i^{\widetilde{U}_{\gamma_l}} \times_{\widetilde{S}_b} \mbA^{0|1}_{\widetilde{S}^\circledS} \stackrel{\sigma_i^{\widetilde{U}_\gamma}}{\migi} \widetilde{U}_{\gamma_l}  \times_{\widetilde{S}_b} \mbA^{0|1}_{\widetilde{S}^\circledS} \stackrel{(\widetilde{\eta}_{\gamma_l}^\circledS)^{-1}}{\migi} \widetilde{Y}^{\circledS}_{\gamma_l}
 \end{align}
   ($l =1, 2$) are compatible (in an evident sense) via $\nu^{\circledS \mr{log}}_{\gamma_1,  \gamma_2}$. 
\end{itemize}
\vspace{1mm}

If a triple $(\gamma_1, \gamma_2, \gamma_3) \in I^{\times 3}$ satisfies that
\begin{align}
Y_{\gamma_1, \gamma_2, \gamma_3}^{\circledS \mr{log}} := Y_{\gamma_1}^{\circledS \mr{log}} \times_{X^{\circledS \mr{log}}} Y_{\gamma_2}^{\circledS \mr{log}}  \times_{X^{\circledS \mr{log}}} Y_{\gamma_3}^{\circledS \mr{log}}  \neq \emptyset,
\end{align}
 then 
there exists uniquely an element 
\begin{align}
\nu_{\gamma_1, \gamma_2, \gamma_3}^\ddagger  \in & \Gamma ((Y_{\gamma_1, \gamma_2, \gamma_3})_b, (\mcT^\mcD_{X^{\circledS \bigstar \text{-}\mr{log}}/S^{\circledS \mr{log}}} \otimes_{f_b^{-1}(\mcO_{S^\circledS})} f_b^{-1}(\mcI))_b) \\
 \big( & = \Gamma ((Y_{\gamma_1, \gamma_2, \gamma_3})_t, (\mcT^{\underline{\mcD}}_{\underline{X}^{\circledS \bigstar \text{-}\mr{log}}/S_t^{\mr{log}}} \otimes_{\underline{f}_t^{-1}(\mcO_{S_t})} \underline{f}^{-1}_t(\mcI))_b)\big) \notag
\end{align}
 such that
\begin{align}
\nu_{\gamma_1, \gamma_2}^{\circledS \mr{log}} \circ \nu_{\gamma_2, \gamma_3}^{\circledS \mr{log}} \circ \nu_{\gamma_3, \gamma_1}^{\circledS \mr{log}} = \mr{id}_{Y^{\circledS \mr{log}}_{\gamma_1, \gamma_2, \gamma_3}} \boxplus^\ddagger \nu_{\gamma_1, \gamma_2, \gamma_3}^\ddagger.
\end{align}
The collection of elements $\{\nu_{\gamma_1, \gamma_2, \gamma_3}^\ddagger  \}_{\gamma_1, \gamma_2, \gamma_3}$ determines an element 
\begin{align}
\nu^\ddagger \in H^2 (\underline{X}_b, (\mcT^{\underline{\mcD}}_{\underline{X}^{\circledS \bigstar \text{-}\mr{log}}/S_t^{\mr{log}}} \otimes_{\underline{f}_t^{-1}(\mcO_{S_t})} \underline{f}^{-1}_t(\mcI))_b).
\end{align}
But, since $S^\circledS$ is affine and $\mr{dim}(\underline{X}_b/S_t) =1$,
  we have 
  \begin{align}
   H^2 (\underline{X}_b, \mcT^{\underline{\mcD}}_{\underline{X}^{\circledS \bigstar \text{-}\mr{log}}/S_t^{\mr{log}}} \otimes_{\underline{f}_b^{-1}(\mcO_{S^\circledS})} \underline{f}^{-1}_b(\mcI)) = 0
   \end{align}
    (in particular, $\nu^\ddagger  =0$).
Thus, after possibly replacing 
$\mbU_I$ with its  refinement and 
replacing each $\nu_{\gamma_1, \gamma_2}^{\circledS \mr{log}}$ with a suitable  isomorphism $\widetilde{\mbY}_{\gamma_1}^{\circledS \mr{log}} |_{Y_{\gamma_1, \gamma_2}^{\circledS \mr{log}}}   \isom \widetilde{\mbY}_{\gamma_2}^{\circledS \mr{log}}  |_{Y_{\gamma_1, \gamma_2}^{\circledS \mr{log}}}$,
the log superschemes $\{ \widetilde{Y}_\gamma^{\circledS \mr{log}} \}_{\gamma \in I}$
 may be glued together to a log supersmooth superstack $\widetilde{X}^{\circledS \mr{log}}$ over $\widetilde{S}^{\circledS \mr{log}}$.
 For each $i \in \{ 1, \cdots, r\}$,  the morphisms $\{ (\sigma_i^{\widetilde{U}_\gamma} \times \mr{id}) \circ (\widetilde{\eta}_\gamma^\circledS )^{-1} : \Sigma_i^{\widetilde{U}_\gamma} \times_{\widetilde{S}_b} \mbA^{0|1}_{\widetilde{S}^\circledS} \migi \widetilde{Y}^{\circledS}_\gamma \}_{\gamma \in I}$ may be glued together to
  a closed immersion $\widetilde{\sigma}^\circledS_i :  \mbA^{0|1}_{\widetilde{S}^\circledS} \migi \widetilde{X}^\circledS$ over $\widetilde{S}^\circledS$ extending   $\sigma_i^\circledS$, for which 
  the collection of data $(\widetilde{X}^{\circledS \mr{log}}/\widetilde{S}^{\circledS \mr{log}}, \{ \widetilde{\sigma}_i^\circledS \}_{i=1}^r)$ forms an $r$-pointed log twisted $(1|1)$-curve of genus $g$ over $\widetilde{S}^{\circledS \mr{log}}$.
   Moreover, $\{ \widetilde{\mcD}_\gamma \}_{\gamma \in I}$
  may be glued together to a superconformal structure  $\widetilde{\mcD}$ on 
 $\widetilde{X}^{\circledS \bigstar \text{-} \mr{log}}/\widetilde{S}^{\circledS \mr{log}}$ extending  $\mcD$.
The  collection of data 
\begin{align}
{^{\S_1} \widetilde{\mfX}}^{\circledS \bigstar} := (\widetilde{X}^{\circledS \mr{log}}/\widetilde{S}^{\circledS \mr{log}}, \{ [\widetilde{\sigma}^\circledS_i] \}_{i=1}^r, \widetilde{\mcD})
\end{align}
 forms a stable log twisted $\text{SUSY}_1$ curve  of type $(g,r, \lambda)$ over $\widetilde{S}^{\circledS \mr{log}}$ which restricts to ${^{\S_1} \mfX}^{\circledS \bigstar}$.
 Consequently, $\mr{Def}_{\widetilde{S}^{\circledS \mr{log}}} ({^{\S_1} \mfX}^{\circledS \mr{log}})$ is nonempty.

 
Also, by considering the above discussion and a usual discussion in deformation theory,
$\mr{Def}_{\widetilde{S}^{\circledS \mr{log}}} ({^{\S_1} \mfX}^{\circledS \mr{log}})$ admits a structure of affine space modeled on 
$H^1 (\underline{X}_b, (\mcT^{\underline{\mcD}}_{\underline{X}^{\circledS \bigstar \text{-}\mr{log}}/S_t^{\mr{log}}} \otimes_{\underline{f}_t^{-1}(\mcO_{S_t})} \underline{f}^{-1}_t(\mcI))_b)$ as described in (\ref{E003}).


Finally,  the above argument implies the remaining portion of the proposition.
Indeed, the group of superconformal automorphisms of  an arbitrary ${^{\S_1} \widetilde{\mfX}}^{\circledS \mr{log}} \in \mr{Def}_{\widetilde{S}^{\circledS \mr{log}}} ({^{\S_1} \mfX}^{\circledS \mr{log}})$ 
 is canonically isomorphic to 
 $H^0 (\underline{X}_b, (\mcT^{\underline{\mcD}}_{\underline{X}^{\circledS \bigstar \text{-}\mr{log}}/S_t^{\mr{log}}} \otimes_{\underline{f}_t^{-1}(\mcO_{S_t})} \underline{f}^{-1}_t(\mcI))_b)$.
If $(\mcL^\bigstar, \eta^\bigstar)$ denotes the spin structure on $(\underline{X}_t/S_t , \{[ \underline{\sigma}_i ]\}_{i=1}^r)$ corresponding to $\underline{\mcD}$ (cf. Proposition \ref{P66}), then we obtain  a sequence of isomorphisms
\begin{align} \label{e896}
& \ (\mcT^{\underline{\mcD}}_{\underline{X}^{\circledS \bigstar \text{-}\mr{log}}/S_t^{\mr{log}}} \otimes_{\underline{f}_b^{-1}(\mcO_{S_t})} \underline{f}^{-1}_t(\mcI))_b  \\
 \isom & \ 
  (\underline{\mcD}^{\otimes 2} \otimes_{\underline{f}_t^{-1}(\mcO_{S_t})} \underline{f}^{-1}_t(\mcI))_b \notag  \\
 \isom & \   ((\mcO_{\underline{X}^\circledS}\otimes \mcL^{\bigstar\vee})^{\otimes 2} \otimes_{(\underline{f}_t)^{-1}(\mcO_{S_t})} \underline{f}_t^{-1}(\mcI))_b \notag  \\
  \isom & \ (\mcT_{\underline{X}^{\bigstar \text{-}\mr{log}}/S_t^\mr{log}} \otimes_{\underline{f}_t^*(\mcO_{S_t})}\underline{f}_t^*(\mcI_b)) \oplus (\mcL^{\bigstar \vee} \otimes_{\underline{f}_t^{*}(\mcO_{S_t})}\underline{f}_t^{*}(\mcI_f)), \notag
\end{align}
where the first isomorphism follows from Proposition \ref{p010} and  the second isomorphism follows from (\ref{E30}).
Here, note that the natural morphisms
\begin{align}
H^0(\underline{X}_b, \mcT_{\underline{X}^{\bigstar \text{-}\mr{log}}/S_t^\mr{log}})  \otimes_{\Gamma (S_t, \mcO_{S_t})} \Gamma (S_t, \mcI_b) \migi H^0(\underline{X}_b, (\mcT_{\underline{X}^{\bigstar \text{-}\mr{log}}/S_t^\mr{log}} \otimes_{\underline{f}_t^*(\mcO_{S_t})}\underline{f}_t^*(\mcI_b)))
\end{align}
and 
\begin{align}
H^0 (\underline{X}_b, \mcL^{\bigstar \vee})  \otimes_{\Gamma (S_t, \mcO_{S_t})} \Gamma (S_t, \mcI_f)  \migi H^0 (\underline{X}_b, (\mcL^{\bigstar \vee} \otimes_{\underline{f}_t^{*}(\mcO_{S_t})}\underline{f}_t^{*}(\mcI_f)))
\end{align}
are surjective.
On the other hand, the fact that $\mr{deg} (\mcT_{\underline{X}^{\bigstar \text{-}\mr{log}}/S_t^\mr{log}} ) <0$ and $\mr{deg} (\mcL^{\bigstar \vee})<0$ implies the equalities  $H^0(\underline{X}_b, \mcT_{\underline{X}^{\bigstar \text{-}\mr{log}}/S_t^\mr{log}}) = H^0 (\underline{X}_b, \mcL^{\bigstar \vee}) =0$.
Hence, we have 
\begin{align}
& \  H^0 (\underline{X}_b, (\mcT^{\underline{\mcD}}_{\underline{X}^{\circledS \bigstar \text{-}\mr{log}}/S_t^{\mr{log}}} \otimes_{\underline{f}_t^{-1}(\mcO_{S_t})} \underline{f}^{-1}_t(\mcI))_b) \\
= 
& \  H^0(\underline{X}_b, (\mcT_{\underline{X}^{\bigstar \text{-}\mr{log}}/S_t^\mr{log}} \otimes_{\underline{f}_t^*(\mcO_{S_t})}\underline{f}_t^*(\mcI_b))) \oplus H^0 (\underline{X}_b, (\mcL^{\bigstar \vee} \otimes_{\underline{f}_t^{*}(\mcO_{S_t})}\underline{f}_t^{*}(\mcI_f))) \notag \\
= & \  0. \notag 
\end{align}
This implies that  there is no nontrivial superconformal automorphism of any ${^{\S_1} \widetilde{\mfX}}^{\circledS \bigstar} \in \mr{Def}_{\widetilde{S}^{\circledS \mr{log}}} ({^{\S_1} \mfX}^{\circledS \mr{log}})$.
This completes the proof of Proposition \ref{P001}.
\end{proof}

\vspace{5mm}

\subsection{Deformations of morphisms} \label{S52}
\leavevmode\\
\vspace{-4mm}

We keep the notation in the previous subsection.
Moreover, 
let 
$T^{\circledS \mr{log}}$ be an fs   log superscheme which is 
  log supersmooth over $S_0$ (of some relative superdimension) and  
${^{\S_1} \mfY}^{\circledS \bigstar} := (f'^{\circledS \mr{log}} : Y^{\circledS \mr{log}} \migi T^{\circledS \mr{log}}, \{ [\sigma'_i ] \}_{i=1}^r, \mcD')$  a stable log twisted $\text{SUSY}_1$ curve of type $(g,r, \lambda)$ over $T^{\circledS \mr{log}}$ such that $\mcK\mcS({^{\S_1} \mfY}^{\circledS \bigstar})$ is an isomorphism.
Suppose that we are given a morphism $s^{\circledS \mr{log}}  : S^{\circledS \mr{log}} \migi T^{\circledS \mr{log}}$ of log superschemes  via  which the base-change of ${^{\S_1} \mfY}^{\circledS \bigstar}$ is isomorphic to ${^{\S_1} \mfX}^{\circledS \bigstar}$ over $S^{\circledS \mr{log}}$.

The following proposition is immediately verified from the various definitions involved, including the affine structures described in Proposition \ref{p0404} (ii) and Proposition \ref{P001}.

\vspace{3mm}
\bpr [cf. ~\cite{LR}, Lemma 2.5] \label{p03} \leavevmode\\
 \ \ \
Suppose that  $\widetilde{S}^{\circledS}$ is affine and that $\mcN_{\widetilde{S}^\circledS} \mcI =0$.
Denote by $\mr{KS}({^{\S_1} \mfY}^{\circledS \bigstar}; \mcI)$ the composite isomorphism
\begin{align} \label{E002}
\Gamma (S_b, (s^{\circledS *}(\mcT_{T^{\circledS \mr{log}}/S_0}) \otimes \mcI)_b) & \isom \Gamma (S_b, (s^{\circledS *}(\mbR^1 f'_{b*} (\mcT^\mcD_{Y^{\circledS \bigstar \text{-} \mr{log}}/T^{\circledS \mr{log}}})) \otimes \mcI)_b) \\
& \isom  \Gamma (S_b, (\mbR^1 f_{b*} (\mcT^\mcD_{\underline{X}^{\circledS \bigstar \text{-} \mr{log}}/S_t^{\circledS \mr{log}}})  \otimes \mcI)_b) \notag \\
& \isom  H^1 (\underline{X}_b, (\mcT^\mcD_{\underline{X}^{\circledS \bigstar \text{-} \mr{log}}/S_t^{\circledS \mr{log}}} \otimes_{\underline{f}_b^{*}(\mcO_{S_t})} \underline{f}^*_b(\mcI))_b), \notag
\end{align}
where the first  isomorphism arises from  $\mcK\mcS({^{\S_1} \mfY}^{\circledS \bigstar})$ and both the second   and third isomorphisms arise from Proposition \ref{Pt66}.
Consider the map of sets
\begin{align}
s^\circledast : \Gamma (\widetilde{S}_b, \mcD ef_{\widetilde{S}^{\circledS \mr{log}}}(s^{\circledS \mr{log}})) \migi \mr{Def}_{\widetilde{S}^{\circledS \mr{log}}} ({^{\S_1} \mfX}^{\circledS \bigstar})
\end{align}
given by pulling-back ${^{\S_1} \mfY}^{\circledS \bigstar}$.
 Then, this map satisfies the equality
 \begin{align} \label{e997}
s^\circledast (\widetilde{s}^{\circledS \mr{log}}  \boxplus^\dagger  \zeta)
  =s^\circledast (\widetilde{s}^{\circledS \mr{log}}) \boxplus^\ddagger \mr{K S}({^{\S_1} \mfY}^{\circledS \bigstar}; \mcI) (\zeta).
 \end{align}
for any $\widetilde{s}^{\circledS \mr{log}} \in \Gamma (\widetilde{S}_b, \mcD ef_{\widetilde{S}^{\circledS \mr{log}}}(s^{\circledS \mr{log}}))$ and $\zeta \in \Gamma (S_b, (s^{\circledS *}(\mcT_{T^{\circledS \mr{log}}/S_0}) \otimes \mcI)_b)$.
 In particular, (since $\mr{KS}({^\S \mfY}^{\circledS \bigstar}; \mcI)$ is an isomorphism)
 $s^\circledast$ is bijective, and hence, $\Gamma (\widetilde{S}_b, \mcD ef_{\widetilde{S}^{\circledS \mr{log}}}(s^{\circledS \mr{log}}))$ is nonempty.
   \epr
\vspace{3mm}

\vspace{3mm}
\bco[cf. ~\cite{LR}, Theorem 2.7] \label{p04} \leavevmode\\
 \ \ \ Suppose that we are given a stable log twisted $\text{SUSY}_1$ curve ${^{\S_1} \widetilde{\mfX}}^{\circledS \bigstar}$ over $\widetilde{S}^{\circledS \mr{log}}$ extending  ${^{\S_1} \mfX}^{\circledS \bigstar}$.
 Then, there exists uniquely an extension  $\widetilde{s}^{\circledS \mr{log}} : \widetilde{S}^{\circledS \mr{log}} \migi T^{\circledS \mr{log}}$ of  $s^{\circledS \mr{log}}$ via which the base-change of ${^{\S_1} \mfY}^{\circledS \mr{log}}$ is isomorphic to   ${^{\S_1} \widetilde{\mfX}}^{\circledS \bigstar}$.
 \eco
\begin{proof}
The assertion may be directly  proved by applying inductively Proposition \ref{p03} to the case where
the pair 
$(S^{\circledS \mr{log}}, \widetilde{S}^{\circledS \mr{log}})$ is taken to be $(\widetilde{S}^{\circledS \mr{log}}_n, \widetilde{S}^{\circledS \mr{log}}_{n+1})$ (where, for each $n \geq 0$, we shall denote by  $\widetilde{S}^\circledS_n$ the strict closed subsuperscheme of $\widetilde{S}^{\circledS \mr{log}}$ determined by $\mcN_{\widetilde{S}^\circledS}^{n+1}$).
\end{proof}

\vspace{5mm}

\subsection{Canonical liftings over complete versal families} \label{S53}
\leavevmode\\
\vspace{-4mm}

\bpr[cf. ~\cite{LR}, Theorem 2.8] \label{p01} \leavevmode\\
 \ \ \ 
Let $\underline{T}^\mr{log}$ be an affine  log smooth  scheme in $\mfS \mfc \mfh_{/S_0}^\mr{log}$,  
and let  ${^{\S_1} \underline{\mfY}}^{\circledS \mr{log}} := (\underline{f}^{\circledS \mr{log}} : \underline{Y}^{\circledS \mr{log}} \migi \underline{T}^\mr{log}, \{[ \underline{\sigma}_{i}] \}_{i=1}^r, \underline{\mcD})$ be   a stable log twisted  $\text{SUSY}_1$ curve of type $(g,r, \lambda)$ over $\underline{T}^\mr{log}$ such that the classical  Kodaira-Spencer map 
\begin{align}
\mcK \mcS (\underline{Y}_b^{\bigstar \text{-}\mr{log}}/\underline{T}^{\mr{log}}) : \mcT_{\underline{T}^{\mr{log}}/S_0} \migi \mbR^1 \underline{f}_{b*} (\mcT_{\underline{Y}_b^{\bigstar \text{-}\mr{log}}/\underline{T}^{\mr{log}}})
\end{align}
 of $\underline{Y}_b^{\bigstar \text{-}\mr{log}} \migi \underline{T}^{\mr{log}}$ (cf. (\ref{E1111}))  is an isomorphism.
  Let us write
  \begin{align} \label{e678}
  \langle \underline{T} \rangle^{\circledS \mr{log}} :=\langle \underline{T}, \mbR^1 \underline{f}_{b *} ((\mcT^{\underline{\mcD}}_{\underline{Y}^{\circledS \bigstar \text{-} \mr{log}}/\underline{T}^{\mr{log}}})_f)^\vee \rangle^{\circledS \mr{log}}.
  \end{align}
Then, there exists  a stable log twisted $\text{SUSY}_1$ curve ${^{\S_1} \mfY}_{\dagger}^{\circledS \bigstar}$ of type $(g,r, \lambda)$ over  $\langle \underline{T} \rangle^{\circledS \mr{log}}$ which restricts to  ${^{\S_1} \underline{\mfY}}^{\circledS \bigstar}$ via $\tau_{\langle \underline{T} \rangle}^{\circledS \mr{log}} : \underline{T}^\mr{log} \migi \langle \underline{T} \rangle^{\circledS \mr{log}}$ and whose Kodaira-Spencer map $\mcK \mcS ({^{\S_1} \mfY}_\dagger^{\circledS \bigstar})$ is an isomorphism.
Moreover, such a stalbe log twisted $\text{SUSY}_1$ curve  is unique in the following sense:
if ${^{\S_1} \mfY}_\dagger^{\circledS \bigstar}$ and ${^{\S_1} \mfY}_\ddagger^{\circledS \bigstar}$
are stable log twisted $\text{SUSY}_1$ curves of type $(g,r, \lambda)$ over $\langle \underline{T} \rangle^{\circledS \mr{log}}$ which restrict to  ${^{\S_1} \underline{\mfY}}^{\circledS \bigstar}$ and whose Kodaira-Spencer maps are isomorphism, then there exists uniquely a superconformal isomorphism ${^{\S_1} \mfY}_\dagger^{\circledS \bigstar}  \isom {^{\S_1} \mfY}_\ddagger^{\circledS \bigstar}$ over $\langle \underline{T} \rangle^{\circledS \mr{log}}$ which restricts to the identity morphism of ${^{\S_1} \underline{\mfY}}^{\circledS \bigstar}$.
 \epr
\begin{proof}
The uniqueness portion  follows from the uniqueness assertion of Corollary \ref{p04}.
We shall prove the existence portion.
For each nonnegative  integer $n$, we shall write  $\langle \underline{T} \rangle_{n}^{\circledS \mr{log}}$ for the strict  closed subsuperscheme of $\langle \underline{T} \rangle^{\circledS \mr{log}}$ corresponding to the ideal $\mcN_{\langle T \rangle^\circledS}^{n+1}$.
Since $\langle \underline{T} \rangle_1^{\circledS \mr{log}}$  is simply $(\underline{T}, \mcO_{\underline{T}} \oplus \mbR^1 \underline{f}_{b *} ((\mcT^{\underline{\mcD}}_{\underline{Y}^{\circledS \bigstar \text{-} \mr{log}}/\underline{T}^{\mr{log}}})_f)^\vee)$, 
we obtain the trivial deformation ${^{\S_1} \mfY}_{1, \mr{triv}}^{\circledS \bigstar}$ of ${^{\S_1} \underline{\mfY}}^{\circledS \bigstar}$ over $\langle \underline{T} \rangle_{1}^{\circledS \mr{log}}$ by pulling-back via  the projection $\langle \underline{T} \rangle_{1}^{\circledS \mr{log}}  \migi \underline{T}^\mr{log}$.
By applying Proposition \ref{p03} and considering the point of  the affine space $\mr{Def}_{\langle \underline{T} \rangle_{1}^{\circledS \mr{log}}} ({^{\S_1} \underline{\mfY}}^{\circledS \bigstar})$ representing ${^{\S_1} \mfY}_{1, \mr{triv}}^{\circledS \bigstar}$  as its origin, we have a canonical composite  bijection
\begin{align} \label{e3402}
& \hspace{8mm}  \mr{Def}_{\langle \underline{T} \rangle_{1}^{\circledS \mr{log}}} ({^{\S_1} \underline{\mfY}}^{\circledS \bigstar})  \\
& \isom H^1 (\underline{Y}_b, (\mcT^{\underline{\mcD}}_{\underline{Y}^{\circledS \bigstar\text{-}\mr{log}}/\underline{T}^\mr{log}}\otimes \underline{f}_b^{*}(   \mbR^1 \underline{f}_{b *} ((\mcT^{\underline{\mcD}}_{\underline{Y}^{\circledS \bigstar \text{-} \mr{log}}/\underline{T}^{\mr{log}}})_f)^\vee))_b) \notag   \\
& \isom H^1 (\underline{Y}_b, (\mcT^{\underline{\mcD}}_{\underline{Y}^{\circledS \bigstar\text{-}\mr{log}}/\underline{T}^\mr{log}})_f \otimes  \underline{f}_b^{*}(   \mbR^1 \underline{f}_{b *} ((\mcT^{\underline{\mcD}}_{\underline{Y}^{\circledS \bigstar \text{-} \mr{log}}/\underline{T}^{\mr{log}}})_f)^\vee))\notag \\
& \isom 
H^1  (\underline{Y}_b,  (\mcT^{\underline{\mcD}}_{\underline{Y}^{\circledS \bigstar\text{-}\mr{log}}/\underline{T}^\mr{log}})_f) \otimes_{H^0(\underline{T}, \mcO_{\underline{T}})} H^1  (\underline{Y}_b,  (\mcT^{\underline{\mcD}}_{\underline{Y}^{\circledS \bigstar\text{-}\mr{log}}/\underline{T}^\mr{log}})_f)^\vee
\notag  \\
& \isom \mr{End}_{H^0(\underline{T}, \mcO_{\underline{T}})}(H^1  (\underline{Y}_b,  (\mcT^{\underline{\mcD}}_{\underline{Y}^{\circledS \bigstar\text{-}\mr{log}}/\underline{T}^\mr{log}})_f)), \notag
\end{align}
where the third bijection follows from Proposition \ref{Pt66}.
If we write ${^{\S_1} \mfY}_{1}^{\circledS \bigstar}$ for the  stable log twisted $\text{SUSY}_1$ curve corresponding to  
\begin{align}
\mr{id}_{H^1  (\underline{Y}_b,  (\mcT^{\underline{\mcD}}_{\underline{Y}^{\circledS \bigstar\text{-}\mr{log}}/\underline{T}^\mr{log}})_f)} \in  \mr{End}_{H^0(\underline{T}, \mcO_{\underline{T}})}(H^1  (\underline{Y}_b,  (\mcT^{\underline{\mcD}}_{\underline{Y}^{\circledS \bigstar\text{-}\mr{log}}/\underline{T}^\mr{log}})_f))
\end{align}
 via (\ref{e3402}), then
its Kodaira-Spencer map turns out to be  an isomorphism.
By Proposition \ref{P001}, ${^{\S_1} \mfY}_{1}^{\circledS \bigstar}$  may be deformed to a  stable log twisted  $\text{SUSY}_1$  curve  
\begin{align}
{^{\S_1} \mfY}_\dagger^{\circledS \bigstar} := (f^{\circledS \mr{log}}_\dagger : Y_\dagger^{\circledS \mr{log}} \migi \langle \underline{T} \rangle^{\circledS \mr{log}}, \{ [\sigma_{\dagger,  i}^\circledS] \}_{i=1}^r, \mcD_\dagger)
\end{align}
 of type $(g,r,\lambda)$ over $\langle \underline{T} \rangle^{\circledS \mr{log}}$.

We shall prove  that the Kodaira-Spencer map $\mcK \mcS ({^{\S_1} \mfY}_\dagger^{\circledS \bigstar})$ is an isomorphism.
To this end,  it suffices to prove that 
its restriction   along the reduced space $\underline{T}$ is an isomorphism.
(Indeed,  by Proposition \ref{Pt66}, $\mbR^1 f^\circledS_{\dagger *}(\mcT^{\mcD_\dagger}_{Y_\dagger^{\circledS \bigstar \text{-} \mr{log}}/\langle \underline{T} \rangle^{\circledS \mr{log}}})$ and $\mcT_{\langle \underline{T} \rangle^{\circledS \mr{log}}/S_0}$ are  locally free of the same rank).
By the definition of $\langle \underline{T} \rangle^{\circledS \mr{log}}$,  the pull-back of $\mcT_{\langle \underline{T} \rangle^{\circledS \mr{log}}/S_0}$ via $\tau^\circledS_{\langle \underline{T} \rangle} : \underline{T} \migi \langle \underline{T} \rangle$
admits a canonical   isomorphism
\begin{align} \label{e3333}
\tau_{\langle \underline{T} \rangle}^{\circledS *}(\mcT_{\langle \underline{T} \rangle^{\circledS \mr{log}}/S_0}) \isom \mcT_{\underline{T}^\mr{log}/S_0} \oplus \mbR^1 \underline{f}_{b *} ((\mcT^{\underline{\mcD}}_{\underline{Y}^{\circledS \bigstar \text{-} \mr{log}}/\underline{T}^{\circledS \mr{log}}})_f).
\end{align}
On the other hand, the pull-back of $\mbR^1 f^\circledS_{\dagger *}(\mcT^{\mcD_\dagger}_{Y_\dagger^{\circledS \bigstar \text{-} \mr{log}}/\langle \underline{T} \rangle^{\circledS \mr{log}}})$ 
admits a canonical composite isomorphism
\begin{align} \label{e3323}
& \ \tau_{\langle \underline{T} \rangle}^{\circledS *}(\mbR^1 f_{\dagger   *}^\circledS(\mcT^{\mcD_\dagger}_{Y_\dagger^{\circledS \bigstar \text{-} \mr{log}}/\langle \underline{T} \rangle^{\circledS \mr{log}}})) \\
 \isom & \ \mbR^1 \underline{f}_{*}^\circledS(\mcT^{\underline{\mcD}}_{\underline{Y}^{\circledS \bigstar \text{-} \mr{log}}/\underline{T}^{\circledS \mr{log}}}) \notag  \\
\isom & \  \mbR^1 \underline{f}_{b*}((\mcT^{\underline{\mcD}}_{\underline{Y}^{\circledS \bigstar \text{-} \mr{log}}/\underline{T}^{\circledS \mr{log}}})_b)  \oplus  \mbR^1 \underline{f}_{b*}((\mcT^{\underline{\mcD}}_{\underline{Y}^{\circledS \bigstar \text{-} \mr{log}}/\underline{T}^{\circledS \mr{log}}})_f)  \notag  \\
    \isom & \ 
\mbR^1 \underline{f}_{b*}(\mcT_{\underline{Y}_b^{\bigstar \text{-} \mr{log}}/\underline{T}^{\mr{log}}})
     \oplus   \mbR^1 \underline{f}_{*}((\mcT^{\underline{\mcD}}_{\underline{Y}^{\circledS \bigstar \text{-} \mr{log}}/\underline{T}^{\circledS \mr{log}}})_f)  \notag 
\end{align}
 where the first isomorphism follows from  Proposition \ref{Pt66} and the third isomorphism  follows from (\ref{E006}).
One may verifies immediately  from the various  definitions involved that  
the pull-back $\tau_{\langle \underline{T} \rangle}^{\circledS *}(\mcK \mcS ({^{\S_1} \mfY}_\dagger^{\circledS \bigstar}))$ of $\mcK \mcS ({^{\S_1} \mfY}_\dagger^{\circledS \bigstar})$
makes the square diagram
\begin{align}
\begin{CD}
\tau_{\langle \underline{T} \rangle}^{\circledS *}(\mcT_{\langle \underline{T} \rangle^{\circledS \mr{log}}/S_0}) @> (\ref{e3333}) > \sim > 
\mcT_{\underline{T}^\mr{log}/S_0} \oplus \mbR^1 f_{b *} ((\mcT^{\underline{\mcD}}_{\underline{Y}^{\circledS \bigstar \text{-} \mr{log}}/\underline{T}^{\circledS \mr{log}}})_f) 
\\
@V \tau_{\langle \underline{T} \rangle}^{\circledS *}(\mcK \mcS ({^{\S_1} \mfY}_\dagger^{\circledS \bigstar}))  V V @VV \mcK \mcS (\underline{Y}_b^{\bigstar \text{-}\mr{log}}/\underline{T}^\mr{log}) \oplus \mr{id} V
\\
\tau_{\langle \underline{T} \rangle}^{\circledS *}(\mbR^1 f_{\dagger   *}^\circledS(\mcT^{\mcD_\dagger}_{Y_\dagger^{\circledS \bigstar \text{-} \mr{log}}/\langle \underline{T} \rangle^{\circledS \mr{log}}})) 
@> \sim > (\ref{e3323}) > \mbR^1 \underline{f}_{b*}(\mcT_{\underline{Y}_b^{\bigstar \text{-} \mr{log}}/\underline{T}^{\mr{log}}})
     \oplus   \mbR^1 \underline{f}_{b *}((\mcT^{\underline{\mcD}}_{\underline{Y}^{\circledS \bigstar \text{-} \mr{log}}/\underline{T}^{\circledS \mr{log}}})_f)
\end{CD}
\end{align}
commute.
Hence, since we have assumed  that $\mcK \mcS (\underline{Y}^{\bigstar \text{-}\mr{log}}/\underline{T}^\mr{log})$ is an isomorphism, $ \tau_{\langle \underline{T} \rangle}^{\circledS *}(\mcK \mcS ({^{\S_1} \mfY}_\dagger^{\circledS \bigstar}))$ is an isomorphism, as desired.
 This completes the proof of Proposition \ref{p01}.
\end{proof}

\vspace{5mm}
\subsection{The proof of Theorem A} \label{S54}
\leavevmode\\
\vspace{-4mm}


In this final section, we shall prove Theorem A, the main result of the present paper.
Since  $({^{\S_1} \overline{\mfM}}_{g,r, \lambda})_t$
 is a smooth Deligne-Mumford stack over $S_0$ (cf. Proposition \ref{P66}), there exists an isomorphism $[\underline{R} \rightrightarrows \underline{U}] \isom ({^{\S_1} \overline{\mfM}}_{g,r, \lambda})_t$ over $S_0$ for some  groupoid $\underline{R} \rightrightarrows \underline{U} := (\underline{U}, \underline{R}, \underline{s}, \underline{t}, \underline{c})$ in $\mfS \mfc \mfh_{/S_0}$ such that both  $\underline{U}$ and $\underline{R}$ are smooth affine schemes  over $S_0$ of relative dimension $3g-3+r$, and both $\underline{s}$ and $\underline{t}$ are \'{e}tale.
Denote by $\pi_{\underline{U}} : \underline{U} \migi ({^{\S_1} \overline{\mfM}}_{g,r, \lambda})_t$ the natural projection (hence $\pi_{\underline{R}} := \pi_{\underline{U}}  \circ \underline{s} = \pi_{\underline{U}} \circ \underline{t}$).
Write $\underline{U}^\mr{log}$ (resp., $\underline{R}^\mr{log}$) for the log  scheme defined to be $\underline{U}$ (resp., $\underline{R}$) equipped with the log structure pulled-back from $({^{\S_1} \overline{\mfM}}_{g,r, \lambda})_t^\mr{log}$.
In particular, $\pi_{\underline{U}}$ (resp., $\pi_{\underline{R}}$) extends to a morphisms $\pi_{\underline{U}}^\mr{log} : \underline{U}^\mr{log} \migi ({^{\S_1} \overline{\mfM}}_{g,r, \lambda})_t^\mr{log}$ (resp., $\pi_{\underline{R}}^\mr{log} : \underline{R}^\mr{log} \migi ({^{\S_1} \overline{\mfM}}_{g,r, \lambda})_t^\mr{log}$)  of log  stacks. 
Moreover, $\underline{s}, \underline{t} : \underline{R} \migi \underline{U}$ extend to morphisms
$\underline{s}^\mr{log}, \underline{t}^\mr{log} : \underline{R}^\mr{log} \migi \underline{U}^\mr{log}$,
and $\underline{c}$ extends to 
 a morphism $\underline{c}^\mr{log} : \underline{R}^\mr{log} \times_{\underline{s}^\mr{log}, \underline{U}^\mr{log}, \underline{t}^\mr{log}} \underline{R}^\mr{log} \migi \underline{R}^\mr{log}$.

Let us write  
\begin{align}
& {^{\S_1} \underline{\mfY}}^{\circledS \bigstar} := (\underline{Y}^{\circledS \mr{log}}/\underline{U}^\mr{log}, \{[\sigma^\circledS_{\underline{Y}, i}] \}_{i=1}^r, \mcD_{\underline{Y}}) \\
(\text{resp.,} \  & {^{\S_1} \underline{\mfX}}^{\circledS \bigstar} := (\underline{X}^{\circledS \mr{log}}/\underline{R}^\mr{log}, \{[\sigma^\circledS_{\underline{X}, i}] \}_{i=1}^r, \mcD_{\underline{X}})) \notag
\end{align}
for the stable   log  twisted $\text{SUSY}_1$ curve over $\underline{U}^\mr{log}$ (resp., $\underline{R}^\mr{log}$) classified by 
$\pi^\mr{log}_{\underline{U}}$ (resp., $\pi^\mr{log}_{\underline{R}}$).
The base-change of ${^{\S_1} \underline{\mfY}}^{\circledS \bigstar}$ via  $\underline{s}^\mr{log}$ and $\underline{t}^\mr{log}$ respectively  are, by definition, isomorphic to ${^{\S_1} \underline{\mfX}}^{\circledS \bigstar}$.
The Kodaira-Spencer morphisms $\mcK \mcS ({^{\S_1} \underline{\mfY}}^{\circledS \bigstar})$, $\mcK \mcS ({^{\S_1} \underline{\mfX}}^{\circledS \bigstar})$ are isomorphisms.
Here, let us define $\langle \underline{U} \rangle^{\circledS \mr{log}}$ (resp., $\langle \underline{R} \rangle^{\circledS \mr{log}}$)  to be the log superscheme obtained from $\underline{U}^{\mr{log}}$ (resp., $\underline{R}^{\mr{log}}$) as in (\ref{e678}), which is  split and log supersmooth over $S_0$ of relative superdimension $3g-3+r | 2g-2+\frac{r}{2}$ (by Proposition \ref{Pt66}).
It follows from Proposition \ref{p01} 
that ${^{\S_1} \underline{\mfY}}^{\circledS \bigstar}$ (resp., ${^{\S_1} \underline{\mfX}}^{\circledS \bigstar}$) may be deformed to a stable log twisted  $\text{SUSY}_1$ curve   
\begin{align}
& {^{\S_1} \mfY}_\dagger^{\circledS \bigstar} := (f^{\circledS \mr{log}}_{Y, \dagger} : Y_\dagger^{\circledS \mr{log}}\migi \langle \underline{U} \rangle^{\circledS \mr{log}}, \{ [\sigma_{Y, \dagger,  i}^\circledS] \}_{i=1}^r, \mcD_{Y, \dagger}) \\
(\text{resp.,} \ & {^{\S_1} \mfX}_\dagger^{\circledS \bigstar} := (f^{\circledS \mr{log}}_{X, \dagger} : X_\dagger^{\circledS \mr{log}}\migi \langle \underline{R} \rangle^{\circledS \mr{log}}, \{ [\sigma_{X, \dagger,  i}^\circledS] \}_{i=1}^r, \mcD_{X, \dagger})) \notag
\end{align}
 over 
$\langle \underline{U} \rangle^{\circledS \mr{log}}$ (resp., $\langle \underline{R} \rangle^{\circledS \mr{log}}$) whose Kodaira-Spencer map is an isomorphism.
Hence, by  Corollary \ref{p04}, 
there exists   morphisms
$\langle \underline{s} \rangle^{\circledS \mr{log}}, \langle \underline{t} \rangle^{\circledS \mr{log}} : \langle \underline{R} \rangle^{\circledS \mr{log}} \migi \langle \underline{U} \rangle^{\circledS \mr{log}}$
via  which the base-changes of ${^{\S_1} \mfY}_{\dagger}^{\circledS \bigstar}$ are isomorphic to ${^{\S_1} \mfX}_\dagger^{\circledS \bigstar}$ and 
which make  the square  diagrams
\begin{align} \label{dg01}
\xymatrix{
 \underline{R}^\mr{log} \ar[r]^{\underline{s}^\mr{log}} \ar[d]_{\tau^{\circledS \mr{log}}_{ \langle \underline{R} \rangle}}&  \underline{U}^\mr{log}  \ar[d]^{\tau^{\circledS \mr{log}}_{ \langle \underline{U} \rangle}}  \\
 \langle \underline{R} \rangle^{\circledS \mr{log}} \ar[r]_{\langle \underline{s} \rangle^{\circledS \mr{log}}} & \langle \underline{U} \rangle^{\circledS \mr{log}}
} \hspace{10mm}
\xymatrix{
 \underline{R}^\mr{log} \ar[r]^{\underline{t}^\mr{log}} \ar[d]_{\tau^{\circledS \mr{log}}_{ \langle \underline{R} \rangle}}&  \underline{U}^\mr{log}  \ar[d]^{\tau^{\circledS \mr{log}}_{ \langle \underline{U} \rangle}}  \\
 \langle \underline{R} \rangle^{\circledS \mr{log}} \ar[r]_{\langle \underline{t} \rangle^{\circledS \mr{log}}} & \langle \underline{U} \rangle^{\circledS \mr{log}}
}
\end{align}
commute.
Moreover, we obtain a morphism 
\begin{align}
\langle \underline{c}\rangle^{\circledS \mr{log}} :  \langle \underline{R} \rangle^{\circledS \mr{log}} \times_{ \langle \underline{s} \rangle^{\circledS \mr{log}},  \langle \underline{U} \rangle^{\circledS \mr{log}},  \langle \underline{t} \rangle^{\circledS \mr{log}}}  \langle \underline{R} \rangle^{\circledS \mr{log}} \migi  \langle \underline{R} \rangle^{\circledS \mr{log}}
\end{align}
extending 
 the morphism $\underline{c}^{\circledS \mr{log}}$.
The uniqueness  assertion in Corollary \ref{p04} implies that  the collection of data 
\begin{align}
\langle \underline{R} \rangle^{\circledS} \rightrightarrows \langle \underline{U} \rangle^{\circledS} := (\langle \underline{U} \rangle^{\circledS}, \langle \underline{R} \rangle^{\circledS}, \langle \underline{s} \rangle^{\circledS}, \langle \underline{t} \rangle^{\circledS},  \langle \underline{c} \rangle^{\circledS})
\end{align} 
forms a groupoid  in $\mfS \mfc \mfh_{/S_0}^{\circledS}$.

For $\Box = \underline{s}$ or $\underline{d}$,
we shall denote by 
\begin{align}
d \langle \Box \rangle^\circledS : \mcT_{\langle \underline{R} \rangle^{\circledS \mr{log}}/S_0}  \migi  \langle \Box \rangle^{\circledS *} (\mcT_{\langle \underline{U} \rangle^{\circledS \mr{log}}/S_0})
\end{align}
 the differential  of $\langle \Box \rangle^\circledS$
   relative to $S_0$.
Then,  $d \langle \Box \rangle^\circledS$ is  an isomorphism  since 
the square diagram
\begin{align}
\xymatrix{
\mcT_{\langle \underline{R} \rangle^{\circledS \mr{log}}/S_0} \ar[r]^{d \langle \Box \rangle^\circledS }   \ar[d]^{\wr}_{\mcK \mcS ({^{\S_1} \mfX}_{\dagger}^{\circledS \bigstar})}& \langle \Box \rangle^{\circledS *} (\mcT_{\langle \underline{U} \rangle^{\circledS \mr{log}}/S_0})  \ar[d]_{\wr}^{\langle \Box \rangle^{\circledS *}(\mcK \mcS ({^{\S_1} \mfY}_{\dagger}^{\circledS \bigstar}))}
\\
\mbR^1 f^\circledS_{X, \dagger} (\mcT^{\mcD_{X, \dagger}}_{X_{\dagger}^{\circledS \bigstar \text{-} \mr{log}}/ \langle \underline{R} \rangle^{\circledS \mr{log}}}) \ar[r]^{\hspace{-5mm} \sim} & \langle \Box \rangle^{\circledS *} (\mbR^1 f^\circledS_{Y, \dagger} (\mcT^{\mcD_{Y, \dagger}}_{Y_{\dagger}^{\circledS \bigstar \text{-} \mr{log}}/ \langle \underline{U} \rangle^{\circledS \mr{log}}}))
}
\end{align}
is commute and cartesian (where  the lower horizontal arrow is isomorphism by  Proposition \ref{Pt66}).
Hence, for each $n \geq 0$, the morphism 
 $\mr{gr}^n_{\langle \underline{U} \rangle^\circledS} \migi \mr{gr}^n_{\langle \underline{R} \rangle^\circledS} $ induced by $d \langle \Box \rangle^\circledS$ is an isomorphism.
 It follows immediately that   $\langle \underline{s} \rangle_b, \langle \underline{t} \rangle_b  : \langle \underline{R} \rangle_b \migi \langle \underline{U} \rangle_b$ are \'{e}tale  (since $\underline{s}$ and $\underline{t}$ are \'{e}tale) and that two morphisms 
 \begin{align}
 (\langle \underline{s} \rangle^\circledS, \beta^\circledS_{\langle \underline{R} \rangle_b}), (\langle \underline{t} \rangle^\circledS, \beta^\circledS_{\langle \underline{R} \rangle_b}) : \langle \underline{R} \rangle^{\circledS} \migi \langle \underline{U}^\circledS \rangle \times_{\langle \underline{U} \rangle_b} \langle \underline{R} \rangle_b
 \end{align}
  are  isomorphisms.
Thus, both $\langle \underline{s} \rangle^\circledS$ and  $\langle \underline{t} \rangle^\circledS$ are super\'{e}tale.
 By Proposition \ref{p0207},  $[\langle \underline{R} \rangle^{\circledS} \rightrightarrows \langle \underline{U} \rangle^{\circledS}]$ forms  a supersmooth Deligne-Mumford superstack over $S_0$ of relative superdimension $3g-3+r | 2g-2+\frac{r}{2}$; it is superproper over $S_0$ since $({^{\S_1} \overline{\mfM}}_{g,r, \lambda})_t$ is proper  over $S_0$ in the classical sense. 

Moreover, the log structures of the various constituents  in $\langle \underline{R} \rangle^{\circledS} \rightrightarrows \langle \underline{U} \rangle^{\circledS}$
 gives rises to a log structure on 
 the superstack $[\langle \underline{R} \rangle^{\circledS} \rightrightarrows \langle \underline{U} \rangle^{\circledS}]$.
Let us write $[\langle \underline{R} \rangle^{\circledS} \rightrightarrows \langle \underline{U} \rangle^{\circledS}]^\mr{log}$ for the resulting log superstack
and write
$\pi^{\circledS \mr{log}}_{\langle \underline{U} \rangle} : \langle \underline{U} \rangle^{\circledS \mr{log}} \migi {^{\S_1} \overline{\mfM}}_{g,r,\lambda}^{\circledS \mr{log}}$
 for 
the classifying morphism 
 of  ${^\S \mfY}_{\dagger}^{\circledS \bigstar}$.
 Then, $\pi^{\circledS \mr{log}}_{\langle \underline{U} \rangle}$  factors through
a morphism 
\begin{align}
\Theta^{\circledS \mr{log}} : [\langle \underline{R} \rangle^{\circledS} \rightrightarrows \langle \underline{U} \rangle^{\circledS}]^\mr{log} \migi {^{\S_1} \overline{\mfM}}_{g,r,\lambda}^{\circledS \mr{log}}.
\end{align}
To complete the proof of Theorem A, it suffices to prove that $\Theta^{\circledS \mr{log}}$ is  an isomorphism.

Consider  the surjective portion.
Let  $S^{\circledS \mr{log}}$ be an object in $\mfS \mfc \mfh_{/S_0}^{\circledS \mr{log}}$ and  $s^{\circledS \mr{log}} : S^{\circledS \mr{log}} \migi {^{\S_1} \overline{\mfM}}_{g,r,\lambda}^{\circledS \mr{log}}$   a morphism of log superstacks,  which  induces a morphism $s_t^\mr{log}  : S_t^\mr{log}  \migi ({^{\S_1} \overline{\mfM}}_{g,r,\lambda})_t^{\mr{log}}$. 
There exists a strict \'{e}tale covering $\pi_{\underline{s}}^\mr{log} : \underline{S}^{' \mr{log}} \migi S_t^\mr{log}$ of $S_t^\mr{log}$ and a morphism ${\underline{s}'}^{\mr{log}} :  \underline{S}^{' \mr{log}} \migi \underline{U}^\mr{log}$ satisfying that
$s_t^{\mr{log}} \circ \pi_{\underline{s}}^\mr{log} \cong   \pi^\mr{log}_{\underline{U}} \circ {\underline{s}'}^{\mr{log}}$.
By Proposition \ref{p0607}, there exists 
a strict super\'{e}tale morphism $\pi_{s}^{\circledS \mr{log}} : S^{' \circledS \mr{log}} \migi S^{\circledS \mr{log}}$  which fits into the following cartesian square diagram
\begin{align}
\xymatrix{
\underline{S}^{'\mr{log}} \ar[r]^{\pi_{\underline{s}}^\mr{log}} \ar[d] \ar@{}[rd]|{\Box} & S_t^\mr{log} \ar[d]^{\tau^{\circledS \mr{log}}_S} \\
S^{' \circledS \mr{log}} \ar[r]_{\pi_{s}^{\circledS \mr{log}}} &  S^{\circledS \mr{log}}.
}
\end{align}
(In particular, the left-hand vertical arrow coincides with $\tau^{\circledS \mr{log}}_{S'}$.)
By Corollary \ref{p04}, the morphism ${\underline{s}'}^{\mr{log}}$ extends   to a morphism ${s'}^{\circledS \mr{log}} : S^{' \circledS \mr{log}}  \migi \langle \underline{U} \rangle^{\circledS \mr{log}}$.
The uniqueness assertion of Corollary \ref{p04} implies that   $\pi_{\langle \underline{U} \rangle}^{\circledS \mr{log}} \circ {s'}^{\circledS \mr{log}} \cong  s^{\circledS \mr{log}} \circ \pi_s^{\circledS \mr{log}}$.
This shows  the subjectivity of  $\Theta^{\circledS \mr{log}}$. 

The injectivity portion follows from  an argument technically similar to the above discussion.
This completes the proof of Theorem A.

\end{document}